\documentclass{amsart}
\usepackage[utf8]{inputenc}

\usepackage{hyperref}
\usepackage{amsmath,amsfonts,amssymb,amscd,amsthm,amsbsy,amsxtra,bbm,bm, epsf,calc,comment,appendix}
\usepackage{color}
\usepackage{datetime}
\usepackage{latexsym}
\usepackage[english]{babel}
\usepackage{enumerate}
\usepackage{graphicx}
\usepackage{epsfig}

\def\XXint#1#2#3{{\setbox0=\hbox{$#1{#2#3}{\int}$ }
\vcenter{\hbox{$#2#3$ }}\kern-.6\wd0}}

\def\R{{ \mathbb{R}}}
\def\Q{\mathbb{Q}}

\def\II{{\rm I\kern-0.5exI}}
\def\III{{\rm I\kern-0.5exI\kern-0.5exI}}

\newcommand{\norm}[1]{\lVert #1 \rVert}

\newcommand{\RR}{\mathbb{R}}

\newcommand{\tp}{\tilde{p}}

\DeclareMathOperator*{\argmin}{argmin}

\DeclareMathOperator{\loc}{\textup{loc}}

\DeclareSymbolFont{bbold}{U}{bbold}{m}{n}
\DeclareSymbolFontAlphabet{\mathbbold}{bbold}

\newcommand{\bn}{\bar{n}}

\newcommand{\vp}{\varphi}
\newcommand{\tr}{tr}
\DeclareMathOperator{\sgn}{sgn}
\DeclareMathOperator{\dist}{dist}

\newcommand{\xs}{x_{\sigma}}

\newcommand{\tilT}{\widetilde{T}}

\newcommand{\spt}{\textup{spt}}

\numberwithin{equation}{section}
\newtheorem{theorem}{Theorem}[section]
\newtheorem{lemma}[theorem]{Lemma}
\newtheorem{prop}[theorem]{Proposition}
\newtheorem{cor}[theorem]{Corollary}

\theoremstyle{remark}
\newtheorem{remark}[theorem]{Remark}

\theoremstyle{definition}
\newtheorem{definition}[theorem]{Definition}
\begin{document}
\title{Free boundary regularity for tumor growth with nutrients and diffusion}

\author{Carson Collins, Matt Jacobs and Inwon Kim}

\begin{abstract}
In this paper, we study a tumor growth model where the growth is driven by nutrient availability and the tumor expands according to Darcy's law with a mechanical pressure resulting from the incompressibility of the cells. Our focus is on the free boundary regularity of the tumor patch that holds beyond topological changes.
%,starting from general initial data. 
A crucial element in our analysis is establishing the regularity of the \textit{hitting time} $T(x)$, namely the first time the tumor patch reaches a given point.  We achieve this by introducing a novel Hamilton-Jacobi-Bellman (HJB) interpretation of the pressure, which is of independent interest.  The HJB structure is obtained by viewing the model as a limit of the Porous Media Equation (PME) and building upon a new variant of the AB estimate. Using the HJB structure, we establish a new Hopf-Lax type formula for the pressure variable.  Combined with barrier arguments, the formula allows us to show that $T$ is $C^{\alpha}$ with $\alpha=\alpha(d)$, which translates into a mild nondegeneracy of the tumor patch evolution.    Building on this and obstacle problem theory, we show that the tumor patch boundary is regular in $\R^d\times (0,\infty)$ except on a set of Hausdorff dimension at most $d-\alpha$.  On the set of regular points, we further show that the tumor patch is locally $C^{1,\alpha}$ in space-time.  This conclusively establishes that instabilities in the boundary evolution do not amplify arbitrarily high frequencies.  %As a consequence, the associated pressure gradient is well-defined, uniformly nonzero, and has continuous directions on the set of regular points.

\end{abstract}

\maketitle

\section{Introduction}

 In this paper, we consider the following tumor growth model:
\begin{equation}\label{first}
\partial_t \rho - \nabla \cdot (\rho \nabla p ) = n\rho, \quad  p(1-\rho)=0, \quad \rho\leq 1,
\end{equation}
where $\rho$ denotes the density of tumor cells, $p$ denotes the pressure, and $n$ is a nutrient variable that evolves according to the diffusion equation
\begin{equation}\label{second}
\partial_t n - \Delta n = -n \rho.
\end{equation}
The form of the pressure-density relation reflects the incompressibility of the tumor cells, namely the pressure variable $p$ acts as the {\it Lagrange multiplier} for the constraint $\rho\leq 1$.   In short, the system (\ref{first}-\ref{second}) describes a cell growth system where the growth rate is mediated by nutrient availability and the tumor region expands according to Darcy's law
with a mechanical pressure driven by the incompressibility of the cells.    Models of the form (\ref{first}-\ref{second}) have been extensively studied by both the mathematical and biological communities with various different assumptions on the growth term and density pressure coupling \cite{BYRNE2003567, Preziosi2008, Ranft20863, maury14, pqv}, to name just a few.  Nonetheless, many mathematical questions remain outstanding, in particular, those regarding the long-time behavior of the tumor boundary region.

Our focus on the specific source term $n\rho $ is due to the fact that the model (\ref{first}-\ref{second}) generates particularly interesting behavior of the tumor patch despite the apparent simplicity of the coupling between the tumor and nutrient.   It is well-known in the biology literature (through numerical and physical experiments) that the tumor patch generated by this model exhibits a fingering instability (c.f. the discussion in \cite{kitsunezaki97}, \cite{MVS02}, \cite{maury14}, \cite{golden22}, \cite{jkt_nutrient}).  In particular, it has been unclear whether this fingering phenomenon occurs at some discrete scale or whether it leads to an immediate or eventual loss of regularity in the tumor boundary. Investigating this behavior will be the main goal of this paper.  %\textcolor{magenta}{In this paper, we conclusively establish that the fingering does not lead to an immediate loss of regularity and hence must be occurring at some discrete scale. }

%\textcolor{blue}{should we not state something about what we establish here in rought words? like " This is what we address in this paper. Namely we will show that the loss of regularity occurs only in small scales." }

Although the tumor system nearly corresponds to that of the classical {\it Hele-Shaw flow}, a mathematically rigorous study of the boundary behavior has remained elusive, due to the difficulties presented by the source term $\rho n$.  In the classical setting, which we will call the {\it injection problem}, the Hele-Shaw flow is given with no source (namely $n=0$) and with a fixed boundary from which the flow is injected at a given rate.   For the injection problem, the global structure of the boundary $\partial \{\rho=1\}$ is well understood by now, mainly through comparison principle type arguments \cite{CJK, CJK2, dong21} or via connections to the obstacle problem \cite{Baiocchi1973, monneau, figalli_serra, figalli_generic}.  

For our problem, the comparison approach is immediately ruled out, as the full system (\ref{first}-\ref{second}) does not have comparison (though note that the individual equations when considered separately do have comparison principles).  As such, we shall proceed via the obstacle problem analysis. However, there is a highly nontrivial roadblock that must be overcome.  Indeed, the source term $\rho n$ necessarily depends on the space-time geometry of the free boundary, while for the injection case, the source is concentrated at a fixed boundary that is safely away from the free boundary. This makes the analysis of the tumor system considerably more difficult, as the influence of the source term cannot be ignored when blowing up the problem at free boundary points (the fundamental technique for the obstacle problem approach).  In particular, to use the obstacle problem toolbox, one must first establish the regularity of the {\it hitting time } $T(x)$,  which records the first time that the tumor patch reaches the point $x$ (ignoring the regularity issues, one can formulate $T(x)$ as $\inf\{t>0: \rho(t,x)=1\}$, see equation (\ref{eq:hitting_time}) for a more careful definition).  This is essentially equivalent to establishing a quantitative non-degeneracy property for the tumor expansion speed, a highly nontrivial task. 

To establish the regularity of $T(x)$ we first derive a novel Hopf-Lax type estimate for the pressure (c.f. Theorem  \ref{thm:main_1}).  To the best of our knowledge, such Hopf-Lax type formulas have not previously appeared in the Hele-Shaw literature, perhaps in part due to the difficulty of controlling the time derivative of $p$.   We get around this by viewing equation (\ref{first}) as the incompressible limit of the Porous Media Equation (PME).   
Given some parameter $\gamma\in (1,\infty)$, the PME analogue of (\ref{first}) is the equation 
\begin{equation}\label{eq:pme_intro}
\partial_t \rho_{\gamma}-\nabla \cdot (\rho_{\gamma}\nabla p_{\gamma})=\rho_{\gamma} n_{\gamma}, \quad p_{\gamma}=\rho_{\gamma}^{\gamma},
\end{equation}
where $n_{\gamma}$ will solve (\ref{second}) with $\rho$ replaced by $\rho_{\gamma}$, and (\ref{first}) can be recovered by sending $\gamma\to\infty$ (see for instance \cite{pqv, perthame_david, jacobs_2021}).   Since the pressure-density coupling $p_{\gamma}=\rho_{\gamma}^{\gamma}$ is explicit for PME, one can rewrite (\ref{eq:pme_intro}) solely in terms of the pressure, namely,
\begin{equation}\label{eq:pme_pressure_intro}
\partial_t p_{\gamma}-|\nabla p_{\gamma}|^2-\gamma p_{\gamma}(\Delta p_{\gamma}+n_{\gamma})=0.
\end{equation}
Interestingly, we ignore the parabolic structure of this equation and instead focus on the Hamilton-Jacobi-Bellman (HJB) structure of the first two terms.  We then build upon the recent improved versions of the Aronson-Ben\'ilan estimate introduced in \cite{jacobs_lagrangian} to show that the positive part of $u_{\gamma}:=-\gamma(\Delta p_{\gamma}+n_{\gamma})$ is uniformly bounded with respect to $\gamma$ in a BMO type space, implying that our limiting $p$ must be a supersolution to the HJB equation
\begin{equation}\label{eq:p_super}
\partial_t p-|\nabla p|^2+ pu_+\geq 0
\end{equation}
where $u:=\lim_{\gamma\to\infty} u_{\gamma}.$
From here we finally obtain the Hopf-Lax formula by adapting the techniques of \cite{cardaliaguet_graber_mfg} for HJB equations with unbounded coefficients.    It is highly intriguing to speculate whether it is possible to obtain (\ref{eq:p_super}) or Hopf-Lax estimates directly from (\ref{first}), however, we will not consider this line of inquiry further in this work.

%Although it has been well-known that (\ref{first}) can be obtained as a limit of PME (see for instance \cite{pqv, perthame_david, jacobs_2021}),  the previous work in this direction has not succeeded in using this limit to extract information about $\partial_t p$, due to the difficulty of controlling the term $\gamma p_{\gamma}(\Delta p_{\gamma}+n)$ as $\gamma\to\infty$.  Here, we build upon the novel Aronson-Ben\'ilan (AB) type estimates introduced in \cite{jacobs_lagrangian}, to show that the positive part of $u_{\gamma}:=-\gamma(\Delta p_{\gamma}+n)$ can be bounded in BMO uniformly in $\gamma$.  As a result, taking the incompressible limit in  (\ref{eq:pme_pressure_intro}), we find that our pressure in (\ref{first}) satisfies
%\begin{equation}\label{eq:p_super}
%\partial_t p-|\nabla p|^2\geq -p u_+
%\end{equation}
%where $u_+=\lim_{\gamma\to\infty} [u_{\gamma}]_+$.    Using the Hamilton-Jacobi structure of (\ref{eq:p_super}) we derive the Hopf-Lax bound for the pressure, where we also adapt some ideas from  \cite{cardaliaguet_graber_mfg} to handle the fact that $u_+\notin L^{\infty}$.

 Once we have the Hopf-Lax formula, we combine this with a powerful barrier-type argument to prove that for any point $x\notin \spt(\rho_0)$ and any sufficiently small radius $r>0$ there exists an explicit time $t_r(x)<T(x)$ such that the tumor patch does not occupy any point in $B_r(x)$.  From here, it will follow that the hitting time is H\"{o}lder continuous with an exponent that depends on the dimension only.  With the H\"{o}lder continuity of $T$ in hand, we can turn to the obstacle problem formulation to address the regularity of the free boundary. Here, the novelty in our analysis lies in establishing the global space-time regularity of the free boundary, with data that is far less regular than the typical injection problems that have previously been considered.

Ultimately, through the obstacle problem analysis, we are able to show that the free boundary is regular except at topological singularities, which are unavoidable for general initial data. This conclusively demonstrates that the observed instabilities for the system (\ref{first}-\ref{second}) do not amplify arbitrarily high frequencies and must occur at some fixed scale.   In particular, we show that the tumor patch boundary is regular in $\R^d\times (0,\infty)$ except on a relatively closed set of Hausdorff dimension at most $d-\alpha$ for some $\alpha\in (0,1)$ depending only on the dimension.  On the set of regular points, we further show that the tumor patch is $C^{1,\alpha}$ in space, locally uniformly in time. It then follows that the associated pressure gradient at regular boundary points is well-defined and uniformly positive in space-time. Moreover, the direction of the pressure gradient on the set of regular points is continuous in space-time.

In the remainder of the introduction, we give a more complete explanation of the obstacle formulation of our problem and the connection to the hitting time.  We then %give a brief overview of our estimates on the hitting time using the Hopf-Lax formula and barrier supersolutions.  Finally, we 
summarize our main results and give a roadmap for the rest of the paper.  

\subsection{The obstacle problem and the hitting time}

To better understand the aforementioned difficulties and the importance of the hitting time, let us describe some properties of the tumor patch and formally introduce the obstacle problem associated to (\ref{first}-\ref{second}).  
Since our main interest is the regularity properties of the tumor patch, throughout the paper, we will assume that 
\begin{equation}\label{initial}
\rho(x,0)\hbox{ is a characteristic function and } n(x,0) \hbox{ is uniformly positive.}
\end{equation}  
 %\textcolor{blue}{We are making this assumption so we don't have to deal with new patches nucleating from regions where $0<\rho<1$. The obstacle problem can't say anything about the regularity when these sets nucleate (for the same reason that obstacle problem can't say anything about initial data.) }
Under these assumptions, $\rho$ will remain a characteristic function for all times and $t\mapsto \rho(x,t)$ will be nondecreasing for a.e. $x\in \RR^d$.

Transitioning to the obstacle problem formulation, if we integrate the pressure variable in time,
\begin{equation}\label{w_def}
w(x,t):=\int_0^t p(x,s) ds,
\end{equation}
the new variable $w$, the so-called \textit{Baoicchi transform}, will satisfy an obstacle problem \cite{Baiocchi1973}.  Since the density is nondecreasing in time, the relation $(1-\rho)p=0$ implies that $(1-\rho)w=0$.  Using the patch property for the density, this coupling can be upgraded to the even stronger relation that the sets $\{w>0\}$ and $\{\rho=1\}$ coincide spacetime almost everywhere (c.f. Lemma \ref{lem:rho_w_positive_sets_agree}).  This key relation can then be combined with the time integral of (\ref{first}) to see that $w$ solves the elliptic obstacle problem 
\begin{equation}\label{eq:w_obst_eqn}
        \Delta w = (1 - \rho_0 -\eta)\chi_{\{w>0\}}, 
    \end{equation}  
    where $\eta(x,t) :=\int_{0}^t\rho(x,s) n(x,s)\, ds$ (c.f. Lemma \ref{obstacle}).  

The main challenge in analyzing (\ref{eq:w_obst_eqn}) is the presence of the term $\eta$, which is absent in the obstacle formulation of the classical injection case (due to local regularity results, $\rho_0$ does not affect the free boundary regularity at positive times  away from the support of $\rho_0$).  Since $\rho$ is a characteristic function, it is not clear whether $\eta$ has any nice regularity.  This is crucial, as obstacle problem regularity theory breaks down without Dini continuity of the coefficients (see \cite{blank}).  Hence, one must hope  
that the time integral induces some smoothing effect.  At the very least, this can only happen if the tumor boundary is strictly expanding. Indeed, if any part of the free boundary stagnates in time, then $\eta$ will become discontinuous across that portion of the boundary.  Note that such stagnation would correspond to a jump in the values of the hitting time function $T$ introduced earlier. Hence, the smoothness of $\eta$ and $T$ are highly intertwined.  In fact, it will turn out that we can express $\eta$ solely in terms of the hitting time $T$ and $n$.  

To see the connection between $\eta$ and $T$, we need to first give a proper definition of the hitting time.
Recall that the hitting time $T(x)$ records the first time that the tumor patch arrives at a point $x$.  We will formally define it using $w$, the most regular variable at our disposal.    Given a point $x\in \RR^d$ we set
\begin{equation}\label{eq:hitting_time}
   T(x):= \inf \{t>0: w(x,t)>0\}.
\end{equation}
Since the positivity set of $w$ coincides almost everywhere with the tumor patch, we have $\rho(x,t)=\sgn_+(t-T(x))$ almost everywhere. Hence, $\eta$ can be rewritten in terms of $T$ and $n$ as
\begin{equation}\label{eq:eta_alternate}
    \eta(x,t)=\sgn_+(t-T(x))\int_{T(x)}^t n(x,s)\, ds.
\end{equation}
From the above formula, we now see that the spatial regularity of $\eta$ is more or less equivalent to the regularity of $T$ and $n$.

Note that generically $T$ is at best Lipschitz continuous, as it is easy to cook up a scenario where two different parts of the tumor patch collide with different velocities.  In addition, topological changes of the tumor boundary can cause the pressure to suddenly jump with highly nonlocal effects.  For instance, the merger of two portions of the boundary can cause far away parts of the boundary to instantaneously start moving faster.   Since $n$ is much better than Lipschitz continuous, it is $T$ that will determine the regularity of $\eta$.  While we are inclined to believe that the Lipschitz continuity of $T$ is true, our methods are only able to show that $T$ is H\"{o}lder continuous with a dimensionally dependent exponent. Nonetheless, the H\"{o}lder continuity is sufficient for us to deduce free boundary regularity using the obstacle problem approach. However, let us note that we are forced to work in a much lower regularity regime than what is typically considered in the obstacle problem literature, requiring us to develop new arguments.
%Hence, even in the most optimistic possible case, the resulting obstacle problem will have much worse regularity than what is typically considered in the literature \textcolor{blue}{(CITATION NEEDED)}.  While we believe that the Lipschitz continuity of the hitting time is probably true, our arguments are only able to show that 

%Without expending too much effort, one can show that $T$ is continuous using soft arguments (c.f. Proposition \ref{prop:hitting_time_continuous}). However, obtaining the quantitative continuity that we need is much harder.  It is here where our new arguments come in. 
    
\begin{comment}
    To see this, note that by integrating (\ref{first}) in time, it is clear that $w$ solves the equation
\begin{equation}\label{eq:first_time_integral}
    \rho-\Delta w=\rho_0+\int_0^t \rho n.
\end{equation}
Since $\rho$ is nondecreasing, we can also characterize the tumor patch as the set where $w>0$, allowing us to
\end{comment}

  %\textcolor{blue}{I am a bit confused here, isn't the characteristic function $\chi_{w>0}$ in (1.5) already require $t>T(x)$in the $\eta$ expression, eliminating the need for the sign term?} \textcolor{red}{Yes, the sign term is only so that $\eta$ is defined consistently with the ``consumed nutrient" function $\int n \rho\,dt$ which arises in the derivation of the equation; $\eta$ can be freely redefined on the set where $w = 0$. That said, I think we discussed stating \eqref{eqn:simple_w_eqn} instead here and introducing the obstacle formulation later? }

\subsection{Main results}
We are now ready to present the main results of our paper.  All of our results will use the following mild assumptions on the initial data.
\begin{enumerate}
   \item[(A1)] $\rho(\cdot,0)\in L^1(\RR^d)\cap \textup{BV}(\RR^d)$ and $\rho(x,0)\in \{0,1\}$ for almost every $x\in \RR^d$.
\item[(A2)] $n(\cdot,0)\in W^{1,\infty}(\RR^d)$ and there exists $c>0$ such that $n(x,0)\geq c$ for all $x\in \RR^d$. 
\end{enumerate}

%Roughly, the first half of the paper will be building towards the goal of showing that $T$ is H\"{o}lder continuous, obtaining the HJB structure of the pressure and the Hopf-Lax formula as important independent results along the way.  In the second half of the paper, we will use the H\"{o}lder continuity of $T$ and the obstacle problem structure to study the regularity of the free boundary.
The main results of the first half of the paper are the HJB structure and Hopf-Lax formula for the pressure, along with the H\"{o}lder continuity of the hitting time.
\begin{theorem}\label{thm:main_1}
The following holds for the unique weak solution $p$ to the system \eqref{first}-\eqref{second}.

\begin{itemize}
\item[(a)] [Cor \ref{cor:p_hjb}] $p$ solves, in the sense of weak solutions,
$$
\partial_t p - |\nabla p|^2 +pu_+ \geq 0,
$$
 where for any $\tau > 0$ there exists $b=b(\tau, d) > 0$ such that $bu_+ e^{bu_+}  \in L^1([0, \tau ] ; \R^d)$.\\
 
\item[(b)] [Prop.~\ref{prop:hj_estimate}] 
 Given points $(x_1, t_1)$, $(x_0, t_0)$ with $t_0<t_1$ and any decreasing function $\lambda\in L^1([0, t_1-t_0])$, there exist constants $C=C(t_1,d)$ and $b=b(t_1, d)$ such that
 \begin{equation*}
p(x_0, t_0)\leq e^{\Lambda(t_1-t_0)}\Big( p(x_1, t_1)+\frac{|x_1-x_0|^2}{4\int_{0}^{t_1-t_0} e^{\Lambda(s)}\, ds} +C(t_1-t_0)^{7/10}e^{-\lambda(t_1-t_0)}\Big)
\end{equation*}
where $\Lambda(t):=\frac{5}{4b}\int_0^t \lambda(s)\, ds+\frac{t}{b}\log(1+\frac{C}{t})$. \\ 
\item[(c)][Theorem \ref{thm:hitting_time_holder}] $T$ is locally H\"{o}lder continuous on the set $\{x\in \RR^d: 0<T(x)<\infty\}$ with an exponent that depends only on the dimension.
\end{itemize}
\end{theorem}

Let us note that Theorem \ref{thm:main_1} parts (a) and (b) represent a significant improvement to our understanding of the Hele-Shaw equation. In particular, any control on the time derivative of the pressure has been previously missing in the literature. Furthermore, the delicate control that we obtain from the Hopf-Lax formula in part (b) is completely new and unexpected.   

As we mentioned earlier, we establish the HJB structure by first going through the PME (\ref{eq:pme_gamma}).  For the classic PME without a source term, bounds on the negative part of $\Delta p_{\gamma}$ are known through the celebrated Aronson-Benilan estimate \cite{ab}.  In the presence of a source term, AB-type bounds on quantities taking a similar form to $u_{\gamma}=-\gamma(\Delta p_{\gamma}+n)$ have been studied in the literature \cite{pqv, gpsg, perthame_david, Bevilacqua2022, jacobs_lagrangian}, however except for \cite{jacobs_lagrangian}, these bounds do not scale well with respect to $\gamma$.  We adapt the arguments from \cite{jacobs_lagrangian} to show that $[u_{\gamma}]_+$ can be bounded uniformly with respect to $\gamma$ in a BMO-type space.  Once we have the uniform control on $u_{\gamma,+}$ we can pass to the limit in (\ref{eq:pme_pressure}) to obtain the result (a). A direct derivation of (a) from the Hele-Shaw flow or the meaning of the singular limit $u=\lim_{\gamma\to\infty}u_{\gamma}$ in terms of the Hele-Shaw flow remains open.  

To obtain (b), we cannot take the usual approach to proving Hopf-Lax type formulas (i.e. differentiating $p$ along paths) due to the potential unboundedness of $u_+$.  To overcome this, we adapt the approach developed in \cite{cardaliaguet_graber_mfg}, which handles unbounded coefficients by averaging over paths
indexed by the unit ball.  Our calculation is somewhat different however, as we can exploit the specific structure of $pu_+$ to decompose $pu_+\leq \lambda p+p(u-\lambda)_+$ for some scalar $\lambda\geq 0$.  By choosing $\lambda$ appropriately we can force $p(u-\lambda)_+$ to be small while using a Gronwall argument to handle $\lambda p$.  This allows us to obtain a much more favorable error term in our Hopf-Lax formula compared to \cite{cardaliaguet_graber_mfg}.

Although Theorem 1.1 (a) and (b) are stated for our particular system (\ref{first}-\ref{second}), equivalent results can be proved for more general tumor growth models where the growth term $\rho n$ is replaced by $\rho G$ for some general growth rate $G$.  In particular, our arguments only need $G\geq c(\tau)>0$ along with some control on $[\partial_t G]_-$.

As a consequence of the Hopf-Lax formula, we obtain the $C^{\alpha}$ regularity of $T$ where $\alpha=\alpha(d)$.  We do this by combining the formula with a novel barrier type argument.  Given a point $(x, T(x))$ on the free boundary and some time $t_0\in (0,T(x))$, 
we use the Hopf-Lax formula and the values of $p$ at time $T(x)$ to construct an explicit supersolution $\psi$ that dominates $p$ on $(t_0, T(x))\times\RR^d$.  The key is that the Hopf-Lax formula allows us to choose the values in such a way that $\psi$ is zero in a neighborhood of $x$ up until the hitting time $T(x)$.  Since we are able to explicitly calculate and invert $t(r)=\inf\{t\in (t_0, T(x)): \sup_{y\in B_r(x)} \psi(y, t)>0\}$ we obtain an upper bound on $\sup_{y\in B_r(x)} T(x)-T(y)$, which implies the H\"{o}lder continuity of $T$.

Some remarks on the previous literature for hitting times are in order.
Quantitative regularity of the hitting time for PME has been obtained in \cite{CF80} for the classical PME and in \cite{paul21} for the PME with a source term and drift.  Nevertheless, both of these results obtain estimates that blow up as $\gamma$ tends to infinity, due to the lack of a uniform AB estimate on $u_{\gamma}$.  As a result, their approaches are not suitable for our problem.  Let us also note that these papers used a rather different approach that did not involve the Hopf-Lax approach that we use here.  Estimates on the hitting time for a simpler version of (\ref{first}) where $n$ is replaced by a decreasing function of $p$, were obtained in \cite{pqm} for dimensions $d\leq 3$.  Their proof strongly relies on the specific structure of their growth term, which allows them to relate the H\"{o}lder continuity of $T$ to that of the pressure through a clever trick.  Again, this approach is not applicable to our problem.  Although we also focus on a specific source term, our method is much more general and can be applied to other instances of the Hele-Shaw or Porous Media equation.

\medskip

The remaining analysis in the paper is devoted to the study of the obstacle problem \eqref{eq:w_obst_eqn}, based on the $C^{\alpha}$ regularity of $T$. 
We build on the low-regularity obstacle problem analysis of Blank \cite{blank} to establish the space-time regularity of the tumor patch. A crucial fact we use is that the solution of the obstacle problem with $C^{\alpha}$ data has a unique blow-up limit at each point, allowing us to decompose the boundary into a regular part and a singular part (the regular points have blow up limits that look like half-planes). A direct application of this dichotomy yields that the boundary has locally finite $H^{d-1}$ measure for each time, as mentioned for instance in \cite{pqm}.
However, this standard description lacks the geometric information of the free boundary over time.  Indeed, the main novelty of our obstacle problem analysis is that we are able to stitch together information from each time $t$ to obtain regularity of the full space time boundary $\Gamma:=\{(x, T(x): x\in \RR^d\}\subset \RR^d\times (0,\infty)$.
In particular, we show that $\Gamma$ is regular in space-time outside of a set of at most Hausdorff dimension $d-\alpha$, and its outward normal is H\"{o}lder continuous in space-time.  While space-time analysis of the singular set has been carried out before for the injection problem (\cite{monneau}, \cite{figalli_serra}, \cite{figalli_generic}), these results have utilized smoothness (at least $C^4$) of the fixed boundary data in an essential way. A more general time-varying source term was considered in \cite{serfaty_serra}, but only for a short range of time that ensures that no topological singularity occurs during the evolution.  

Our results are summarized in the following Theorem. 
\begin{theorem}
Let $\Gamma$ denote the space-time boundary set of the tumor region i.e. $\Gamma=\{(x,T(x)): x\in \RR^d\}$
\begin{itemize}

\item[(a)][Prop. \ref{prop:regular_open}, Prop. \ref{prop:singular_manifold}] The set $\{ 0 < T(x) < \infty \}\subset \R^d$ decomposes as $R\cup \Sigma$, where the set $R$ of regular points is open in $\R^d$ and the set $\Sigma$ of singular points is locally contained in a $C^1$ manifold of dimension $d-1$. 

\item[(b)][Prop. \ref{prop:regular_final_boundary_regularity}]
 At any $x\in R$, the free boundary near $(x,T(x))\in\Gamma$ can be locally represented as a graph $\{x_n = f(x', t)\}$ where $f$ is $C^{1,1}$ in $x'$ and Lipschitz in time. 

\item[(c)] 
[Prop. \ref{prop:regular_linear_growth}, Cor.\ref{cor:regular_points_t_lipschitz} ] $p(\cdot,T(x))$ has linear growth at $x\in R$, with locally uniform growth rates. In particular $T$ is Lipschitz in $R$.

%namely for $y$ sufficiently close to $x$ and for sufficiently small $r$ we have $\sup_{B_r(y)} p(\cdot, T(y)) \sim r$. 

\item[(d)][Prop. \ref{prop:unit_normal_holder}] The map $\nu: R\to \mathcal{S}^d$, where $\nu(x)$ denotes the spatial outward normal of $\Gamma$ at $(x, T(x))$,  is H\"{o}lder continuous. In particular $\nabla p(x, T(x))$ is well-defined for $x\in R$ and has continuous direction.

\end{itemize}
\end{theorem}
Note that, while Theorem 1.2 (d) yields the continuity of the direction of $\nabla p$ on $\Gamma$, we cannot expect the same for $|\nabla p|$: this can be easily seen from examples where a topological change occurs far away from the given free boundary point.
In terms of the quadratic blow up limit of $w$, Theorem 1.2 (a) and (d) yield its continuity at free boundary points, along $\Sigma$ and along $R$. As a consequence, $D^2 w(x, T(x))$ exists on $\Sigma$ and exists in a one-sided sense on $R$, and is continuous on each set (though not necessarily on their union).

\medskip

Our last theorem discusses the Hausdorff measure of the free boundary in space-time coordinate. Let us introduce the notation
\begin{equation}\label{domain}
\Omega_t := \{w(\cdot,t)>0\}, \quad \Gamma_t := \partial\Omega_t.
\end{equation}

\begin{theorem}[Corollary \ref{cor:hausdorff_measure}]
\begin{itemize}  
    \item[(a)] The free boundary $\partial\{w>0\}$ has Hausdorff dimension d in $(x,t)$-coordinates. \\
    \item[(b)]  $Graph(R) = \{(x,t): x\in \Gamma_t, x\in R\}$ is relatively open with locally finite $H^d$ measure.\\
    \item[(c)]
     $Graph(\Sigma)= \{(x,t): x\in \Gamma_t, x\in \Sigma\}$
      has
locally finite $H^{d-\alpha}$ measure.
\end{itemize}
\end{theorem}

\medskip

 Let us mention that we expect $T$ to be Lipschitz for all points,  not just in $R$.  For instance, in the classical setting with constant Dirichlet fixed boundary data, the Lipschitz continuity of $T$ was shown by \cite{monneau} by a simple comparison principle. The remaining challenge in our setting lies in the analysis of the singular points. This is an intriguing question as the blow-up profile of the tumor patch at these points suggests that the evolution at these points should be non-degenerate in general. In fact, one might even expect the gradient of $T$ to vanish at these points.   Nonetheless, accurately capturing the hitting time behavior near singular points appears to be out of reach for the moment. It would also be interesting to improve upon our estimate of the singular set, using a generic notion of initial data. While it seems plausible, new ideas seem to be necessary to obtain such a result.  

\medskip

The rest of the paper is organized as follows. In Section 2, we review the basic properties of the system (\ref{first}-\ref{second}) and the connection to the obstacle problem. In Section 3, we establish the HJB structure of the pressure along with the Hopf-Lax formula.  In Section 4, we construct a barrier supersolution using the Hopf-Lax formula which allows us  to establish the H\"{o}lder regularity of the hitting time map $T(x)$. Section 5 builds on the regularity of $T$ and the existing obstacle theory to investigate the global regularity of the free boundary $\Gamma:= \partial\{w(x,t)>0\}$.

  \section{Basic properties of the system}

Here we recall the notion and show basic properties of solutions for \eqref{first}-\eqref{second}.  We first introduce the weak notion of our solutions, parallel to those introduced in \cite{jkt_nutrient} for a similar model. 

  \begin{definition}\label{def:weak_sln}
      A triple $(\rho, p, n)$ is a weak solution to \eqref{first}-\eqref{second} for initial data $\rho_0\in L^1(\R^d)\cap BV(\R^d)$ and $n_0\in L^\infty(\R^d)\cap BV(\R^d)$ if for any $\tau > 0$,
      \begin{enumerate}[(i)]
          \item $p(1-\rho) = 0$ in $\mathcal{D}'(\R^d\times [0, \tau])$
          \item For any $\psi\in H^1(\R^d\times [0,\tau])$ vanishing at time $\tau$,
          \begin{equation}\label{eq:first_dist_eqn}
              \int_0^{\tau} \int_{\R^d} \nabla \psi\cdot \nabla p - \rho \partial_t \psi\,dx\,dt = \int_{\R^d} \psi(x, 0)\rho_0(x)\,dx + \int_0^{\tau} \int_{\R^d} \psi n \rho \,dx\,dt.
          \end{equation}
          \item 
          \[ \partial_t n - \Delta n = -\rho n\hbox{ in }\mathcal{D}'(\R^d\times[0, \tau]),\quad n(x, 0) = n_0(x) \]
          \item We have $\rho\in C([0, \tau]; L^1(\R^d))\cap L^\infty([0, \tau]; BV(\R^d))$, $p\in L^2([0, \tau]; H^1(\R^d))$, and $n\in L^\infty(Q_{\tau})\cap L^\infty([0, \tau]; BV(\R^d))$.

      \end{enumerate}
      Here, $BV(\R^d)$ is the space of (not necessarily integrable) functions with finite total variation.
  \end{definition}

We also record a few useful properties of a weak solution.

\begin{lemma}\label{lem:rho_regularity}\label{lem:n_regularity}\label{lem:p_dist_eqn}
For any $\tau > 0$, $\varepsilon\in (0,1)$, we have
\begin{enumerate}[(i)]
    \item $\rho\in C^{0,1}_t L^1_x([0,\tau]; \R^d)$.
    \item The support of $\rho$ in $\R^d\times [0, \tau]$ is compact, and $0\leq \rho \leq 1$.
    \item For a.e. $x\in \R^d$, $\rho(x, \cdot)$ is increasing in time.
    \item $n\in L^{\infty}([0,\tau]\times\RR^d),\; \partial_t n \in \textup{BMO}([0,\tau]\times\RR^d), \;D^2 n \in \textup{BMO}([0,\tau]\times\RR^d).$
    \item We have 
        \begin{equation}\label{eq:p_dist_eqn}
            p(\Delta p + n) = 0\hbox{ in }\mathcal{D}'(\R^d\times [0, \tau])
        \end{equation}
        Also, $p\in L^\infty(\R^d\times [0, \tau])$.
\end{enumerate}
\end{lemma}
\begin{proof}
    Statements (i) and (ii) are proved in \cite{jkt_nutrient}, Theorem 2.2 and Proposition 3.6.
    
(iii) follows from Lemma 3.11 of \cite{jkt_nutrient}, which provides comparison for the equation
\[ \partial_t \rho^i - \nabla\cdot(\rho^i \nabla p^i) = f^i \]
Namely, if $\rho^0(x,0) \leq \rho^1(x,0)$, $f^0 \leq f^1$, and $p^i (1 - \rho^i) = 0$, then $\rho^0(x,t)\leq \rho^1(x,t)$. For our system, since $n\rho \geq 0$, comparison to the system with $p^0 = f^0 = 0$ implies that $\rho(x, t_0) \leq \rho(x, t_1)$ for any $t_0 \leq t_1$ and a.e. $x$. The measure zero set where this fails depends on $t_0, t_1$, so we conclude by applying this with a countable basis of intervals.

Item (iv) follows from parabolic estimates for the heat equation with $L^{\infty}$ coefficients (see e.g. \cite{ogawa22}).

The distributional equation $p(\Delta p + n) = 0$ is also proved in \cite{jkt_nutrient}, Theorem 2.2. The $L^\infty$ bound for $p$ follows from the compact support of $\rho$, and the boundedness of $n$; one can take a sufficiently large paraboloid supersolution to $\Delta p = -n$ to get an upper bound.
\end{proof}

From the definition of the weak solution, we can derive that $w$ satisfies an elliptic equation at each time.
  
\begin{lemma}\label{lem:simple_w_eqn}\label{lem:w_regularity}
    We have $w(1 - \rho) = 0$ a.e. in space-time. For each $t > 0$, $w$ solves
    \begin{equation}\label{eqn:simple_w_eqn}
       \Delta w(x, t) = \rho(x,t) - \rho(x,0) - \int_0^t n(x,s)\rho(x,s)\,ds\hbox{ in }\mathcal{D}'(\R^d) 
    \end{equation}
    In particular, for any $\tau > 0$ and $\varepsilon \in (0,1)$, we have $w\in C^{0,1}(\R^d\times [0,\tau])\times L^\infty_t C^{1,1-\varepsilon}([0, \tau], \R^d)$.
\end{lemma}
\begin{proof}
    For the first statement, we simply note that the monotonicity of $\rho$ in time from Lemma \ref{lem:rho_regularity} implies
    \[ 0 \leq w(x,t)(1 - \rho(x, t)) = \int_0^t p(x,s)(1 - \rho(x,t))\,ds \leq \int_0^t p(x,s)(1-\rho(x,s))\,ds \equiv 0 \]

    For the second statement, we first note that both $\rho$ and the function $\int_0^t n\rho\,ds$ are continuous in time into any $L^p$ with $p < \infty$; this follows from the weak solution definition, since we have $\rho\in C_t L^1$ with values in $[0, 1]$ and compact support for bounded time intervals, while $n$ is spacetime continuous from Lemma \ref{lem:n_regularity}. Then to derive a distributional equation for $w$, we consider Definition \ref{def:weak_sln}(ii) with $\psi$ of the form $\psi(x, t) = \varphi(x)\chi(t)$, where $\varphi\in C^\infty_c(\R^d)$ and $\chi\in C^\infty_c([0, \tau))$ with $\chi\equiv 1$ near 0. Using that $w = \int p\,dt$, we obtain
    \[ \int_{\R^d} \chi \nabla \varphi\cdot \nabla w(\cdot, \tau)\,dx - \int_{\R^d}\int_0^{\tau} \partial_t \chi \varphi \rho \,dx\,dt = \int_{\R^d}\varphi(x)\rho_0(x)\,dx + \int_{\R^d} \varphi \int_0^{\tau} \chi n \rho \,dt\,dx \]
    If we take for $\chi$ a sequence of cutoffs valued in $[0,1]$ and converging pointwise to the indicator of $[0, T)$, then we can apply the aforementioned time continuity of $\rho$ and $\int_0^t n\rho\,ds$ to obtain the limiting equation
    \[ \int_{\R^d} \nabla \varphi \cdot \nabla w(\cdot, \tau) + \varphi \rho(\cdot, \tau)\, dx = \int_{\R^d} \varphi (\rho_0 + \eta(\cdot, \tau))\,dx \]
    from which we conclude \eqref{eqn:simple_w_eqn}.

    Since $w(x, t) = \int_0^t p(x, s)\,ds$, the upper bound for $p$ implies that $w$ is Lipschitz in time uniformly in space. This will improve to Lipschitz in spacetime once we have $w\in L^\infty_t C^{1,1-}_x$.
        
    Since $n$ is bounded on $\R^d\times [0,\tau]$, $\eta$ is bounded on $\R^d\times [0,\tau]$. Then since $\Delta w = (1 - \eta)\chi_{\{w > 0\}}$ is uniformly bounded in $L^\infty$, and up to time $\tau$, $w(\cdot, t)$ is compactly supported in $\overline{\Omega_\tau}$, it follows that for any $p\in (1, \infty)$, $\Delta w(\cdot, t)\in L^p(\R^d)$. Calderon-Zygmund estimates then give $w(\cdot, t)\in W^{2,p}(\R^d)$, and thus $w(\cdot, t)\in C^{1,1-\varepsilon}(\R^d)$ for any $\varepsilon > 0$, uniformly in $t\in [0, \tau]$.
\end{proof}

From now on we make the assumptions (A1) and (A2) on our initial data $\rho_0, n_0$.

\begin{lemma}\label{lem:nutrient_lower_bound}
    Let $\bar{n}(t):=\inf_{x\in \RR^d} n(t,x)$  For any $t>0$, $\bn(t)\geq e^{-t} \bn(0).$
\end{lemma}
\begin{proof}
   Suppose that $\tilde{n}$ satisfies the equation $\partial_t \tilde{n}-\alpha \Delta \tilde{n}=-\tilde{n}$ with constant initial data $\tilde{n}(0,x)=\bn(0).$  The comparison principle for the heat equation implies that $\tilde{n}\leq n$ almost everywhere.  If we define $\tilde{N}=e^{t}\tilde{n}$, then $\tilde{N}$ satisfies $\partial_t \tilde{N}-\alpha\Delta\tilde{N}=0$ with initial data $\tilde{N}(0,x)=\bn(0)$.  Hence, $\tilde{N}(t,x)=\tilde{N}(0,x)=\bar{n}(0)$ and thus it follows that $\tilde{n}(t,x)=e^{-t}\bn(0)$, which implies the result.
\end{proof}

The following characterization of the pressure variable replaces the formal description of $p$ solving the elliptic problem $-\Delta p = n$ in $\{\rho=1\}$ with zero Dirichlet data, to avoid ambiguity rising from potentially irregular boundary of $\{\rho=1\}$. The argument is similar to ones that have previously appeared in the literature \cite{pqv, santambrogio_crowd_motion, guillen_kim_mellet, jacobs_2021}, here we include a proof since our setting is slightly different.

\begin{lemma}\label{lem:variational}
    For almost every time $t$, the pressure is a solution to the variational problem
    \[
    p(\cdot,t)=\argmin_{\vp(\cdot)(1-\rho(\cdot,t))=0, \; \vp\geq 0} \, \int_{\RR^d} \frac{1}{2}|\nabla \vp(x)|^2-\vp(x) n(x,t)\, dx.
    \]
\end{lemma}
\begin{proof}
Given $\epsilon>0$ define $p_{\epsilon}(x,t):=\frac{1}{\epsilon}\int_{t-\epsilon}^t p(x,s)\, ds$ where we set $p(x,s)=p(x,0)$ if $s<0$ and $p^{\epsilon}(x,t):=\frac{1}{\epsilon}\int_{t}^{t+\epsilon} p(x,s)\, ds$.  Fix a time $t_0$ such that $p_{\epsilon}(t_0,\cdot), p^{\epsilon}(t_0,\cdot)$ converge to $p(t_0,\cdot)$ in $H^1(\RR^d).$
Choose some nonnegative function $\vp\in H^1(\RR^d)$ such that $\vp(x)(1-\rho(x,t_0))=0$ for almost every $x\in \RR^d$ (note that space integrals of $\rho(\cdot,t_0)$ against functions in $H^1(\RR^d)$ are well defined at any time $t_0$ since $\partial_t \rho\in L^2([0,T];H^{-1}(\RR^d))$, which itself is a consequence of the continuity equation and $p\in L^2([0,T];H^1(\RR^d))$ ).

 Integrating equation (\ref{first}) from time $t_0-\epsilon$ to $t_0$, dividing by $\epsilon$, and integrating against $\vp$ we see that 
 \[
 \int_{\RR^d} \vp(x)\frac{\rho(x,t_0)-\rho(x,t_0-\epsilon)}{\epsilon}+\nabla \vp(x)\cdot \nabla p_{\epsilon}(x,t_0)\, dx=\int_{\RR^d} \vp(x) \frac{1}{\epsilon}\int_{t_0-\epsilon}^{t_0} \rho(x,s) n(x,s)\, ds\, dx
 \]
 The condition $\vp(x)(1-\rho(x,t_0))=0$ implies that $ \vp(x)\frac{\rho(x, t_0)-\rho(x, t_0-\epsilon)}{\epsilon}= \vp(x)\frac{1-\rho(x, t_0-\epsilon)}{\epsilon}$. Combined with the constraint $\rho\leq 1$, we can conclude that
  \[
 \int_{\RR^d}\nabla \vp(x)\cdot \nabla p_{\epsilon}(x,t_0)\, dx\leq \int_{\RR^d} \vp(x) \frac{1}{\epsilon}\int_{t_0-\epsilon}^{t_0} \rho(x,s)n(x,s)\, ds\, dx.
 \]
 Applying the same logic to the time integral over the interval $[t_0,t_0+\epsilon]$, we find that
  \[
 \int_{\RR^d}\nabla \vp(x)\cdot \nabla p^{\epsilon}(x,t_0)\, dx\geq \int_{\RR^d} \vp(x) \frac{1}{\epsilon}\int^{t_0+\epsilon}_{t_0}\rho(x,s) n(x,s)\, ds\, dx.
 \]
 Sending $\epsilon\to 0$ we can conclude that 
 \begin{equation}
 \int_{\RR^d} \nabla \vp(x)\cdot \nabla p(x,t_0)=\int_{\RR^d} \vp(x)\rho(x, t_0)n(t_0,x)=\int_{\RR^d} \vp(x)n(x, t_0)
 \end{equation}
 where the final equality follows from the fact that $\vp(x)(1-\rho(x,t_0))=0$ almost everywhere. The above equation is the Euler-Lagrange equation for the variational problem,  thus, combined with the strong convexity of the variational problem, we see that $p$ solves the variational problem at every time $t_0$ where $p_{\epsilon}(\cdot, t_0), p^{\epsilon}(\cdot, t_0)$ converge to $p(\cdot, t_0)$ in $H^1(\RR^d).$  Since this must hold for almost every $t_0\in [0,T]$ we are done. 
\end{proof}

A straightforward consequence of the previous Lemma is the following Lemma which gives a crude comparison between the pressure values at different times.  We will obtain a much sharper comparison property in Section 3 when we establish the Hopf-Lax type formula for the pressure. 
\begin{lemma}\label{lem:containment}
    Fix some time $\tau>0$. Given almost any times $s,t\in [0,\tau]$ such that $s<t$, there exists a constant $C(\tau)$ such that
    \[
    p(s,x)\leq C(\tau) p(t,x).
    \]
    In particular, this implies
    \[
    \{x\in \RR^d: w(s,x)>0\}\subset \{x\in \RR^d: p(s,x)>0\}\subset  \{x\in \RR^d: w(t,x)>0\}.
    \]
\end{lemma}
\begin{proof}
  By Lemma \ref{lem:nutrient_lower_bound}, there exists a constant $C(\tau)>0$ such that  $n(x,s)< C(\tau) n(x,t)$ for all $x\in \RR^d$.  Since $\rho$ is increasing with respect to time and $\rho\leq 1$, we know that $\vp(x)(1-\rho(x,t))=0$ for any nonnegative function $\vp(x)$ such that $\vp(x)(1-\rho(x,s))=0$.  Let us choose $\vp(x)=(p(s,x)-C(\tau)p(t,x))_+$.  It then follows from Lemma \ref{lem:variational} that 
  \[
   \int_{\RR^d} \nabla (p(x,s)-C(\tau)p(x,t))_+ \cdot \nabla p(x,t)\, dx=\int_{\RR^d} (p(x,s)-C(\tau)p(x,t))_+n(x,t)\, dx, 
   \]
   and
   \[
   \int_{\RR^d} \nabla (p(x,s)-C(\tau)p(x,t))_+ \cdot \nabla p(x,s)\, dx=\int_{\RR^d} (p(x,s)-C(\tau)p(x,t))_+n(x,s)\, dx.
  \]
  Hence, 
  \[
   \int_{\RR^d} \nabla (p(x,s)-C(\tau)p(x,t))_+ \cdot \nabla (p(x,s)-C(\tau)p(x,t))\, dx=\int_{\RR^d} (p(x,s)-C(\tau)p(x,t))_+(n(x,s)-C(\tau)n(x,t))\, dx.
  \]
  The left-hand side of the above equation is nonnegative while the right-hand side of the equation is nonpositive.  This is only possible if $(p(x,s)-C(\tau)p(x,t))_+=0$ almost everywhere. 
\end{proof}

\begin{lemma}\label{lem:rho_w_positive_sets_agree}
    Up to a set of measure zero, for any $t > 0$ we have $\{x\in \RR^d: \rho(x,t)=1\}=\{x\in \RR^d: w(x,t)>0\}$.
\end{lemma}
\begin{proof}
    This is nearly Lemma 4.6 of \cite{jkt_nutrient}, except that in the diffusion case we lack an explicit formula for the nutrient. Nevertheless, we proceed along the same lines. 

    From Lemma \ref{lem:simple_w_eqn}, we have $w(1-\rho) = 0$, and thus
    \[ \{x\in \RR^d: w(x,t)>0\}\subset \{x\in \RR^d: \rho(x,t)=1\} \]
    Thus, we must show that the set $A_t := \{ x : \rho(x, t) = 1, w(x,t) = 0 \}$ has measure zero.
    
    For this, we observe that $\Delta w$ vanishes a.e. where $w$ vanishes, and thus \eqref{eqn:simple_w_eqn} implies that 
    \[ \rho(x, 0) + \int_0^t n(x,s)\rho(x, s)\,ds = \rho(x,t) = 1 \hbox{ a.e. on }A_t\]
    From the pressure equation \eqref{eq:p_dist_eqn}, any interior point of $\{ \rho(x, 0) = 1 \}$ has positive pressure at every positive time, and our assumptions on the initial data provide that the boundary of this set has zero measure.
    
    Thus, we need only consider the case where $\int_0^t n\rho\,ds = 1$. Since the nutrient is uniformly positive due to Lemma \ref{lem:nutrient_lower_bound}, this occurs for at most one time for a given $x$. On the other hand, since the nutrient is uniformly bounded, this function is continuous in time, and $x$ must be in $A_t$ for an open set of times before $\int_0^t n\rho\,ds = 1$ is satisfied. It follows directly that $A_t$ (and, in fact, $\bigcup_t A_t$) is null.
    
\end{proof}

\begin{lemma}\label{lem:alt_hitting_time}
    \[ T(x) := \inf \{ t \geq 0 : w(x, t) > 0 \} = \inf \{ t \geq 0 : \rho(x, t) = 1\} \hbox{ for a.e. } x\in \R^d \]
\end{lemma}
\begin{proof}
    For conciseness, write $\tilT(x) = \inf \{ t\geq 0 : \rho(x, t) = 1 \}$.

    Suppose that for some $x$, we have $\tilT(x) < T(x)$. Then we have $w(x, t) = 0$ for all $t < T(x)$. Then, since Lemma \ref{lem:rho_regularity} gives that $\rho$ is monotone in time except for a null set of $x$ which we ignore, we have $x\in \{ y : w(y, t) = 0, \rho(y, t) = 1 \}$ for all $t\in (\tilT(x), T(x))$. Then any such $x$ is contained in $\bigcup_{t\in \Q\cap (0, \infty)} \{ y : w(y,t) = 0, \rho(y, t) = 1 \}$, which is null by Lemma \ref{lem:rho_w_positive_sets_agree}.

    The other direction is similar. Suppose instead that for some $x$, we have $T(x) < \tilT(x)$. Then $w(x, t) > 0$ for all $t > T(x)$, while monotonicity implies that for a.e. $x$ we have $\rho(x, t) = 0$ for all $t < \tilT(x)$. Then any such $x$ is contained in $\bigcup_{t\in \Q\cap (0, \infty)} \{ y : w(y, t) > 0, \rho(y, t) = 0 \}$, which is null by Lemma \ref{lem:rho_w_positive_sets_agree}.
\end{proof}

\begin{lemma}\label{obstacle}
    $w(\cdot,t)$ solves the obstacle problem
  \eqref{eq:w_obst_eqn}.
\end{lemma}
\begin{proof}
    This follows from applying Lemma \ref{lem:alt_hitting_time} to Lemma \ref{lem:simple_w_eqn}. We have
    \[\rho(x, t) = \chi_{\{ T < t\}}(x) = \chi_{ \{ w(\cdot, t) > 0\} }(x)\] 
    for a.e. $x$ and all $s\neq T(x)$. Thus, to obtain \eqref{eq:w_obst_eqn} we modify \eqref{eqn:simple_w_eqn} by replacing the occurrence of $\rho$ in $\int_0^t n\rho\,ds$ with $\chi_{\{ T < t\}}$ and the other occurrence of $\rho$ with $\chi_{\{w(\cdot, t) > 0\}}$.
\end{proof}

Recalling the notation of \eqref{domain}, we additionally define
\begin{equation}\label{domain:2}
\Omega_\infty := \{0\leq T(x) <\infty\}=\bigcup_{t > 0} \Omega_t, \quad \mathcal{O} :=\{0<T(x)<\infty\} = \Omega_{\infty}\setminus \overline{\Omega_0}.
\end{equation}

We now prove that the hitting time $T$ is continuous. While this justifies the characterization of the level sets of $T$ as the free boundary $\Omega_t$, it also is an important first step that initiates the regularity analysis of $T$ in section 4. The main idea will be to show that a discontinuity must result in a point $x_0$ and times $t_0 < t_1$ such that $x_0$ is in $\partial \Omega_t$ for $t_0 \leq t \leq t_1$. Then $w(\cdot, t_1) - w(\cdot, t_0)$ is a positive superharmonic function on $\Omega_{t_0}$, so one would like to apply the Hopf lemma to draw a contradiction between $w(x_0, t_0) = w(x_0, t_1) = 0$ and $\nabla w(x_0, t_0) = \nabla w(x_0, t_1) = 0$. Unfortunately, $\Omega_{t_0}$ does not a priori have the regularity needed to apply the Hopf lemma, so we must first use obstacle problem techniques to shift to a setting where we do have such regularity. The key tool in doing so will be the quadratic blowup of $w$ at free boundary points:

\begin{lemma}[\cite{blank} Corollary 2.5]\label{lem:quad_blowup_cpt}
    Let $u$ be a nonnegative solution to $\Delta u = f\chi_{\{u > 0\}}$, for some $f$ which is strictly positive and bounded near the free boundary $\partial \{ u > 0\}$. Then if $x_0$ is a free boundary point, then the quadratic blowup sequence $r^{-2}w(r(x - x_0) + x_0, t)$ is compact in $C^{1,\alpha}(B_1(x_0))$ as $r\to 0^+$. Moreover, if $f$ is continuous at $x_0$, then the subsequential limits solve $\Delta v = f(x_0)\chi_{\{ v > 0\}}$.
\end{lemma}

The subsequential limit enjoys better geometry, due to the following property of global solutions to the constant-source obstacle problem:

\begin{lemma}[\cite{caffarelli98} Corollary 7]\label{lem:global_obst_sln_cvx}
    A nonnegative solution to $\Delta u = \chi_{\{ u > 0\}}$ on $\R^d$ is convex.
\end{lemma}

\begin{prop}\label{prop:hitting_time_continuous}
 \begin{itemize}
\item[(a)]    $T$ is continuous.  
\item[(b)] $x\in \partial\Omega_t$ if and only if $x\in \mathcal{O}$ and $t = T(x)$, for all $x\in\R^d$ and $t>0$.
\end{itemize}
\end{prop}
\begin{proof}

First, we verify that the $\overline{\Omega_t}$ are continuous from above, in the sense that for any $t$ we have:
\begin{equation}\label{eq:set_continuity}
    \overline{\Omega_t} = \bigcap_{\varepsilon > 0} \overline{\Omega_{t + \varepsilon}}
\end{equation}
The forward inclusion is trivial by the monotonicity of $w$. For the reverse inclusion, we suppose for contradiction that there exists $x \in \bigcap_{\varepsilon > 0} \overline{\Omega_{t + \varepsilon}} \setminus \overline{\Omega_t}$. Let $r$ be sufficiently small that $B_r(x)\subset \{ w(\cdot, t) = 0 \}$. From \eqref{eqn:simple_w_eqn}, $\Delta w(\cdot, t + \varepsilon) \geq 1 - \varepsilon \|n_0\|_{\infty}$ on $B_r(x)\cap \{ w(\cdot, t + \varepsilon) > 0 \}$, and by assumption we have $x\in \overline{\Omega_{t + \varepsilon}}$ for all $\varepsilon > 0$. It follows by quadratic nondegeneracy for the obstacle problem (Lemma \ref{lem:obst_nondegen}) that if $\varepsilon < \frac{\|n_0\|_{\infty}}{2}$, then
\[ \sup_{B_r(x)} w(\cdot, t + \varepsilon) \geq Cr^2 \]
uniformly in $\varepsilon$. Since $w(\cdot, t) \equiv 0$ on $B_r(x)$, we get a contradiction with the Lipschitz continuity of $w$ in time by shrinking $\varepsilon$.

\medskip

Now, we introduce 
\begin{equation}\label{eq:t_0_defn}
 T_0(x) := \inf \{ t > 0 : x\in \overline{\Omega_t} \}   
\end{equation}
It is immediate that $T_0(x) \leq T(x) = \inf \{ t > 0 : x\in \Omega_t \}$. We claim that $x\in \partial \Omega_t$ if and only if $t\in [T_0(x), T(x)]$. It is clear that for $t < T_0(x)$, $x\notin \overline{\Omega_t}$, and for $t > T(x)$, $x\in \Omega_t$. Since $x\in \Omega_{T(x) + \varepsilon}$ for every $\varepsilon > 0$, \eqref{eq:set_continuity} gives $x\in \overline{\Omega_{T(x)}}$. On the other hand, by continuity of $w$ and minimality of $T$, we have $x\notin \Omega_{T(x)}$, so $x\in \partial \Omega_{T(x)}$. Monotonicity implies that $x$ is a boundary point for all $t \leq T(x)$ for which $x\in \overline{\Omega_t}$. We have $x\in \overline{\Omega_{T_0(x)}}$ by using \eqref{eq:set_continuity} with the definition of $T_0$, so we get the claim.

For purely topological reasons related to how each is defined, $T_0$ is lower semicontinuous and $T$ is upper semicontinuous. To check lower semicontinuity of $T_0$, let $(x_n)$ be a sequence converging to $x$ with $\liminf T_0(x_n) := t$. Then for any $\varepsilon > 0$, the $x_n$ are eventually in $\overline{\Omega_{t + \varepsilon}}$, and thus $x\in \overline{\Omega_{t + \varepsilon}}$. It follows that $T_0(x) \leq t$. To check upper semicontinuity of $T$, we note that $x\in \Omega_{T(x) + \varepsilon}$ for all $\varepsilon > 0$. Any sequence $x_n$ converging to $x$ eventually has $T(x_n) \leq T(x) + \varepsilon$ since $\Omega_{T(x) + \varepsilon}$ is open, so we conclude that $\limsup T(x_n) \leq T(x)$.

Therefore, both parts of the proposition will follow if we can show that $T_0 \equiv T$. For this, we will need the following useful property:
\begin{equation}\label{eq:t_t0_continuity}
    \lim_{\Omega_{T_0(x)}\ni x_n \to x} T(x_n) = T_0(x)
\end{equation}
To see this, we note that if $(x_n)$ is such a sequence, then $T(x_n) \leq T_0(x)$ for each $n$. On the other hand, by minimality of $T_0$, $x\notin \overline{\Omega_{T_0(x) - \varepsilon}}$, and so we eventually have $T(x_n) \geq T_0(x) - \varepsilon$ for any $\varepsilon > 0$.

\begin{comment}
First, we observe that from the definition of $T$, we get a property which is stronger than lower semicontinuity, for $x\in \mathcal{O}$: 

\begin{equation}\label{continuity_0}
T(x) = \liminf_{\mathcal{O}\ni x_n\to x} T(x_n) = \lim_{\overline{\Omega_{T(x)}}\ni x_n\to x} T(x_n)
\end{equation}

Indeed, for any sequence $(x_n)$ in $\R^d$ with $x_n\to x$, we must have $T(x) \leq \liminf T(x_n)$. Otherwise, for any time $t > \liminf T(x_n)$ there would be a sequence of points in $\Omega_t$ converging to $x$, giving $x\in \overline{\Omega_t}$, and $T(x)\leq t$. On the other hand, we have $x\in \overline{\Omega_{T(x)}}$ by definition of $T$ \textcolor{blue}{shouldn't we explain a bit more, we are using that these sets do not jump over time due to continuity of $\rho$ in $L^1$} \textcolor{red}{I now realize this was originally written for a slightly different $T$ with the opposite semicontinuity: $T(x) = \inf \{ t : x\in \overline{\Omega_t}\}$... I will make the fixes.} and any sequence $(x_n)$ in $\overline{\Omega_{T(x)}}$ with $x_n\to x$ satisfies $T(x_n) \leq T(x)$ for each $n$. It follows that $\limsup T(x_n) \leq T(x)$ for sequences in $\overline{\Omega_{T(x)}}$.

\end{comment}

\medskip
Finally, we proceed to the proof that $T_0 = T$. Suppose that $x_0\in \partial\Omega_t$ with $T_0(x_0) < T(x_0)$. We will use the obstacle problem theory to compare blowups of $w$ at $(x_0, T_0(x_0))$ and at $(x_0,t)$ with $t>T_0(x_0)$ to derive a contradiction. First let us ensure that the blow-up profiles are well-defined. Due to (\ref{eq:w_obst_eqn}) it follows that
$$
0\leq \eta(x,t) \leq (t - T(x))_+\|n_0\|_\infty
$$
We also get a continuity estimate. Assuming $T(y) \leq T(x)$, we have either $T(y) \leq t \leq T(x)$, giving
\[ |\eta(x, t) - \eta(y, t)| = \eta(x, t) \leq \|n_0\|_\infty |T(x) - T(y)| \]
or else $t \leq T(y) \leq T(x)$, giving
\begin{align*}
    |\eta(x,t) - \eta(y,t)| &= \left|\int_{T(x)}^t n(x,s) - n(y,s)\,ds - \int_{T(y)}^{T(x)} n(y,s)\,ds\right|
    \\&\leq |t - T(x)|\sup_{T(x)\leq s\leq t} |n(x,s) - n(y,s)| + \|n_0\|_\infty |T(x) - T(y)|
\end{align*}
Thus, using the nutrient regularity from Lemma \ref{lem:n_regularity} and the result of \eqref{eq:t_t0_continuity}, we conclude that $\eta(\cdot, T_0(x_0))$ restricted to $\Omega_{T_0(x_0)}$ is continuous at $x_0$, and in a sufficiently small neighborhood of $x_0$, we can ensure that it is less than $\frac{1}{2}$. Then Lemma \ref{lem:quad_blowup_cpt} gives that the family of rescalings $x\mapsto r^{-2}w(r(x - x_0) + x_0, T_0(x_0))$ are compact as $r\to 0$ in $C^{1,\alpha}_{loc}$, and their subsequential limits are nonzero global solutions of 
\begin{equation}\label{obstacle_00}
\Delta u = \chi_{\{u > 0\}}.
\end{equation}

\medskip

Now, choose $\tau\in (T_0(x), T(x))$, sufficiently small such that we still have $\eta(x, \tau) < \frac{1}{2}$ in some neighborhood of $x_0$. By taking a further subsequence, the discussion above yields a sequence $r_n\to 0$ such that 
$$
r_n^{-2}w(r_n(x - x_0) + x_0, T(x_0))\to u \hbox{ and } r_n^{-2}w(r_n(x - x_0) + x_0, \tau)\to v,
$$
 for some $u, v$ in $C^{1,\alpha}_{loc}(\R^d)$. Unlike with $u$, $\eta(\cdot, \tau)$ restricted to $\Omega_\tau$ is not known to be continuous at $x_0$, since we do not yet know that $T$ is continuous. In particular, we do not know that $v$ solves a constant Laplacian obstacle problem. However, we do get $u(x_0) = v(x_0) = 0$ and $\nabla u(x_0) = \nabla v(x_0) = 0$ from the convergence, and we also have that $v\geq u$ since $w(\cdot, t) \geq w(\cdot, T_0(x_0))$.

We will apply the Hopf lemma to $v-u$ in the domain $U: = \{ u > 0\}$. First observe that from the definition of $u$, we have $r_n(x - x_0) + x_0 \in \Omega_{T(x_0)}$ if $x\in U$ and if $n$ is sufficiently large depending on $x$. We have checked above that $\eta(\cdot, \tau)$ restricted to $\Omega_{T_0(x_0)}$ is continuous at $x_0$. Thus, for any $x\in U$, we have $\eta(r_n(x - x_0) + x_0, \tau) \to \eta(x_0, \tau)$, and so we conclude that 
$$\Delta v = 1 - \eta(x_0, \tau) \hbox{ in } U.$$ Comparing this equation to \eqref{obstacle_00}, it follows that $v - u$  satisfies 
\[ \Delta (v - u) = - \eta(x_0, \tau) \leq - (\tau - T_0(x_0))\inf_{t\in [T_0(x_0), \tau]}n(x_0, t) < 0 \]
with the last inequality following from the nutrient lower bound in Lemma \ref{lem:nutrient_lower_bound} and the assumption that the initial nutrient is bounded away from 0. This implies that $v - u$ is strictly superharmonic inside $U$, and so our previous observation that $v - u\geq 0$ by the monotonicity in time of $w$ improves to $v - u > 0$ inside $U$. Lastly let us observe that, from Lemma \ref{lem:global_obst_sln_cvx}, the complement of $U$ is convex and so $U$ satisfies the interior ball condition at $x_0$. Putting together the above information, the Hopf lemma applied at $x_0$ implies that $\nabla v(x_0) - \nabla u(x_0) \neq 0$, which is a contradiction. It follows that $T_0 = T$, so we finish.
\end{proof}
\begin{remark}
    An important consequence of Proposition \ref{prop:hitting_time_continuous} is that the spacetime interface is exactly the graph of $T$ on $\mathcal{O}$. In other words,
    \begin{equation}\label{eq:graph}
        \{ (x, t) : t\in (0. \infty), x\in \partial \Omega_t \} = \mathrm{Graph}_T(\mathcal{O}) := \{ (x, T(x)) : x\in \mathcal{O} \}
    \end{equation}
    This also means that the interface is a $d$-dimensional topological manifold, and the regularity of its parametrization in $d$ spatial variables is exactly that of $T$. We will use the notation $\mathrm{Graph}_T$ with subsets of $\mathcal{O}$, which may be understood in this light as projections of the spacetime interface into $\R^d$.
\end{remark}

Finally, in light of the regularity of $w$ and $T$, we note a natural way to standardize $\rho$ on measure zero sets.  A corresponding standardization of the pressure will need to wait until the next Section, due to the need to preserve certain delicate structures.

\begin{lemma}
    $\partial \Omega_t$ has zero measure in $\R^d$ for all $t > 0$. The weak solution $(\rho, p, n)$ can be taken such that $\rho$ is upper semicontinuous in space and time, with $\{ \rho(\cdot, t) = 1 \} = \overline{\Omega_t}$ for each $t$. In particular, the support of $p$ is then contained in $\{ \rho = 1 \}$.
\end{lemma}
\begin{proof}
By Proposition \ref{prop:hitting_time_continuous}, for any $t 
 >0$ and $\varepsilon\in (0, t)$, we have $\partial \Omega_t\subset \Omega_{t + \varepsilon}\setminus \Omega_{t - \varepsilon}$. By Lemma \ref{lem:rho_w_positive_sets_agree}, up to measure zero sets we can replace the right-hand-side with $\{ \rho(\cdot, t + \varepsilon) = 1 \}\setminus \{ \rho(\cdot, t - \varepsilon) = 1 \}$, and by time continuity of $\rho$ in $L^1$, the measure of this set goes to 0 with $\varepsilon$. Thus, $\partial \Omega_t$ has zero measure.
\\
Then we claim that $(\rho, p, n)$ with $\rho$ redefined as $\chi_{\overline{\{ w > 0\}}}$ and $p$ redefined to vanish outside $\overline{\{ w > 0\}}$ remains a weak solution as defined in Definition \ref{def:weak_sln}. Indeed, since this changes $\rho, p$ by measure zero sets for each time, equation \eqref{eq:first_dist_eqn} is unaffected. On the other hand, we have $p (1 - \rho) = 0$ by construction. To check that $\rho$ is spacetime upper semicontinuous, we only need to show that the set $\{ \rho = 1 \} := \{ (x, t) : t\geq 0, x\in \overline{\Omega_t}\}$ is closed. Let $(x_n, t_n)$ be a sequence in this set converging to $(x, t)$. Since the $t_n$ converge, they are bounded, and so the sequence is contained in an $\overline{\Omega_\tau}$ for $\tau$ sufficiently large. Then $T(x) < \infty$, so we have $T(x_n)\to T(x)$ by continuity, and since $T(x_n) \leq t_n$ for each $n$ by Proposition \ref{prop:hitting_time_continuous}, we have $T(x) \leq t$. It follows that $x\in \overline{\Omega_t}$, so we conclude.
\end{proof}

%\textcolor{blue}{In the next section, we will employ an alternative, integral-oriented way to define an upper semi-continuous version of $p$. This definition carries the advantage that the estimates we obtain for $p_{\gamma}$ is preserved.}

\section{AB estimates and the Hopf-Lax bound}\label{sec:hjb}

In this section, we will show that there exists a nonnegative function $u_+$ such that the pressure is a super solution to the following Hamilton-Jacobi equation
\begin{equation}\label{eq:hjb}
\partial_t p-|\nabla p|^2\geq pu_+.
\end{equation}
We will then use (\ref{eq:hjb}) to obtain a Hopf-Lax type formula for the pressure. 
 In particular, given a fixed time $t_0$, this will allow us to give lower bounds for the pressure at times $t>t_0$ and upper bounds for the pressure at times $t<t_0$ in terms of $p(t_0,\cdot)$.  This will give us a very precise way of constructing pressure  super solutions that lead to powerful barrier-type arguments and eventually H\"{o}lder regularity of the hitting times (c.f. Section \ref{sec:barrier}). 

Let us emphasize that to the best of our knowledge, the Hopf-Lax type bounds we obtain have not previously appeared in the literature for Hele-Shaw type equations and they require some highly nontrivial efforts to obtain. 
First, to establish (\ref{eq:hjb}), we go through the Porous Media Equation (PME) and use the fact that our solution $(\rho, p)$ can be obtained as the incompressible limit of solutions $(\rho_{\gamma}, p_{\gamma}, n_{\gamma})$ of the PME-nutrient system
\begin{equation}\label{eq:pme_gamma}
\partial_t \rho_{\gamma}-\nabla \cdot (\rho_{\gamma}\nabla p_{\gamma})=\rho_{\gamma}n_{\gamma}, \quad p_{\gamma}=\rho_{\gamma}^{\gamma},
\end{equation}
\begin{equation}\label{eq:nutrient_gamma}
\partial_t n_{\gamma} -\Delta n_{\gamma}=-\rho_{\gamma} n_{\gamma}
\end{equation}
as the scalar parameter $\gamma$ is sent to infinity, and where the nutrient variable from our original system is held fixed.  The advantage of the PME system is that it is possible to use the relation $p_{\gamma}=\rho_{\gamma}^{\gamma}$ to rewrite (\ref{eq:pme_gamma}) solely in terms of the pressure variable $p_{\gamma}$, which yields the equation
\begin{equation}\label{eq:pme_pressure}
    \partial_t p_{\gamma}-|\nabla p_{\gamma}|^2-\gamma p_{\gamma}(\Delta p_{\gamma}+n)=0.
\end{equation}
The main difficulty in obtaining (\ref{eq:hjb}) is to show that as $\gamma\to\infty$, $u_{\gamma}=-\gamma(\Delta p_{\gamma}+n)$ converges to a meaningful limit object $u$, whose positive part can be controlled. For the classic PME without a source term, bounds on the negative part of $\Delta p_{\gamma}$ are known through the celebrated Aronson-Benilan estimate \cite{ab}.  In the presence of a source term, AB-type bounds on quantities taking a similar form to $\gamma(\Delta p_{\gamma}+n)$ have been studied in the literature \cite{pqv, gpsg, perthame_david, jacobs_lagrangian}, however except for \cite{jacobs_lagrangian}, these bounds do not scale well with respect to $\gamma$.  We adapt the arguments from \cite{jacobs_lagrangian} to show that $[u_{\gamma}]_+$ can be bounded uniformly with respect to $\gamma$ in BMO-type spaces.  Note that we are unable to get $L^{\infty}$ bounds on $u_{\gamma}$ essentially because two key quantities in the estimate $\partial_t n$ and $\nabla n\cdot \nabla p$ are not in general bounded in $L^{\infty}$.  It would be interesting to see whether equation (\ref{eq:hjb}) could be obtained directly from the original system without going through PME, but we leave this question to a future work.

Once we have obtained equation (\ref{eq:hjb}), there is still significant work required to obtain a Hopf-Lax type control for $p$.  Here the difficulty is that $u_+$ is not bounded in $L^{\infty}$.  Noting that $p$ should satisfy
\[
\partial_t p-|\nabla p|^2+pu_+\geq 0,
\]
the derivative of $p$ along an arbitrary path $x(t)$ gives
\[
\frac{d}{dt} p(x(t),t)=\partial_t p(x(t), t)+x'(t)\cdot \nabla p(x(t),t)\geq
\]
\[
|\nabla p(x(t),t)|^2-p(x(t),t)u_+(x(t),t)+x'(t)\cdot \nabla p(x(t),t)\geq -p(x(t),t)u_+(x(t),t)-\frac{1}{4}|x'(t)|^2.
\]
Unfortunately, without $L^{\infty}$ control on $u_+$ it is not clear that time integrals of the final quantity will be well-defined.  This prevents the usual approach to proving Hopf-Lax type bounds.

To overcome this, we adapt the approach developed in \cite{cardaliaguet_graber_mfg}, which handles unbounded coefficients by instead considering an average over paths
indexed by the unit ball.  Our calculation is somewhat different however, as we can exploit the specific structure of $pu_+$ to write $pu_+=\lambda p+p(u-\lambda)_+$ for some scalar $\lambda\geq 0$.  By choosing $\lambda$ appropriately we can force $p(u-\lambda)_+$ to be small while using a Gronwall argument to handle $\lambda p$.  This allows us to obtain a much more favorable error term in our Hopf-Lax formula compared to \cite{cardaliaguet_graber_mfg} (c.f. Proposition \ref{prop:hj_estimate}).

We begin with the aforementioned result (well-known) that says we can approximate our system (\ref{first}-\ref{second}) with a sequence of smooth solutions to PME.
\begin{prop}[see e.g. \cite{pqv, gpsg, jacobs_lagrangian}]
There exists a sequence of smooth solutions $(\rho_{\gamma}, p_{\gamma}, n_{\gamma})$  to the PME-nutrient system (\ref{eq:pme_gamma}-\ref{eq:nutrient_gamma}) with initial data $(\rho_{0,\gamma}, n_{0,\gamma})$ such that for any $\tau>0$ we have that $\rho_{\gamma}$ converges strongly in $L^1([0,\tau]\times\RR^d)$, $p_{\gamma}, n_{\gamma}$ converge strongly in $L^2([0,\tau];H^1(\RR^d))$ to the unique solution $(\rho, p, n)$ to the system (\ref{first}-\ref{second}) with initial data $(\rho_0, n_0)$ as $\gamma\to\infty$.  Furthermore, one may choose $\rho_{0,\gamma}$ such that $u_{\gamma,+}(\cdot,0)$ is bounded in $L^{\infty}(\RR^d)$ uniformly in $\gamma$.
\end{prop}

Next, we record the following simple result for solutions to the heat equation with $L^{\infty}$ source.

\begin{lemma}\label{lem:dt_n_d2_n}
There exists some $b_0>0$ such that $n$ satisfies the bound
    \[
    \exp(b_0\frac{|\partial_t n|+|\Delta n|}{n})-1\in L^1([0,\tau]\times\RR^d)
    \]
\end{lemma}
\begin{proof}
Thanks to Lemma \ref{lem:n_regularity} we know that $\partial_t n$ and $D^2 n$ are bounded in BMO.
    Since we also have
    \[
    \frac{1}{2}\norm{\partial_t n}_{L^2([0,\tau]\times\RR^d)}^2+\norm{\nabla n}_{L^2(\{\tau\}\times\RR^d)}^2\leq \norm{\nabla n}_{L^2(\{0\}\times\RR^d)}^2+\frac{1}{2}\norm{n}_{L^{\infty}([0,\tau]\times\RR^d)}\norm{\rho}_{L^{2}([0,\tau]\times\RR^d)},
    \]
    the BMO bound implies the existence of a constant $c>0$ such that $\exp(c|\partial_t n|)-1, \exp(c|\Delta n|)-1 \in L^1([0,T]\times\RR^d)$.  Following the logic of Lemma \ref{lem:nutrient_lower_bound}, it follows that $n$ is uniformly bounded from below on any time interval.  Hence, there must be an appropriate choice of $b_0$ where the result holds. 
\end{proof}

\begin{prop}\label{prop:ab}
If $(p_{\gamma}, n_{\gamma})$ is a smooth solution to the system (\ref{eq:nutrient_gamma}-\ref{eq:pme_pressure})
for some $\gamma\in (1,\infty)$,  then for any $\tau>0$  there exists $b>0$ that only depends on $\tau$ such that $(b[u_{\gamma}]_+-1)\exp(b[u_{\gamma}]_+)+1$ is uniformly bounded in  $L^1([0,\tau]\times\RR^d)$ with respect to $\gamma$, where $u_{\gamma} := -\gamma(\Delta p_{\gamma} + n_{\gamma})$.
\end{prop}

\begin{proof}

If we differentiate $\frac{1}{\gamma} u_{\gamma}$ with respect to time, we get
\[
\partial_t \frac{1}{\gamma} u_{\gamma}=-\partial_t n_{\gamma}-\Delta \partial_t p_{\gamma}=
-\partial_t n_{\gamma}-\Delta(|\nabla p|^2-p_{\gamma}u_{\gamma})
\]
Expanding the Laplacian, we see that
\[
\partial_t \frac{1}{\gamma} u_{\gamma}=-\partial_t n_{\gamma}-2|D^2p_{\gamma}|^2-2\nabla \Delta p_{\gamma}\cdot \nabla p_{\gamma}+2\nabla p_{\gamma}\cdot \nabla u_{\gamma}+p_{\gamma}\Delta u_{\gamma}+u_{\gamma}\Delta p_{\gamma}.
\]
Noting that $-\Delta p_{\gamma}=n_{\gamma}+\frac{1}{\gamma} u_{\gamma}$ we can rewrite the previous line as
\[
\partial_t \frac{1}{\gamma} u_{\gamma}=2\nabla n_{\gamma}\cdot \nabla p_{\gamma}-\partial_t n_{\gamma}-2|D^2p_{\gamma}|^2+2(1+\frac{1}{\gamma})\nabla p_{\gamma}\cdot \nabla u_{\gamma}+p_{\gamma}\Delta u_{\gamma}-n_{\gamma} u_{\gamma}-\frac{1}{\gamma} u_{\gamma}^2.
\]
Hence, eliminating $|D^2p|^2$, we can conclude that
\begin{equation}\label{eq:u_ineq}
n_{\gamma}u_{\gamma} +  \frac{1}{\gamma} (\partial_t u_{\gamma}+u_{\gamma}^2)\leq 2\nabla n_{\gamma}\cdot \nabla p_{\gamma}-\partial_t n_{\gamma}+2(1+\frac{1}{\gamma})\nabla p_{\gamma}\cdot \nabla u_{\gamma}+p_{\gamma}\Delta u_{\gamma}.
\end{equation}

Now let $f:\RR\to\RR$ be a $C^2$ convex function such that $f'\geq 0$ everywhere and $f=0$ on $(-\infty,0]$.     If we integrate (\ref{eq:u_ineq}) against $f'(u_{\gamma})$ on $[0,\tau]\times \RR^d$ we find that
\begin{multline}\label{eq:f_u_1}
\int_{\RR^d\times \{\tau\}} \frac{1}{\gamma} f(u_{\gamma})+\int_{\RR^d\times [0,\tau]} n_{\gamma}u_{\gamma}f'(u_{\gamma})+\frac{1}{\gamma}u_{\gamma}^2f'(u_{\gamma})\leq \\
\int_{\RR^d\times \{0\}} \frac{1}{\gamma}f(u_{\gamma})+\int_{\RR^d\times [0,\tau]} f'(u_{\gamma})\big(2\nabla n_{\gamma}\cdot \nabla p_{\gamma}-\partial_t n_{\gamma}\big)+2(1+\frac{1}{\gamma})\nabla p_{\gamma}\cdot \nabla \big(f(u_{\gamma})\big)+p_{\gamma}f'(u_{\gamma})\Delta u_{\gamma} .
\end{multline}
Noting that $\Delta \big(f(u_{\gamma})\big)=f'(u_{\gamma})\Delta u_{\gamma} +f''(u_{\gamma})|\nabla u_{\gamma}|^2$, we can integrate by parts in (\ref{eq:f_u_1}) to obtain
\begin{multline}\label{eq:f_u_2}
\int_{\RR^d\times \{\tau\}} \frac{1}{\gamma} f(u_{\gamma})+\int_{\RR^d\times [0,\tau]} n_{\gamma}u_{\gamma}f'(u_{\gamma})+\frac{1}{\gamma}u_{\gamma}^2f'(u_{\gamma})+p_{\gamma}f''(u_{\gamma})|\nabla u_{\gamma}|^2\leq \\
\int_{\RR^d\times \{0\}} \frac{1}{\gamma}f(u_{\gamma})+\int_{\RR^d\times [0,\tau]} f'(u_{\gamma})\big(2\nabla n_{\gamma}\cdot \nabla p_{\gamma}-\partial_t n_{\gamma}\big)-(1+\frac{2}{\gamma})f(u_{\gamma})\Delta p_{\gamma} 
\end{multline}
We then integrate by parts in $\nabla n_{\gamma}\cdot \nabla p_{\gamma}$ to get 
\begin{multline}\label{eq:f_u_3}
\int_{\RR^d\times \{\tau\}} \frac{1}{\gamma} f(u_{\gamma})+\int_{\RR^d\times [0,\tau]} n_{\gamma}u_{\gamma}f'(u_{\gamma})+\frac{1}{\gamma}u_{\gamma}^2f'(u_{\gamma})+p_{\gamma}f''(u_{\gamma})|\nabla u_{\gamma}|^2\leq \\
\int_{\RR^d\times \{0\}} \frac{1}{\gamma}f(u_{\gamma})-\int_{\RR^d\times [0,\tau]} f'(u_{\gamma})\big(2p_{\gamma}\Delta n_{\gamma}+\partial_t n_{\gamma}\big)+p_{\gamma} f''(u_{\gamma})\nabla u_{\gamma}\cdot \nabla n_{\gamma}+(1+\frac{2}{\gamma})f(u_{\gamma})\Delta p_{\gamma} .
\end{multline}
Once again using $-\Delta p_{\gamma}=n_{\gamma}+\frac{1}{\gamma} u_{\gamma}$ and using the quadratic Young's inequality on $\nabla u_{\gamma}\cdot \nabla n_{\gamma}$ we find that
\begin{multline}\label{eq:f_u_4}
\int_{\RR^d\times \{\tau\}} \frac{1}{\gamma} f(u_{\gamma})+\int_{\RR^d\times [0,\tau]} n_{\gamma}u_{\gamma}f'(u_{\gamma})+\frac{1}{\gamma}u_{\gamma}^2f'(u_{\gamma})+\frac{1}{2}p_{\gamma}f''(u_{\gamma})|\nabla u_{\gamma}|^2\leq \\
\int_{\RR^d\times \{0\}} \frac{1}{\gamma}f(u_{\gamma})+\int_{\RR^d\times [0,\tau]} (1+\frac{2}{\gamma})f(u_{\gamma})(n_{\gamma}+\frac{1}{\gamma}u_{\gamma}) -f'(u_{\gamma})\big(2p_{\gamma}\Delta n_{\gamma}+\partial_t n_{\gamma}\big)+\frac{1}{2}p_{\gamma} f''(u_{\gamma})| \nabla n_{\gamma}|^2.
\end{multline}
Next, to help compare the left and right hand sides, we divide and multiply by multiples of $n_{\gamma}$ to get
\begin{multline}\label{eq:f_u_5}
\int_{\RR^d\times \{\tau\}} \frac{1}{\gamma} f(u_{\gamma})+\int_{\RR^d\times [0,\tau]} n_{\gamma}u_{\gamma}f'(u_{\gamma})+\frac{1}{\gamma}u_{\gamma}^2f'(u_{\gamma})+\frac{1}{2}p_{\gamma}f''(u_{\gamma})|\nabla u_{\gamma}|^2\leq \\
\int_{\RR^d\times \{0\}} \frac{1}{\gamma}f(u_{\gamma})+\int_{\RR^d\times [0,\tau]} (1+\frac{2}{\gamma})f(u_{\gamma})(n_{\gamma}+\frac{1}{\gamma}u_{\gamma}) -\frac{n_{\gamma}}{2}f'(u_{\gamma})\big(\frac{4p_{\gamma}\Delta n_{\gamma}+2\partial_t n_{\gamma}}{n_{\gamma}}\big)+n_{\gamma} f''(u_{\gamma})\frac{p_{\gamma}| \nabla n_{\gamma}|^2}{2n_{\gamma}}.
\end{multline}
Using the identity $uf'(u)-f(u)=f^*(f'(u))$ and applying Young's inequality to $-\frac{n_{\gamma}}{2}f'(u_{\gamma})\big(\frac{2p_{\gamma}\Delta n_{\gamma}+\partial_t n_{\gamma}}{n_{\gamma}}\big)$, we get
\begin{multline}\label{eq:f_u_6}
\int_{\RR^d\times \{\tau\}} \frac{1}{\gamma} f(u_{\gamma})+\int_{\RR^d\times [0,\tau]} (n_{\gamma}+\frac{1}{\gamma}u_{\gamma})f^*(f'(u_{\gamma}))+\frac{1}{2}p_{\gamma}f''(u_{\gamma})|\nabla u_{\gamma}|^2\leq \\
\int_{\RR^d\times \{0\}} \frac{1}{\gamma}f(u_{\gamma})+\int_{\RR^d\times [0,\tau]} \frac{2}{\gamma}f(u_{\gamma})(n_{\gamma}+\frac{1}{\gamma}u_{\gamma}) +\frac{n_{\gamma}}{2}f^*(f'(u_{\gamma}))+\frac{n_{\gamma}}{2}f\big(\frac{4p_{\gamma}\Delta n_{\gamma}+2\partial_t n_{\gamma}}{n_{\gamma}}\big)+n_{\gamma} f''(u_{\gamma})\frac{p_{\gamma}| \nabla n_{\gamma}|^2}{2n_{\gamma}},
\end{multline}
which is finally in a form that will allow us to estimate.

Fix some $b\leq \frac{b_0}{4\max\big(1,\sup_{\gamma}\norm{p_{\gamma}}_{L^{\infty}([0,\tau]\times\RR^d)}\big)} $ where $b_0$ is the constant from Lemma \ref{lem:dt_n_d2_n}.  If we choose $f$ such that $f$ grows like $\exp(b u)$ at infinity, then $f^*(f'(u))$ grows like $b u\exp(bu)$ at infinity, and hence $f^*(f'(u))$ dominates both $f(u)$ and $f''(u)$ at infinity.   PME has finite propagation in time (uniform in $\gamma$) \cite{vazquez_book}, thus, there exists a radius $R=R_{\tau}>0$ sufficiently large such that $(\rho_{\gamma}, p_{\gamma})$ is supported in $B_R$ independently of $\gamma$.  Recalling that $u_{\gamma}=-\gamma(\Delta p_{\gamma}+n_{\gamma})$ and $n_{\gamma}\geq 0$, it follows that $f(u_{\gamma}), f'(u_{\gamma}), f''(u_{\gamma})$ are all supported on $B_R$ independently of $\gamma$ as well.
Since we are integrating functions with uniformly bounded support and $f\big(\frac{4p_{\gamma}\Delta n_{\gamma}+2\partial_t n_{\gamma}}{n_{\gamma}}\big)$ is bounded by Lemma \ref{lem:dt_n_d2_n} for our choice of $b$, it follows that the left-hand side of (\ref{eq:f_u_6}) dominates the right-hand side and so the result follows.

\end{proof}

It essentially immediately follows that $p$ is a weak supersolution to the appropriate HJB equation.
\begin{cor}\label{cor:p_hjb} Given any $L^2_{\loc}([0,\infty); L^2(\RR^d))$ weak limit point $u_+$ of the family $u_{\gamma, +}$
    $p$ solves, in the sense of weak solutions,
\begin{equation}\label{eq:p_hjb}
    \partial_t p - |\nabla p|^2  +u_+p\geq 0,
\end{equation}
 where for any $\tau > 0$ there exists $b=b(\tau,d) > 0$ such that $(bu_+-1) e^{bu_+} + 1 \in L^1([0, \tau ] ; \R^d)$.
\end{cor}

Although we now know that $p$ is a supersolution to an HJB equation, it is somewhat annoying to directly obtain the Hopf-Lax formula from (\ref{eq:p_hjb}), due to the fact that $p$ is not continuous.  Instead, we will work towards the Hopf-Lax formula by once again going through the $\gamma$ limit.  Here, we will still need to deal with the difficulty that $u_{\gamma, +}$ is not uniformly bounded in $L^{\infty}$.  
We proceed by adapting an argument from \cite{cardaliaguet_graber_mfg}, which provides a method to obtain Hopf-Lax type formulas for Hamilton-Jacobi equations with unbounded coefficients.  A key difference in our setting is that the right-hand side has the specific form $p_{\gamma} u_{\gamma, +}$. This structure allows us to combine their approach with Gronwall-type estimates to obtain much stronger bounds.

\begin{lemma}\label{lem:hj_estimate} Choose a decreasing nonnegative function $\lambda\in L^1([0,t_1-t_0])$.
Given any $\gamma\in (1,\infty)$ and any points $(x_1,t_1), (x_0, t_0)$ with $t_0<t_1$ there exists a constant $C=C(t_1,d)$ such that
\begin{equation}\label{eq:hj_bound}
p_{\gamma}(x_0, t_0)\leq e^{\Lambda_{\gamma}(t_1-t_0)}\Big( p_{\gamma}(x_1, t_1)+\frac{|x_1-x_0|^2}{4\int_0^{t_1-t_0}e^{\Lambda_{\gamma}(s)}\, ds}+C(t_1-t_0)^{7/10}e^{-\lambda(t_1-t_0)}\Big)
\end{equation}
\end{lemma}
where  $b$ is the constant from Proposition \ref{prop:ab} and
\[
 \Lambda_{\gamma}(t) := \frac{5}{4b}\int_{0}^t \lambda(a)\, da+\frac{1}{b}\int_{0}^t\log(1+\norm{\exp(b u_{\gamma,+})-1}_{L^1(\{t_1-a\}\times\RR^d)}) \, da .
\]

\begin{proof}
Define $\vp_{\gamma}(x,t):=p_{\gamma}(x, t_1-t)$. It then follows that $\vp_{\gamma}$ satisfies the differential inequality
\[
\partial_t \vp_{\gamma}(x,t)+|\nabla \vp_{\gamma}(x,t)|^2\leq  \vp_{\gamma}(x,t) u_{\gamma,+}(x,t_1-t)
\]
almost everywhere.
Define 
$$\bar{\lambda}_{\gamma}(s):=\lambda(s)+\frac{1}{b}\log(\norm{\exp(b u_{\gamma,+})-1}_{L^1(\{a\}\times\RR^d)})$$ and split
\[
\vp_{\gamma}(x,t) u_{\gamma,+}(x,t_1-t)\leq \vp_{\gamma}(x,t)\bar{\lambda}_{\gamma}(t)+ \vp_{\gamma}(x,t)(u_{\gamma}(x,t_1-t)-\bar{\lambda}_{\gamma}(t))_+
\]
Multiplying both sides of the differential inequality by $e^{-\Lambda_{\gamma}(t)}$ we see that
\[
\partial_t (e^{-\Lambda_{\gamma}(t)} \vp_{\gamma})+e^{-\Lambda_{\gamma}(t) }|\nabla \vp_{\gamma}|^2\leq e^{-\Lambda_{\gamma}(t)}\vp_{\gamma}(u_{\gamma}-\bar{\lambda}_{\gamma}(t))_+.
\]
Let $q_{\gamma}(t,x):=e^{-\Lambda_{\gamma}(t)}\vp_{\gamma}(t,x)$, we then have 
\[
\partial_t q_{\gamma}+e^{\Lambda_{\gamma}(t)}|\nabla q_{\gamma}|^2\leq q_{\gamma}(u_{\gamma}-\bar{\lambda}_{\gamma}(t) )_+.
\]

Fix any two points $x_1, x_0\in \RR^d$. We now introduce a family of paths $x_{\sigma}$ in the spirit of the path optimization argument introduced in \cite{cardaliaguet_graber_mfg}. 
For each $\sigma$ in the unit ball $B_1$ let $x_{\sigma}:[0,t_1-t_0]\to \RR^d$ be a path such that $x_{\sigma}(0)=x_1$ and $x_{\sigma}(t_1-t_0)=x_0$.  Consider 
\begin{multline}
\frac{d}{dt}\big[ q_{\gamma}(x_{\sigma}(t),t)-\frac{1}{4}\int_{t_0}^t e^{-\Lambda_{\gamma}(s)}|x_{\sigma}'(s)|^2\, ds\big]=\\
\partial_t q_{\gamma}(x_{\sigma}(t),t)+\nabla q_{\gamma}(x_{\sigma}(t),t)\cdot x_{\sigma}'(t)-\frac{e^{-\Lambda_{\gamma}(t)}}{4}|\xs'(t)|^2\leq q_{\gamma}(\xs(t),t)(u_{\gamma}(,\xs(t)t_1-t)-\bar{\lambda}_{\gamma}(t) )_+
\end{multline}
Thus,
\[
q_{\gamma}(x_0, t_1-t_0)\leq q_{\gamma}(x_1,0)+\frac{1}{4}\int_{0}^{t_1-t_0} e^{-\Lambda_{\gamma}(s)}|x_{\sigma}'(s)|^2+ q_{\gamma}(\xs(s),s)(u_{\gamma}(\xs(s),t_1-s)-\bar{\lambda}_{\gamma}(s))_+\, ds,
\]
It then follows that
\begin{equation}\label{eq:hj_1}
\vp_{\gamma}(x_0, t_1-t_0)e^{-\Lambda_{\gamma}(t_1-t_0)}\leq  \vp_{\gamma}(x_1, 0)+\frac{1}{4}\int_{0}^{t_1-t_0} e^{-\Lambda_{\gamma}(s)}\Big(|x_{\sigma}'(s)|^2+ \vp_{\gamma}(s,\xs(s))(u_{\gamma}(t_1-s,\xs(s))-\bar{\lambda}_{\gamma}(s))_+\Big)\, ds.
\end{equation}
We now assume that $\xs$ has the form 
$$
\xs(s)=\sigma \xi(s)+x_0+z(s)(x_1-x_0),
$$
where $\xi:[0,t_1-t_0]\to[0,1]$ satisfies $\xi(0)=\xi(t_1-t_0)=0$ and $z:[0,t_1-t_0]\to [0,1]$ is an increasing function such that $z(0)=0$ and $z(t_1-t_0)=1$.  For notational simplicity, we will write $\alpha_{\gamma}=\vp_{\gamma}(u_{\gamma}-\bar{\lambda}_{\gamma})_+$.  Averaging (\ref{eq:hj_1}) over $B_1$ we see that 
\[
\vp_{\gamma}(x_0, t_1-t_0)e^{-\Lambda_{\gamma}(t_1-t_0)}\leq  \vp_{\gamma}(x_1, 0)+\frac{1}{4|B_1|}\int_{0}^{t_1-t_0} \int_{B_1} e^{-\Lambda_{\gamma}(s)}\Big(|\sigma|^2|\xi'(s)|^2+|x_1-x_0|^2z'(s)^2+ \alpha_{\gamma}(s, \xs(s))\Big)\, d\sigma\, ds.
\]

The optimality condition for $z$ implies that $(z'(s)e^{-\Lambda_{\gamma}(s)})'=0$, therefore $z'(s)=\frac{e^{\Lambda_{\gamma}(s)}}{\int_{0}^{t_1-t_0} e^{\Lambda_{\gamma}(s)}\, ds}$. Thus, making this choice we see that
\[
\vp_{\gamma}(x_0, t_1-t_0)e^{-\Lambda_{\gamma}(t_1-t_0)}\leq  \vp_{\gamma}(x_1, 0)+\frac{|x_1-x_0|^2}{4\int_0^{t_1-t_0} e^{\Lambda_{\gamma}(s)}\, ds}+\frac{1}{4|B_1|}\int_{0}^{t_1-t_0} \int_{B_1} e^{-\Lambda_{\gamma}(s)}\Big(|\sigma|^2|\xi'(s)|^2+ \alpha_{\gamma}(s, \xs(s))\Big)\, d\sigma\, ds.
\]
Changing variables $y=x_{\sigma}$, it follows that
\[
\frac{1}{|B_1|}\int_{B_1} \alpha_{\gamma}(s,x_{\sigma}(s))d\sigma =\frac{\xi(s)^{-d}}{|B_1|}\int_{B_{\xi(s)}(x_0+z(s)(x_1-x_0))} \alpha_{\gamma}(s,y)dy
\]
where $B_{\xi(s)}(x_0+z(s)(x_1-x_0))$ is the ball of radius $\xi(s)$ centered at $x_0+z(s)(x_1-x_0)$. 

  Using H\"{o}lder's inequality with exponent $2d$, it follows that the above quantity is bounded above by $\xi(s)^{-1/2}\norm{\alpha_{\gamma}}_{L^{2d}(\{s\}\times\RR^d)}.$
Hence, after dropping the good term $e^{-\Lambda_{\gamma}(s)}$ in the last integral we see that
\[
\vp_{\gamma}(x_0, t_1-t_0)e^{-\Lambda_{\gamma}(t_1-t_0)}\leq  \vp_{\gamma}(x_1, 0)+\frac{|x_1-x_0|^2}{4\int_0^{t_1-t_0} e^{\Lambda_{\gamma}(s)}\, ds} +\frac{1}{4}\int_{0}^{t_1-t_0} \big(|\xi'(s)|^2+ \xi^{-1/2}(s)\norm{\alpha_{\gamma}(s,\cdot)}_{L^{2d}(\RR^d)}\big)\, ds.
\]
Fix some $a>0$ and set
\[
\xi(s):=\begin{cases}
    as^{3/4} & \textup{if}\;\; t_0\leq s<(t_1-t_0)/2,\\
    a(t_1-s)^{3/4} & \textup{if}\;\; (t_1-t_0)/2\leq s\leq t_1.
\end{cases}
\]
Using H\"{o}lder's inequality with exponent $2$ on $\int_0^t \xi^{-1/2}(s)\norm{\alpha_{\gamma}(s,\cdot)}_{L^{2d}(\RR^d)}\, ds$, we see that
\[
 \norm{\xi'}_{L^2([t_0,t_1])}^2\leq C (t_1-t_0)^{1/2} a^2, \quad \norm{\xi^{-1/2}}_{L^{2}([0,t_1-t_0])}\leq Ca^{-1/2} (t_1-t_0)^{1/8}
\]
thus,
\[
\vp_{\gamma}(x_0, t_1-t_0)e^{-\Lambda_{\gamma}(t_1-t_0)}\leq  \vp_{\gamma}(x_1, 0)+\frac{|x_1-x_0|^2}{4\int_0^{t_1-t_0} e^{\Lambda_{\gamma}(s)}\, ds} +C((t_1-t_0)^{1/8}a^{-1/2}\norm{\alpha_{\gamma}}_{L^2([0,t_1-t_0];L^{2d}(\RR^d))}+(t_1-t_0)^{1/2}a^2)
\]
Optimizing over $a>0$, we obtain
\begin{equation}\label{eq:hj_almost}
\vp_{\gamma}(x_0, t_1-t_0)e^{-\Lambda_{\gamma}(t_1-t_0)}\leq  \vp_{\gamma}(x_1, 0)+\frac{|x_1-x_0|^2}{4\int_0^{t_1-t_0} e^{\Lambda_{\gamma}(s)}\, ds} +C(t_1-t_0)^{\frac{3}{10}}\norm{\alpha_{\gamma}}_{L^2([0,t_1-t_0];L^{2d}(\RR^d))}^{4/5}
\end{equation}
for a potentially different constant $C>0$.

Finally, it remains to estimate $\norm{\alpha_{\gamma}}_{L^2([0,t];L^{2d}(\RR^d))}$. Recalling that $\alpha_{\gamma}=\vp_{\gamma} (u_{\gamma}-\lambda)_+$, we may write
\[
\norm{\alpha_{\gamma}}_{L^{2d}(\{s\}\times \RR^d))}^{2d}\leq \norm{\vp_{\gamma}}_{L^{\infty}([t_0,t_1]\times\RR^d)}^2 \int_{0}^{\infty} 2dv^{2d-1}|\{x\in \RR^d: u_{\gamma,+}(t_1-s,x)>v+\bar{\lambda}_{\gamma}(s)\}|\, dv
\]
By Chebyshev's inequality, for any strictly increasing function $f:\RR\to\RR$,
\[
\leq \norm{p_{\gamma}}_{L^{\infty}([t_0,t_1]\times\RR^d)}^2\int_{0}^{\infty} 2dv^{2d-1}\frac{\norm{f(u_{\gamma,+})-f(0)}_{L^1(\{t_1-s\}\times\RR^d)}}{f(\bar{\lambda}_{\gamma}(s)+v)-f(0)}\, dv
\]
If we choose $f(a)=\exp(b a)-1$, then we see that
\[
\norm{\alpha_{\gamma}}_{L^{2d}(\{s\}\times \RR^d))}\leq C e^{-b\bar{\lambda}_{\gamma}(s)}\norm{\exp(b u_{\gamma,+})-1}_{L^1(\{t_1-s\}\times \RR^d)}=Ce^{-b\lambda(s)},
\]
where one should note carefully that now $\bar{\lambda}_{\gamma}$ has been replaced by $\lambda$ in the last right-hand term.
Thus, we have the estimate
\[
\norm{\alpha_{\gamma}}_{L^2([t_0,t_1];L^{2d}(\RR^d))}^2\lesssim \int_{0}^{t_1-t_1}e^{-\frac{5}{2}\lambda(s)}\leq (t_1-t_0)e^{-\frac{5}{2}\lambda(t_1-t_0)},
\]
hence,
\[
\norm{\alpha_{\gamma}}_{L^2([t_0,t_1];L^{2d}(\RR^d))}^{4/5}\lesssim  (t_1-t_0)^{2/5}e^{-\lambda(t_1-t_0)}.
\]

Combining our work, we now have 
\begin{equation}\label{eq:hj_almost_2}
\vp_{\gamma}(x_0, t_1-t_0)e^{-\Lambda_{\gamma}(t_1-t_0)}\leq  \vp_{\gamma}(x_1, 0)+\frac{|x_1-x_0|^2}{4\int_0^{t_1-t_0} e^{\Lambda_{\gamma}(s)}\, ds} +C(t_1-t_0)^{7/10}e^{-\lambda(t_1-t_0)},
\end{equation}
for some possibly new constant $C$.  The result follows after replacing $\vp_{\gamma}$ with $p_{\gamma}$ and multiplying both sides by $e^{\Lambda_{\gamma}(t_1-t_0)}$

\end{proof}

Before we can show that the Hopf-Lax formula also holds for the limiting pressure, we first need a Lemma that gives us a pointwise well-defined representative of our weak solution $p$.  The argument is a simple adaptation of a result from \cite{pqm}.
\begin{lemma}\label{lem:p_pointwise}
Suppose that $(\rho, p, n)$ is a weak solution to (\ref{first}-\ref{second}). $p$ can be redefined on a set of measure zero so that
\begin{equation}\label{eq:p_pointwise}
    p(x,t)=\lim_{r\to 0} \frac{1}{r^2|B_r|}\int_{B_r(x)}\int_{t}^{t+r^2} p(y,s)\, ds\, dy
\end{equation}
for all $(x,t)$.  With this definition, $p$ is spacetime upper semicontinuous and for all $x\in \RR^d$ the mapping $t\mapsto p(x,t)$ is continuous from the right.
\end{lemma}
\begin{proof}
    From our control on $u_{\gamma}$ and the relation $\Delta p_{\gamma}=-n_{\gamma}-\frac{1}{\gamma}u_{\gamma}$ it follows that after taking limits, we have 
    \[
    \Delta p\geq -n\geq -n_0
    \]
    in the sense of spacetime distributions.  Therefore, for any $\epsilon>0$, 
    \[
   \Delta \Big(\frac{1}{\epsilon^2}\int_{0}^{\epsilon^2} p(\cdot,t+s)\, ds\Big) \geq -n_0
    \]
    in the sense of space distributions.  Hence, the mean value property for Laplace's equation implies that for all $x\in \RR^d$ and $t>0$ the function
    \[
\phi(r,\epsilon):=\frac{1}{\epsilon^2|B_r|}\int_{B_r(x)}\int_{0}^{\epsilon^2} p(y,t+s)+\frac{n_0}{2d}|y-x|^2\, ds\, dy
    \]
    is non-decreasing with respect to $r$.   We also note that for $0\leq r'\leq r$ we have 
    \[
    \phi(r,r)-\phi(r,r')=\frac{1}{r^2|B_r|}\int_{B_r(x)}\int_{0}^{r^2} p(y,t+s)-p(y,t+(\frac{r'}{r})^2s ) \,ds\, dy
    \]
    \[
    \geq \frac{1}{r^2|B_r|}\int_{B_r(x)}\int_{0}^{r^2} \int_{(\frac{r'}{r})^2s}^s [\partial_t p(y,a)]_-\, da\,  ds\, dy.
    \]
    From Corollary \ref{cor:p_hjb}, it follows that $[\partial_t p]_-$ is bounded in $L^q(\RR^d\times [0,\tau])$ for any $q\in [1,\infty)$.  Therefore,
    \[
     \phi(r,r)-\phi(r,r')\geq -|B_r|^{-1/q}(r^2-r'^2)^{1-1/q}\norm{[\partial_t p]_-}_{L^q([0,\tau]\times\RR^d)}\geq -Cr^{1-(d+1)/q} (r-r')^{1-1/q}
    \]
    for some constant $C>0$. 
    
    By choosing $q>d+1$, we can conclude that there exists a Holder continuous function $g$ such that $r\mapsto \phi(r,r)+g(r)$ is nondecreasing and $g(0)=0$. As a result, $\lim_{r\to 0^+} \phi(r,r)$ must exist for all $(x,t)$. Hence, (\ref{eq:p_pointwise}) is well defined everywhere.  The Lebesgue differentiation theorem also implies that our redefinition only changes $p$ on a set of measure zero. 

    Finally, to see that $p$ is upper semicontinuous, we note that $\lim_{r\to 0} \phi(r,r)=\lim_{r\to 0} \phi(r,r)+g(r)=\inf_{r>0} \phi(r,r)+g(r)$.  Thus, we may write 
    \[
     p(x,t)=\inf_{r>0 } g(r)+ \frac{1}{r^2|B_r|}\int_{B_r(x)}\int_{t}^{t+r^2} p(y,s)\, ds\, dy.
    \]
The infimum over a family of functions always produces an upper semicontinuous function, hence, $p$ is upper semicontinuous.

\end{proof}

At last we obtain the main result of this Section, the Hopf-Lax formula for our limit pressure $p$.
\begin{prop}\label{prop:hj_estimate}
Given any points $(x_1,t_1), (x_0, t_0)$ with $t_0<t_1$ and a decreasing function $\lambda\in L^1([0,t_1-t_0])$, there exists a constant $C=C(t_1,d)$ such that
\begin{equation}\label{eq:hj_bound_limit}
p(x_0, t_0)\leq e^{\Lambda(t_1-t_0)}\Big( p(x_1, t_1)+\frac{|x_1-x_0|^2}{4\int_0^{t_1-t_0}e^{\Lambda(s)}\, ds}+C(t_1-t_0)^{7/10} e^{-\lambda(t_1-t_0)}\Big)
\end{equation}
\end{prop}
where  $b$ is the constant from Proposition \ref{prop:ab} and
\[
\Lambda(t):=\frac{5}{4b}\int_0^t \lambda(s)\, ds+\frac{t}{b}\log(1+\frac{C}{t})
\]
\begin{proof}
Using the formula from Lemma \ref{lem:p_pointwise}, we have 
\[
p(x_0,t_0)=\lim_{r\to 0^+} \frac{1}{r^2|B_r|}\int_{B_r(x)}\int_{t}^{t+r^2} p(y,s)\, ds\, dy.
\]
Choose a point $x_2\in \RR^d$ and $t_2>t_0$ such that $p_{\gamma}(x_2,t_2)$ converges to $p(x_2,t_2)$ along some subsequence $\gamma_k$.
Using the $L^2_tH^1_x$ strong convergence of $p_{\gamma}$ to $p$ and then applying Lemma \ref{lem:hj_estimate}, we have for any $x_2\in \RR^d$ and $t_2>t_0$
\[
p(x_0,t_0)=\lim_{r\to 0^+}\lim_{k\to \infty} \frac{1}{r^2|B_r|}\int_{B_r(x)}\int_{t_0}^{t_0+r^2} p_{\gamma_k}(y,s)\, ds\, dy\leq 
\]
\[
\lim_{r\to 0^+}\lim_{k\to \infty} \frac{1}{r^2|B_r|}\int_{B_r(x)}\int_{t_0}^{t_0+r^2}  e^{\Lambda_{\gamma_k}(t_2-s)}\Big( p_{\gamma_k}(x_2, t_2)+\frac{|x_2-x_0|^2}{4\int_0^{t_2-s}e^{\Lambda_{\gamma_k}(a)}\, da}+C(t_2-t_0)^{7/10} e^{-\lambda(t_2-s)}\Big)  \, dy\, ds. 
\]

Recall that
\[
 \Lambda_{\gamma}(t) := \frac{5}{4b}\int_{0}^t \lambda(a)\, da+\frac{1}{b}\int_{0}^t\log(1+\norm{\exp(b u_{\gamma,+})-1}_{L^1(\{t_1-a\}\times\RR^d)}) \, da .
\]
Applying Jensen's inequality, we have the bound
\[
 \Lambda_{\gamma}(t) \leq \frac{5}{4b}\int_{0}^t \lambda(a)\, da+\frac{t}{b}\log(1+\frac{1}{t}\norm{\exp(b u_{\gamma,+})-1}_{L^1([t_1-t, t_1]\times\RR^d)}) .
\]
Hence, we can find a potentially new constant $C=C(\tau,d)>0$ such that
\[
 \Lambda_{\gamma}(t) \leq   \frac{5}{4b}\int_{0}^t \lambda(a)\, da+\frac{t}{b}\log(1+\frac{C}{t})=\Lambda(t).
\]
for all $\gamma$.  Therefore, 

\[
 p(x_0,t_0)\leq e^{\Lambda(t_2-t_0)}\Big( p(x_2, t_2)+\frac{|x_2-x_0|^2}{4\int_0^{t_2-t_0}e^{\Lambda(s)}\, ds}+C(t_2-t_0)^{7/10} e^{-\lambda(t_2-t_0)}\Big).
\]
Since $p_{\gamma}$ converges pointwise almost everywhere to $p$ along appropriate subsequences, it follows that 
\[
p(x_0,t_0)\leq e^{\Lambda(t_2-t_0)}\Big( p(x_2, t_2)+\frac{|x_2-x_0|^2}{4\int_0^{t_2-t_0}e^{\Lambda(a)}\, da}+C(t_2-t_0)^{7/10} e^{-\lambda(t_2-t_0)}\Big).
\] 
for a dense set of $(x_2, t_2)$ with $t_2>t_0$.  The result now follows from the upper semicontinuity of $p$.
\end{proof}

\section{H\"{o}lder continuity of the Hitting time}\label{sec:barrier}

  We are now going to construct a radial supersolution  which will give an upper bound for the rate of expansion for the tumor, and thus, will yield a lower bound on the arrival times.  Given a point of interest $x_0\in \RR^d$,    our supersolution will be defined on the time-dependent annulus
\[
A(t):=\{t\}\times \{ x: r(t) \leq |x-x_0| \leq m r(t)\} , \quad A=\bigcup_{t\in [0,\epsilon]} A(t),
\]
for some $m>1$ and a function $r(t)\geq 0$ that we will define shortly. 
Given some starting time $t_0$, for each $t, r\geq 0$ we define
\[
\bar{p}(t,r):=\sup_{x\in B_{r}(x_0)} p(t+t_0,x),
\]
where the sup is well defined since $p$ is upper semicontinuous in space.

Now we will construct our supersolution $\psi(t,x)$ by solving
\[
\begin{cases}
    -\Delta \psi(t,x)=\bar{n}_0& \textup{if} \; r(t)<|x-x_0|<mr(t),\\
    \psi(t,x)=0 & \textup{if} \;|x-x_0|\leq r(t),\\
    \psi(t,x)=\bar{p}(t,|x-x_0|) & \textup{if} \; |x-x_0|\geq mr(t).\\
\end{cases}
\]

On $A(t)$, the equation admits the explicit radial solution
\begin{equation}\label{eq:supersolution_psi}
\psi(t,x)=h(t)\Gamma_d(|x-x_0|)-\frac{\bar{n}(0)}{2d}|x-x_0|^2+g(t),
\end{equation}
where $\Gamma_d$ is the fundamental solution of the Laplace equation in dimension $d$, i.e. $\Gamma_d'(r)=r^{1-d}$,
\begin{equation}\label{eq:super_solution_h}
h(t):=\frac{\bar{p}(t,mr(t))+(2d)^{-1}\bar{n}(0)(m^2-1)r(t)^2}{\Gamma_d(mr(t))-\Gamma_d(r(t))},
\end{equation}
and 
\begin{equation}
    g(t):=\frac{\bar{n}_0}{2d}-h(t)\Gamma_d(r(t)).
\end{equation}
Finally, we define $r(t)$ by choosing some initial data $r(0)$ and then solving the ODE

\begin{equation}\label{eq:radius_ode}
r'(t)=-|\nabla \psi(t,y)|
\end{equation}
where the right hand side is evaluated at any point $y$ such that $|y-x_0|=r(t)$. 

We now show that $\psi$ is indeed a supersolution as long as comparison holds at initial time. Due to the lack of regularity for the pressure variable, we establish comparison using the time integrated versions of $\psi$ and $p$.  This creates an annoying issue where it is difficult to establish that the boundary data stays ordered as the annulus moves.  To avoid this problem, we establish comparison by first going through a sequence of supersolutions $\psi_k$, where the $\psi_k$ are defined on modified annuli whose outer radii are taken to be piecewise constant in time. 

\begin{lemma}\label{lem:comparison}
Let $\mu(t,x)$ be the characteristic function of the set $\{x\in \RR^d: |x-x_0|\geq r(t)\}$.  
If $\mu(0,x)\leq \rho(t_0,x)$ for almost every $x\in \RR^d$,  then $p(t_0+t,x)\leq \psi(t,x)$ for almost every $x\in \RR^d$ and almost every time $t\geq 0$
\end{lemma}

\begin{proof}
As we noted above, we will first prove the comparison for a modified sequence of supersolutions $\psi_k$.  The $\psi_k$ will be defined in precisely the same way as $\psi$, except that we will modify the construction of the moving annulus.  Hence, given radii $r_k(t)< R_k(t)$, we define $\psi_k$ by solving
\[
\begin{cases}
    -\Delta \psi_k(t,x)=\bar{n}_0& \textup{if} \; r_k(t)<|x-x_0|<R_k(t),\\
    \psi_k(t,x)=0 & \textup{if} \;|x-x_0|\leq r_k(t),\\
    \psi_k(t,x)=\bar{p}(t,|x-x_0|) & \textup{if} \; |x-x_0|\geq R_k(t).\\
\end{cases}
\]
$r_k(t)$ will be defined as before via the ODE $r_k'(t)=-|\nabla \psi_k(t,y)|$ where $y$ is any point satisfying $|y-x_0|=r_k(t)$.  We then define $R_k$ by setting 
\[
R_k(t):=mr_k(t_{k,j}), \quad \textup{if} \; t\in [t_{k,j}, t_{k, j+1}),
\]
where we inductively define the points $t_{k,j}$ by setting $t_{k,0}=0$ and then taking
\[
t_{k,j+1}:=\inf \{ t\geq t_{k,j}: r_k(t)<(1-\frac{1}{k+1})r_k(t_{k,j})\}.
\]

As before, on the annulus $r_k(t)\leq |x-x_0|\leq R_k(t)$, the $\psi_k$ will admit the explicit radial solutions 
\begin{equation}\label{eq:supersolution_psi_k}
\psi_k(t,x)=h_k(t)\Gamma_d(|x-x_0|)-\frac{\bar{n}(0)}{2d}|x-x_0|^2+g_k(t),
\end{equation}
where 
\begin{equation}\label{eq:super_solution_h_k}
h_k(t):=\frac{\bar{p}(t,R_k(t))+(2d)^{-1}\bar{n}(0)(R_k(t)^2-r_k(t)^2)}{\Gamma_d(R_k(t))-\Gamma_d(r_k(t))},
\end{equation}
and
\begin{equation}\label{eq:super_solution_g_k}
    g_k(t):=\frac{\bar{n}_0}{2d}-h_k(t)\Gamma_d(r_k(t)).
\end{equation}

Let $\Psi_k(t,x)=\int_0^t \psi_k(s,x)\, ds$.  Since $\psi_k$ is clearly Lipschitz in space on $|x-x_0|\leq R_k(t)$, it follows that $\Psi_k$ is Lipschitz in space on $|x-x_0|\leq R_k(t)$ and
\[
\nabla \Psi_k(t,x)=\int_0^t \nabla \psi_k(s,x)\, ds
\]
almost everywhere on $|x-x_0|\leq R_k(t)$.  Define $\tilde{t}_k(r)$ to be the inverse function of $r_k(t)$.  From the definition of $\psi_k$, it follows that
\[
\nabla \Psi_k(t,x)=\int_{\min(t, \tilde{t}_k(|x-x_0|))}^t \nabla \psi_k(s,x)\, ds.
\]
Now if $x$ is a point such that $|x-x_0|<R_k(t)$ and $|x-x_0|<r(0)$, then for each fixed $s\in (\tilde{t}(|x-x_0|), t]$, there exists a neighborhood of $x$ such that $\nabla \psi_k(s,x)$ is differentiable and $-\Delta \psi_k(s,x)=\bar{n}_0$.  Thus, it follows that
\[
-\Delta \Psi_k(t,x)=\sgn_+\big(t-\tilde{t}(|x-x_0|)\big)\tilde{t}'(|x-x_0|)|\nabla \psi_k(\tilde{t}(|x-x_0|),x)|+ \int_{\min(t, \tilde{t}_k(|x-x_0|))}^t \bar{n}_0\, ds.
\]
Since $r'_k(\tilde{t}_k(|x-x_0|))=-|\nabla \psi(\tilde{t}(|x-x_0|), x)|$ and $\tilde{t}_k(r)$ is the inverse of $r_k(t)$, we see that
\[
-\Delta \Psi_k(t,x)=-\sgn_+\big(t-\tilde{t}(|x-x_0|)\big)+ (t- \tilde{t}(|x-x_0|)_+\bar{n}_0.
\]
On the other hand, if $x$ is a point such that $r(0)<|x-x_0|<R_k(t)$, then 
\[
-\Delta \Psi_k(t,x)=t \bar{n}_0.
\]

Let $\mu_k(t,x)$ be the characteristic function of the set $\{(t,x): |x-x_0|\geq r_k(t)\}$ and note that $\mu_k(t,x)=\sgn_+(t-\tilde{t}_k(|x-x_0|))$ and $\int_0^t \mu_k(s,x)\, ds=(t- \tilde{t}_k(|x-x_0|)_+.$
Combining our work from above, we can conclude that for almost every $x$ satisfying $|x-x_0|<R_k(t)$ we have 
\[
-\Delta \Psi_k(t,x)=\mu_k(0,x)-\mu_k(t,x)+\int_0^t \mu_k(s,x) \bar{n}_0,
\]
and for almost all $x\in \RR^d$ we have
 $\Psi_k(t,x)(1-\mu_k(t,x))=0,$ as well as $\psi_k(1-\mu_k(t,x))=0$.

Now let us define the time shifted variables $\tilde{w}(t,x):=\int_0^t p(s+t_0,x)\, ds$, $\tilde{\rho}(t,x):=\rho(t_0+t,x)$, and $\tilde{n}(t,x):=n(t_0+t,x)$.  It then follows that $-\Delta \tilde{w}(t,x)=\tilde{\rho}(t,x)-\tilde{\rho}(0,x)+\int_{t_0}^t \tilde{\rho}(s,x) \tilde{n}(s,x)$ and $(1-\rho(t+t_0,x))w_0=0$ almost everywhere.  For any time $t\in [0,t_{k,1})$, the definition of $\psi_k$ guarantees that $\Psi_k(t,x)\geq w(t,x)$ for all $x$ satisfying $|x-x_0|=R_k(t)=mr(0)$.  Hence, for any $t\in [0, t_1)$ and any increasing $C^1$ function $\eta:\RR\to\RR$ such that $\eta(a)=0$ if $a\leq 0$, we have
\[
\int_{\{|x-x_0|\leq R_k(0)\}} (\tilde{\rho}-\mu_k)\eta(\tilde{w}-\Psi_k)+\eta'(\tilde{w}-\Psi_k)|\nabla (\tilde{w}-\Psi_k)|^2\leq \int_{\{|x-x_0|\leq R_k(0)\}}\int_0^t \eta(\tilde{w}-\Psi_k)(\tilde{\rho}\tilde{n}-\mu_k\bar{n}_0)
\]
Letting $\eta$ approach $\sgn_+$ and using the fact that $\sgn_+(\tilde{w}-\Psi_k)=\sgn_+(\tilde{\rho}-\mu_k)$, we can conclude that
\[
\int_{\{|x-x_0|\leq R_k(0)\}} (\tilde{\rho}-\mu_k)_+\leq \int_{\{|x-x_0|\leq R_k(0)\}}\bar{n}_0\int_0^t(\tilde{\rho}-\mu_k)_+.
\]
Hence, Gronwall's inequality now implies that $\tilde{\rho}(t,x)\leq \mu_k(t,x)$ for all $t\in [0,t_{k,1})$ and almost all $x\in \RR^d$ (recall it is immediate that $\tilde{\rho}\leq \mu_k$ on $|x-x_0|\geq R(0)$ from the definition of $\mu_k$).  The masses of the differences $\mu_k(t,x)-\mu_k(0,x)$ and $\tilde{\rho}(t,x)-\tilde{\rho}(0,x)$ are continuous functions of time, therefore, the ordering $\tilde{\rho}\leq \mu_k$ must hold at time $t_1$.  This allows us to run the above argument on $[t_{k,1}, t_{k,2})$.  Iterating, we conclude that the ordering $\tilde{\rho}\leq \mu_k$ must hold for all times $t$ when $r_k(t)>0$.  

Now we wish to argue that $\liminf_{k\to\infty} r_k(t)\geq r(t)$. Let 
\[
t_*=\inf\{t>0: \liminf_{k\to\infty} r_k(t)<r(t)\},
\]
and note that $\liminf_{k\to\infty} r_k(t_*)=r(t_*)$.  Using the explicit formulas (\ref{eq:supersolution_psi}) and (\ref{eq:supersolution_psi_k}), as well as the upper semicontinuity of $r\mapsto \bar{p}(t,r)$, it follows that
\[
\liminf_{k\to\infty} r_k'(t_*)\geq r'(t_*)
\]
whenever $r(t_*)>0$.  Hence, $r(t)\leq \liminf_{k\to\infty} r_k(t)$ for all times where $r(t)>0$.  This implies that $\tilde{\rho}(t,x)\leq \mu(t,x)$ for all $t$ 
and almost all $x$.

Finally, we note that the ordering $\tilde{\rho}(t,x)\leq \mu(t,x)$ implies that
for almost every time $t$
\[
(p-\psi)_+(\Delta p+n)=0, \quad (p-\psi)_+(\Delta \psi+\bar{n}_0)=0.
\]
distributionally. Thus, for any $T>0$
\[
\int_{Q_T} |\nabla (p-\psi)_+|^2=\int_{Q_T} (p-\psi)_+(n-\bar{n}_0)
\]
which is only possible if $(p-\psi)_+=0$ almost everywhere.

\end{proof}

We can now use this barrier supersolution to get bounds on the H\"{o}lder continuity of the hitting time. The key 
is to use our Hopf-Lax estimate from \ref{prop:hj_estimate} to ensure that the supersolution arrives at the point of interest at the correct time.

\begin{theorem}\label{thm:hitting_time_holder}

$T$ is locally H\"{o}lder continuous on the set $\{x\in \RR^d: 0<T(x)<\infty\}$.  In particular, for any $x_1\in \RR^d$ such that $T(x_1)\in (0,\infty)$, we have
\begin{equation}\label{eq:holder_bound}
   \sup_{y\in B_R(x_1)} T(x_1)-T(y)\lesssim R^{\alpha_d}
    \end{equation}
    for all $R>0$ sufficiently small, where
    \begin{equation}\label{eq:holder_exponent}
            \alpha_d:=\begin{cases}
    \frac{2}{e} & \textup{if}\; d=2,\\
   2(\frac{2}{d})^{\frac{d}{d-2}} & \textup{if}\; d>2.\\
    \end{cases}
    \end{equation}
    
\end{theorem}
\begin{proof}
Let $\epsilon>0$ be a small value that we will choose later.
Let $$\delta=\delta(\epsilon):=\inf\{R>0: \sup_{y\in B_R(x_1)} T(x_1)- T(y)\geq \epsilon\}.$$
Since $T$ is continuous at $x_1$, it follows that $\lim_{\epsilon\to 0}\delta(\epsilon)=0$.  Let $t_0=T(x_1)-\epsilon$ and $t_1=T(x_1)$.  Thanks to the super solution that we have constructed above, we know that 
\[
\inf_{y\in B_{r(t)}(x_1)} T(y)\geq t_0+t,
\]
which implies 
\[
\sup_{y\in B_{r(t)}(x_1)} T(x_1)-T(y)\leq (t_1-t_0-t).
\]
Hence, if we can provide lower bounds on $r(t)$ in terms of $t$, we can get a H\"{o}lder estimate for $T$ at $x_1$.  In particular, a bound of the form $(t_1-t_0-t)^{1/\alpha}\lesssim r(t)$ will imply that
 $\sup_{y\in B_R(x_1)} T(x_1)-T(y)\lesssim R^{\alpha}$.

To bound $r(t)$ from below, we must consider the ODE (\ref{eq:radius_ode}),  which can be simplified to $$r'(t)=-|h(t)||\Gamma_d'(r(t))|-\frac{\bar{n}(0)}{d}r(t).$$  Noting that in any dimension there exists a function $\xi_d(m)$ such that $\frac{|\Gamma_d'(r(t))|}{|\Gamma_d(m(r(t))-\Gamma_d(r(t))|}=r(t)^{-1}\xi_d(m)$, it follows from the structure of $h$ and the ODE that there exists some constant $K>0$ such that
\begin{equation}\label{eq:ode_1}
r'(t)+Kr(t)\geq -\frac{\bar{p}(t, m r(t))\xi_d(m)}{r(t)}.
\end{equation}

Now we want to estimate $\bar{p}(t, m r(t))$.
To do so, we will apply the bounds from Lemma \ref{lem:hj_estimate}, choosing to evaluate $p$ at $(x_1, t_1)$ and leaving the choice of $\lambda\in L^1([0,t_1-t_0])$ until later. With these choices, we see that
$$
\begin{array}{lll}
\bar{p}(t,mr(t))&=&\sup_{x\in B_{m r(t))}(x_1)} p(t+t_0,x)\\
&\leq& 
\sup_{x\in B_{m r(t))}(x_1)} 
H(t)|x-x_1|^2 +F(t)=m^2r^2 H(t)+F(t),\label{eq:p_upper_bound}
\end{array}
$$
where we have defined 
\begin{equation}\label{eq:H_and_F}
H(t):=e^{\Lambda(t_1-t_0-t)}(4\int_{0}^{t_1-t_0-t} e^{\Lambda(s)}\, ds)^{-1}, \quad   F(t):=C(t_1-t_0-t)^{7/10}e^{-\lambda(t_1-t_0-t)+\Lambda(t_1-t_0-t)}  
\end{equation}  for notational convenience.

Returning to equation (\ref{eq:ode_1}) and applying the upper bound on $\bar{p}$ obtained above, we have 
\[
r'(t)+r(t)(K+m^2\xi_d(m) H(t))\geq -\xi_d(m)\frac{F(t)}{r(t)}
\]
Multiplying both sides by $2r(t)$ and defining $z(t)=r(t)^2$, we get 
\begin{equation}\label{eq:z_ode_1}
z'(t)+z(t)(2K+2m^2\xi_d(m) H(t))\geq -\xi_d(m)F(t).
\end{equation}

Now we choose $m$ by optimizing $m^2\xi_d(m)$.  Define $$\xi_d:=\inf_{m>1}\frac{m^2}{2}\xi_d(m).$$ One can then check that $\xi_d=(\frac{d}{2})^{\frac{d}{d-2}}$ and $\argmin \frac{m^2}{2}\xi_d(m)=(\frac{d}{2})^{\frac{1}{d-2}}$ (where these should be understood in a limiting sense when $d=2$).  Thus, we have 
\begin{equation}\label{eq:z_ode_2}
z'(t)+z(t)(2K+4\xi_d H(t))\geq -dF(t).
\end{equation}
Let $\bar{H}(t)=\int_0^t 4H(s)\, ds$.  Multiplying both sides of (\ref{eq:z_ode_2}) by $\exp(2Kt+\xi_d\bar{H}(t))$  and integrating in time, we can conclude that
\begin{equation}\label{eq:z_ode_3}
z(t)e^{2Kt+\xi_d\bar{H}(t)}\geq z(0) -d\int_0^t F(s)e^{2Ks+\xi_d \bar{H}(s)}\, ds.
\end{equation}

Now we need to provide upper bounds on $\exp(\xi_d\bar{H}(t))$. To do so, we will need to make a choice for $\lambda$.  Fix some $\theta>0$ and set
\[
\lambda(s)=\theta+s^{-1/2}.
\]
We then have 
\[
\Lambda(t)=\frac{5}{4b}(\theta t+2t^{1/2})+\frac{t}{b}\log(1+\frac{C}{t}).
\]

Using the above estimates, we see that
\[
4H(t)\leq \frac{\exp\Big(\frac{5}{4b}(\theta (t_1-t_0-t)+2(t_1-t_0-t)^{1/2})+(t_1-t_0-t)\log(1+C/(t_1-t_0-t))\Big)}{\int_0^{t_1-t_0-t} e^{\frac{5}{4b}\theta s}\, ds}=
\]
\[
\frac{\frac{5\theta}{4b}\exp\Big(2(t_1-t_0-t)^{1/2})+(t_1-t_0-t)\log(1+C/(t_1-t_0-t))\Big)}{1- e^{-\frac{5}{4b}\theta (t_1-t_0-t)}}.
\]
Since $(t_1-t_0-t)\leq (t_1-t_0)=\epsilon$, we can assume that $\epsilon$ is sufficiently small that 
\[
4H(t)\leq \frac{5\theta}{4b(1- e^{-\frac{5}{4b}\theta (t_1-t_0-t)})}+\frac{20\theta(t_1-t_0-t)^{1/2}}{4b(1- e^{-\frac{5}{4b}\theta (t_1-t_0-t)})}
\]
Hence, for some possibly new constant $C>0$ independent of $\epsilon$ and $\theta$ we get
\[
\bar{H}(t)\leq \log\big(\frac{e^{\frac{5\theta}{4b}\epsilon}-1}{e^{\frac{5\theta}{4b}(t_1-t_0-t)}-1}\big)+C(1+\epsilon^{3/2}\theta).
\]
Thus, 
\[
\exp(\xi_d\bar{H}(t))\lesssim \big(\frac{e^{\frac{5\theta}{4b}\epsilon}-1}{e^{\frac{5\theta}{4b}(t_1-t_0-t)}-1}\big)^{\xi_d}\exp(C\epsilon^{3/2}\theta)
\]

Now we move to estimating $F(s)\exp(2Ks+\xi_d \bar{H}(s))$. From our choice of $\lambda$, it is clear that $\lambda(t)\geq 2\Lambda(t)$ for all $t$ sufficiently small. Thus, it follows that
\[
F(s)\lesssim (t_1-t_0-s)^{7/10}e^{-\lambda(t_1-t_0-s)/2},
\]
hence, we have the bound
\[
F(s)\exp(2Ks+\xi_d \bar{H}(s))\lesssim \big(\frac{e^{\frac{5\theta}{4b}\epsilon}-1}{e^{\frac{5\theta}{4b}(t_1-t_0-s)}-1}\big)^{\xi_d}(t_1-t_0-s)^{7/10}\exp\Big(\frac{1}{2}\big(4Ks+(2C\epsilon^{3/2}-1)\theta-(t_1-t_0-s)^{-1/2})\big)\Big)
\]
Therefore, once $\epsilon$ is small enough that $2C\epsilon^{3/2}+\frac{5}{4b}\epsilon<1$ we can choose $\theta$ large enough that
\[
z(0) -d\int_0^{\epsilon} F(s)e^{2Ks+\xi_d \bar{H}(s)}\, ds\geq z(0)/2,
\]
and from there we can conclude that
\[
z(0)\lesssim z(t)e^{2Kt+\xi_d\bar{H}(t)}
\]
for all $t\in [0, \epsilon]$.
This implies that
\[
z(0)(\theta(t_1-t_0-t))^{\xi_d}\lesssim  z(t)=r(t)^2.
\]
Hence, $r(0)(t_1-t_0-t)^{\xi_d/2}\lesssim r(t).$  The result now follows from the fact that $\alpha_d=2/\xi_d$.

\end{proof}

 \section{Results from Obstacle Problem Theory}

In this section, we use techniques from the theory of the obstacle problem to study the local behavior of the interface. The main technique here is the quadratic blowup, which classifies free boundary points into regular points, where the zero set is asymptotically a half-space, and singular points, where the zero set is asymptotically lower dimensional. With sufficiently regular source term, the blowup limit approximates the solution at a uniform scale, and we can use this to extract information on the local geometry of the positive set. The regularity of the source term in the equation satisfied by $w$ is governed by the regularity of the nutrient and the regularity of the hitting time. Since the nutrient enjoys parabolic regularity as in Lemma \ref{lem:n_regularity}, H\"older continuity of the hitting time leads to  H\"older continuous source, which is enough to control the blowup limit at both types of free boundary points.

Thus, using that $T\in C^{0,\alpha}_{\loc}(\mathcal{O})$ for the $\alpha\in (0,1)$ from Theorem \ref{thm:hitting_time_holder}, we show that the regular points form an open set of full measure in the spacetime interface, on which $T$ improves to locally Lipschitz and the spatial interface evolves as a locally $C^{1,1-}$ graph. The scale and bounds for which this regularity is achieved can be quantified in terms of the H\"older seminorm of $T$ and the scale at which the zero set achieves sufficiently large density near the regular point. We also show H\"older regularity of the unit normal to the interface in spacetime, using the spatial regularity of the interface and the monotonicity of its expansion. Under the stronger assumption that no singular points occur for some time interval, this lets us improve $T$ to $C^{1,1/2-}_{\loc}$ on the corresponding region.

As for singular points, we show that they form a relatively closed set in $\mathcal{O}$ contained in a $C^1$ manifold of dimension $d-1$. This improves on the standard obstacle problem result that the singular points at a fixed time are contained in a $C^1$ manifold of dimension $d-1$ and implies that the worse case, where the singular points have positive $d-1$ Hausdorff measure for some time, only occurs for at most countably many times. Under the additional assumption that $T$ is Lipschitz up to the singular set, we show a stronger generic regularity result which gives that for a.e. time the singular points have $d-2$ Hausdorff measure 0. In dimension 2, this would imply that the times with singular points have zero measure, as a relatively closed subset of $(0,\infty)$.

We note that we are not currently able to prove that $T$ is Lipschitz up to singular points. It is not clear to what extent the obstacle problem can be leveraged to understand the geometry of the patch at times just before a singular point occurs, in order to prove nondegeneracy of the pressure. Nondegeneracy at later times is also uncertain, but appears more tractable since the blowup is available. For example, suppose one knew, for a singular point $x_0$ and all sufficiently small $r$, that the set $\{ w(\cdot, T(x_0)) = 0 \}\cap B_r(x_0)$ is contained in a strip of width $Cr^{1 + \alpha}$, for some $C, \alpha$ depending on the source term. Assuming such a strip condition, then the Hopf lemma could be applied to establish nondegeneracy of $p$ near $x_0$ at times $t\geq T(x_0)$. This strip condition has been proven by \cite{figalli_serra} for singular points in the $(d-1)$-dimensional stratum, albeit using methods which require much stronger regularity than $C^{0,\alpha}$ source. A related result on the rate of convergence of the quadratic blowup at singular points in dimension 2 has been proven by \cite{colombo}. As far as we are aware, it is not currently known whether or in what sense this strip condition may hold for the obstacle problem with $C^{0,\alpha}$ source.

Nevertheless, we expect that $T$ is indeed Lipschitz, as it seems unlikely that the pressure sometimes becomes degenerate, but only at instances of merging or topological change. It is clear that we cannot hope for better than Lipschitz, since $T$ cannot be differentiable when two pieces of the boundary collide while traveling at different speeds. We also discuss examples in Remark \ref{rem:regular_points_nondifferentiable} which show that Lipschitz continuity of $T$ is sharp at regular points without the additional assumption to give global control over singular points.

For the obstacle problem, the $C^{1,\alpha}$ regularity of the free boundary at regular points and the $C^1$ manifold covering singular points are well-known for H\"older source (see appendix for more details). Thus, the main challenge in lies in the analyzing the time-indexed family of obstacle problems satisfied by $w(\cdot, t)$ to control these properties in time. We note that such parameterized families have now been studied extensively for the constant source obstacle problem with varying fixed boundary data, mainly with the goal of understanding generic behavior of singular points (\cite{monneau}, \cite{figalli_generic}). Our problem differs in that we must contend with a varying low regularity source and no fixed boundary data, which rules out many of the techniques typically used. In particular, the results of \cite{figalli_generic}, including that the singular set has $(d-4)$-Hausdorff measure zero, do not appear to be in reach with even Lipschitz source. Our approach draws from arguments in \cite{monneau} to establish the $C^1$ manifold property for the singular set. However, whereas comparison arguments with the fixed boundary data allow Monneau to prove directly that the hitting time is Lipschitz, we must work harder to get lower regularity for $T$. Finally, the analysis of the regular set for the time-parameterized family, to our knowledge, is new. The main facts we make use of are the H\"{o}lder continuity of $T$, the spacetime continuity of $w$, the $L^1$ time-continuity of $\rho$, and the monotonicity of $\rho$ and $w$ in time.

We remark that Proposition \ref{prop:hitting_time_continuous} shows that the interface strictly expands, and in space-time is exactly the graph of $T$. Therefore, regularity improvements to $T$ correspond exactly to regularity of the space-time interface as a $d$-dimensional manifold. As a result, we will generally not consider the space-time perspective directly, preferring to work with the subsets of $\R^d$ traced out by the moving interface.

\bigskip

Now, we proceed to study the local situation at the free boundary. As we noted above, the new regularity of $T$ from Theorem \ref{thm:hitting_time_holder} feeds back into the obstacle problem satisfied by $w$ through the dependence of $\eta$, as defined in (\ref{eq:w_obst_eqn}), on $T$. We state this precisely below:

\begin{lemma}\label{lem:eta_regularity}
    Up to $C^{1,1-}$, $\eta(\cdot, t)$ has the same spatial regularity as $T$ on $\overline{\{ w(\cdot, t) > 0\}}$. In particular, for any $\tau > 0$, we have $\eta\in L^\infty_t C^{0,\alpha}_x([0, \tau]; \R^d)$.
\end{lemma}
\begin{proof}
    The first part follows immediately from the $L^\infty_t C^{1,1-}$ regularity of $n$, from Lemma \ref{lem:n_regularity}. The second part follows from the $C^{0,\alpha}$ regularity of $T$.
\end{proof}

The exact regularity of $\eta$ is relevant for determining the spatial regularity of the free boundary near regular points. However, for most results in this section, we only require H\"{o}lder continuity to give uniqueness of the quadratic blowup limit introduced in Lemma \ref{lem:quad_blowup_cpt}, and to give spatial equicontinuity of $\eta(\cdot, t)$. The uniqueness of the blowup limit and its subsequent characterization is best expressed as the following dichotomy, originally due to Caffarelli:

\begin{lemma}\label{lem:caffarelli_dichotomy}
Let $u$ be a solution of the obstacle problem $\Delta u = f\chi_{\{u>0\}}$ in $\R^d$ with $f$ positive and $C^{0,\alpha}$ near 0. If $0\in \partial \{ u > 0 \}$, then one of the following holds:
\begin{enumerate}
    \item $\{ u = 0\}$ has density $\frac{1}{2}$ at 0, and the quadratic rescalings $r^{-2}u(rx)$ converge in $C^{1,1-}(B_1)$ to $\frac{f(0)}{2}(x\cdot e)^2_+$ for some unit vector $e$.
    \item $\{ u = 0\}$ has density 0 at 0, and the quadratic rescalings $r^{-2}u(rx)$  converge in $C^{1,1-}(B_1)$ to $\frac{f(0)}{2}x\cdot D^2 u(0) x$, where $D^2 u(0)$ exists in the classical sense and is a positive semidefinite matrix with trace 1.
\end{enumerate}
Points of the first type are called regular points, and points of the second type are called singular points.
\end{lemma}

This dichotomy was proven in \cite{caffarelli98} for the constant source obstacle problem, with the note that minor modifications could extend the proof to the H\"{o}lder continuous case. An energetic criterion for the dichotomy appears in \cite{weiss}. Careful proofs for the uniqueness of the blowup limit in the H\"{o}lder continuous case are given in \cite{blank} for regular points and \cite{monneau} for singular points.

\medskip
The dichotomy applies to the free boundary of $w(\cdot, t)$ for each $t$. We let $R_t$ denote the regular points of $\partial \Omega_t$, and $\Sigma_t$ denote the singular points of $\partial \Omega_t$, for the obstacle problem solved by $w(\cdot, t)$ at each time. Subsequently, we take 

$$R := \bigcup_{t>0} R_t \hbox{ and } \Sigma := \bigcup_{t>0} \Sigma_t,
$$ so that 
$$
\mathcal{O} =\{0<T(x)<\infty\} = R \cup \Sigma.
$$ Let us also mention that $\{R_t\}_{t>0}$ is a foliation of $R$, and so is $\{\Sigma_t\}_{t>0}$ for $\Sigma$, due to Proposition~\ref{prop:hitting_time_continuous}. We further subdivide singular points into strata by the dimensionality of the zero set; specifically, for $0\leq k\leq d-1$ we denote 
$$
\Sigma^k_t: = \{x\in \partial \Omega_t : \dim \ker D^2 w = k\} \hbox{ and } \Sigma^k := \bigcup_t \Sigma^k_t.
$$

For the obstacle problem, regular points are relatively open in the free boundary (\cite{blank}, Corollary 4.8), and thus singular points form a closed set. This topological control is lost in the union over all times, so a first step is to reestablish that control for our $R$ and $\Sigma$. For this, we use a lemma due to Blank, which allows us to identify regular points by finite-scale behavior.

\begin{lemma}[\cite{blank} Theorem 4.5]\label{lem:modified_regular_pt_criterion}
    Let $u \geq 0$ solve $\Delta u = f\chi_{\{u > 0\}}$ in $B_1(0)$, with $0\in \partial \{ u > 0\}$ and $0 < f\leq 1$. Then there exist universal parameters $\lambda_0, r_0, \tau\in (0,1)$ such that if $\lambda_0 < f$ in $B_1$ and
    $$\frac{|\{x : u(x) = 0\}\cap B_r(x_0)|}{|B_r|} \geq \frac{1}{8}\hbox{ for some } r < r_0,
    $$
    then
    $$\frac{|\{x : u(x) = 0\}\cap B_s(x_0)|}{|B_s|} \geq \frac{3}{8}\hbox{ for all } s < \tau r.
    $$
\end{lemma}
In particular, if the hypothesis of the lemma holds, then $\{ u = 0\}$ has positive density at 0, so 0 is a regular point. This lemma is applicable even when $f$ is less regular than H\"{o}lder, and results in a modified regular-singular dichotomy in that case. Essentially, one may take away that nonuniqueness of the blowup limit in the low regularity setting can occur due to infinite rotation, but not due to any sort of mixing of regular and singular point behavior at different scales. Since the following result only relies on the previous lemma and continuity of $T$, it also holds when $T$ is less regular than H\"{o}lder.

\begin{prop}\label{prop:regular_open}
    $R$ is open. Thus, $\Sigma$ is relatively closed in $\mathcal{O}$.
\end{prop}
\begin{proof}
    Let $\lambda_0, r_0, \tau \in (0,1)$ be the parameters given by Lemma \ref{lem:modified_regular_pt_criterion}. Suppose $x_0\in R$, so that $\Omega_{T(x_0)}$ has density $\frac{1}{2}$ at $x_0$, and thus there exists $r > 0$ such that
    $$\frac{|\{x : w(x, T(x_0)) = 0\}\cap B_r(x_0)|}{|B_r|} \geq \frac{1}{8}
    $$
    Then the result of the lemma is that for all $s < \tau r$,
    \[ \frac{|\{x : T(x) \geq T(x_0)\}\cap B_s(x_0)|}{|B_s|} \geq \frac{3}{8} \]
    Fix $s_0 = \tau r / 2$. By Lemma \ref{lem:rho_regularity}, $|\Omega_t|$ is continuous in $t$, so we can choose $\delta$ for $s_0$ such that
    \[ \frac{|\{x : T(x) > T(x_0) + \delta\}\cap B_{s_0}(x_0)|}{|B_{s_0}|} \geq \frac{5}{16} \]
    Let $s_1 < s_0$ be such that if $|x - x_0| < s_1$, then $|T(x) - T(x_0)| < \delta$. Then for $x_1\in B_{s_1}(x_0)$, we compute
    \begin{align*}
        \frac{|\{x : T(x) > T(x_1)\}\cap B_{s_0 + |x_1 - x_0|}(x_1)|}{|B_{s_0 + |x_1 - x_0|}|}
        &\geq \frac{|\{x : T(x) > T(x_0) + \delta\}\cap B_{s_0}(x_0)|}{|B_{s_0 + |x_1 - x_0|}|}
        \\&\geq \left(\frac{5}{16}\right)\left(1 - \frac{s_0}{s_0 + s_1}\right)^d
    \end{align*}
    Thus, if we take $s_1$ sufficiently small, this last quantity is greater than $\frac{1}{8}$, and all points in $B_{s_1}(x_0)$ are regular.    
\end{proof}

\subsection{Regular points}

We now turn toward understanding the behavior of the interface near regular points. Standard obstacle problem theory (\cite{caffarelli98}, \cite{blank}) gives that for $C^{0,\alpha}$ source, the interface is locally $C^{1,\alpha}$ at regular points, with the scale at which the regularity is achieved depending on the scale at which the zero set is sufficiently large. We discuss the dependence of this regularity in greater detail in the appendix. In particular, for this problem we have:

\begin{prop}\label{prop:regular_point_standard_boundary_regularity}
    $R$ can be covered by open neighborhoods $V$, each with the property that there exist constants $C, r > 0$ such that for each $x\in V$, $B_r(x)\cap\Omega_{T(x)}$ is the intersection of $B_r(x)$ with the lower graph of a $C^{1,\alpha}$ function in some coordinate system (depending on $x$) with seminorm bounded by $C$.
\end{prop}
\begin{proof}
    By Lemma \ref{lem:regular_point_boundary_estimate}, we need only show that if the zero set reaches density sufficiently close to $\frac{1}{2}$ at scale $r$ near $x$, then it does so at the same scale at all points near $x$. This is essentially immediate from the fact that $\Omega_t$ expands monotonically in $t$, with the measure $|\Omega_t|$ Lipschitz as a function of $t$ by Lemma \ref{lem:rho_regularity}.
\end{proof}

The dependence of the coordinate system in Proposition \ref{prop:regular_point_standard_boundary_regularity} is only a minor inconvenience, and we will eventually remove it in Proposition \ref{prop:regular_final_boundary_regularity}. To better understand this dependence, we introduce $\nu(x)$, defined for $x\in R$ as the outward unit normal to $\Omega_{T(x)}$ at $x$. We have spatial regularity of $\nu$ from the obstacle problem; namely, $\nu\in C^{0,\alpha}_{\loc}(R_t)$ for each $t$. Our goal will be to improve this to regularity of $\nu$ on $R$.

The key ingredients will be the regularity of $\Omega_t$ near regular points, and the strictly monotonic expansion of the $\Omega_t$. The essential idea will be that if the tangent planes to $\Omega_{T(x)}$ at $x$ and to $\Omega_{T(y)}$ at $y$ intersect for some points $x,y$ with different hitting times, they must intersect well away from $x$ and $y$ or else we will be able to use the regularity of the interfaces to show that $\partial \Omega_{T(x)}$ and $\partial \Omega_{T(y)}$ intersect, which contradicts monotonicity. This then gives control over the angle at which the tangent planes may intersect in terms of the distance between $x$ and $y$.

\begin{prop}\label{prop:unit_normal_holder}
    The outward unit normal vector $\nu$ to $\Omega_{T(x)}$ at $x$ satisfies $\nu\in C^{0,\alpha/(1 + \alpha)}_{\loc}(R)$.
\end{prop}
\begin{proof}

By Proposition \ref{prop:regular_point_standard_boundary_regularity}, we may cover $R$ with neighborhoods $V$ such that for each $x\in V$ uniformly, $\partial\Omega_{T(x)}\cap V$ is the lower graph in some coordinate system of a function $f_{T(x)}$, with the $f_t$ uniformly bounded in $C^{1,\alpha}$. We will restrict to such a $V$ for the remainder of the proof.

Then as a preliminary step, we can observe continuity of $\nu$ from the regularity and monotonicity of the interface by a purely geometrical argument. Namely, a $C^{1,\alpha}$ domain entertains a uniform interior and exterior cone condition, where the angle of the cone improves toward $\pi$ as we allow its height to approach 0; specifically, the cone in $B_r$ can be taken with angle $2\arccos(C r^\alpha)$, when the $C^{1,\alpha}$ seminorm is $C$. Thus, if $\nu$ were discontinuous at some $x\in R$, we could use compactness to find a sequence $(y_n)$ converging to $x$ with $T(y_n)$ either increasing or decreasing to $T(x)$ and $\nu(y_n)$ converging to some unit vector distinct from $\nu(x)$. Then for $n$ sufficiently large, at a sufficiently small scale, the interior cone at $x$ will intersect with the exterior cone of a $y_n$, or vice versa, and we draw a contradiction with the monotonic expansion of the $\Omega_t$ depending on whether the $T(y_n)$ are decreasing or increasing. 

\medskip

Then, we have checked that $x\to \nu(x)$ is continuous. Now, to obtain a quantitative local continuity estimate in view of the cone regularity we described above, we may restrict attention to $x,y\in V$ with $\frac{1}{2} < \nu(x)\cdot \nu(y) < 1$. Moreover, since the case $T(x) = T(y)$ is managed by the spatial regularity of the interface, we may assume that $T(x) > T(y)$. For notation, we let $r = \nu(x) \cdot \nu(y)$ and use $P_x, P_y$ to refer to the tangent planes to $\Omega_{T(x)}$ at $x$ and $\Omega_{T(y)}$ at $y$ respectively.

Let $v := \frac{\nu(y) - |r|\nu(x)}{\sqrt{1 - r^2}}$ be the projection of $\nu(y)$ into $\nu(x)^\perp$ , scaled to unit norm. Considering the point $x - hv$ for $h > 0$, we compute that its $\nu(y)$ component is $x\cdot \nu(y) - h\sqrt{1 - r^2}$. Thus $x - hv$ reaches $P_y$ precisely when $h = \frac{(x - y)\cdot \nu(y)}{\sqrt{1 - r^2}}$, and in general we have

\begin{equation}\label{eq:dist_Py}
  d(x-hv, P_y) \geq h\sqrt{1 - r^2} - \delta  \quad \hbox{ where } \delta:=|x-y|  
\end{equation}

Now, we apply the $C^{1,\alpha}$ regularity of the interface in $V$. For all $h$ sufficiently small, this regularity implies that $\partial \Omega_{T(x)}$ in $B_h(x)$ is contained in a $Ch^{1 + \alpha}$-neighborhood of $P_x$; in other words, $\Omega_{T(x)}$ and its exterior contain the following halfspaces:
\begin{equation}\label{eq:px_lower_halfspace}
   \{ z\in B_h(x) : (z-x-Ch^{1 + \alpha})\cdot \nu(x)\leq 0\} \subset \Omega_{T(x)}\cap B_h(x) 
\end{equation}
\begin{equation}\label{eq:px_upper_halfspace}
    \{ z\in B_h(x) : (z-x+Ch^{1 + \alpha})\cdot \nu(x)\geq 0\} \subset  (\R^d\setminus\Omega_{T(x)})\cap B_h(x)
\end{equation}
In particular, since $x - hv\in P_x$, it follows that there is a point $\tilde{x}\in \partial \Omega_{T(x)}$ with
\begin{equation}\label{eq:dist_x_tilde}
   |\tilde{x} - (x - hv)| < Ch^{1 + \alpha}. 
\end{equation}
Let $y_1$ be the nearest point in $P_y$ to $\tilde{x}$. We illustrate this with the figure below.

\begin{figure}[h]
    \centering
    \includegraphics{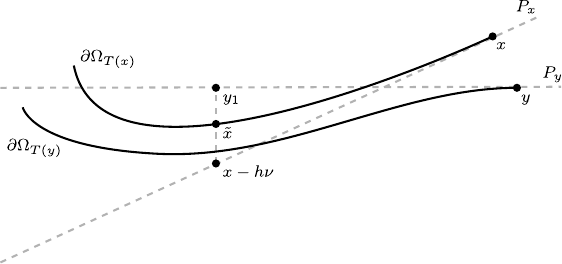}
    \caption{}
    \label{fig:unit_normal_proof}
\end{figure}

As the figure may suggest, $\tilde{x}$ cannot be too far below $y_1$, or else it falls into $\Omega_{T(y)}$, contradicting that $T(x) > T(y)$. Specifically, as in (\ref{eq:px_lower_halfspace}), we can apply the $C^{1,\alpha}$ regularity to to get that for all sufficiently small $r$, $\Omega_{T(y)} \cap B_r(y)$ contains the halfspace $\{ z : (z - y - Cr^{1 + \alpha}\nu(y))\cdot \nu(y) \leq 0 \}\cap B_r(y)$. Since $\tilde{x}\notin \overline{\Omega_{T(y)}}$, $\tilde{x}$ must not be contained in that halfspace, and we have

\begin{equation}\label{eq:halfspace_ineq_x_tilde}
    (\tilde{x} - y_1)\cdot \nu(y) > -C |y_1 - y|^{1 + \alpha}
\end{equation}

The left side here is the signed distance of $\tilde{x}$ to $P_y$. From (\ref{eq:dist_Py}), the signed distance of $x - hv$ to $P_y$ is bounded above by $-(h\sqrt{1 - r^2} - \delta)$, so using (\ref{eq:dist_x_tilde}), we conclude that the left side above is bounded above by $Ch^{1 + \alpha} - (h\sqrt{1 - r^2} - \delta)$. On the other hand, we have
    \[ |y_1 - y| \leq |\tilde{x} - y| \leq |\tilde{x} - (x - hv)| + |(x - hv) - x| + |x - y| \leq C_0 h^{1 + \alpha} + h + \delta \]
When $\delta < h < 1$, this is $O(h)$, and so $|y_1 - y|^{1 + \alpha} \geq -C h^{1 + \alpha}$ for some $C$. Thus, the inequality (\ref{eq:halfspace_ineq_x_tilde}) becomes
    \[ Ch^{1 + \alpha} - (h\sqrt{1 - r^2} - \delta) > - C h^{1 + \alpha}  \]
Rearranging and absorbing constants, this means
    \[ h\sqrt{1 - r^2} - \delta \leq Ch^{1 + \alpha} \]
so that
    \[ |\nu(x) - \nu(y)| \leq \sqrt{1 - r^2} \leq Ch^{\alpha} + \delta h^{-1} \]
Optimizing $h$, we get
    \[ |\nu(x) - \nu(y)| \leq C\delta^{\frac{\alpha}{1 + \alpha}} \]
where the constant depends only on the uniform bound for the $C^{1,\alpha}$ seminorms of the graphs, and on $\alpha$. Hence we conclude.
\end{proof}

As a result of the regularity of $\nu$, we can improve Proposition \ref{prop:regular_point_standard_boundary_regularity} to also have the coordinate system chosen locally uniformly. In other words, near regular points, one can fix a local coordinate system in which the free boundary evolves as a $C^{1,\alpha}$ graph over some time interval.

We now turn toward applying the improved geometry of the patch at regular points to the pressure. Elliptic regularity for $C^{1,\alpha}$ domains implies that the pressure $p(\cdot, t)$ has a well-defined gradient on $R_t$, and the Hopf lemma for $C^{1,\alpha}$ domains implies that $\nabla p(\cdot, t)$ is nonvanishing on $R_t$. However, there is an important limitation here: $x\mapsto \nabla p(x, T(x))$ is not necessarily continuous on $R$, complicating our analysis. This is illustrated with the following example:

\begin{remark}\label{rem:regular_points_nondifferentiable}
    Singular points can exert a nonlocal effect on the pressure gradient.

    For example, if we consider a pressure supported on a strip of width $h$ with zero boundary conditions and constant Laplacian $-1$, then we can see that $|\nabla p| = \frac{h}{2}$ on the boundary, since the solution to the one-dimensional problem with $p(0) = p(h) = 0$ is $p(x) = \frac{1}{2}x(h-x)$. In particular, it follows that if we have a patch which consists of two strips, and those strips merge along a hyperplane at some time, then $\nabla p$ has a jump discontinuity in time at every regular point at the time those strips merge.

    Similar examples can be considered for singular points in each stratum $\Sigma^k$ by examining a cylindrical patch with a cylindrical hole as the radius of the cylinder shrinks to 0. Here, by cylinder we mean the product of $\R^k$ and a $(d-k)$-dimensional ball, for $0\leq k \leq d-1$.
\end{remark}

We note that the discontinuity in the gradient in Remark \ref{rem:regular_points_nondifferentiable} is a jump in magnitude, not in direction. Indeed, the zero boundary condition implies that $\frac{\nabla p(x, T(x))}{|\nabla p(x, T(x))|} = -\nu(x)$ on $R_t$, so the regularity of $\nu$ from Proposition \ref{prop:unit_normal_holder} rules out such discontinuities.  

This example also proves to be an obstacle to higher regularity of $T$; we will later see in Lemma \ref{lem:hitting_time_differentiable} that the derivative of $T$ closely depends on $\nabla p(x, T(x))$ when it exists. As a result, we can hope for $T$ to be at best Lipschitz on $R$.

The key in establishing this regularity will be to obtain a quantitative estimate from the Hopf lemma, to get a locally uniform lower bound for $|\nabla p(x, T(x))|$ on $R$. For this, we require an a priori estimate for the growth of the solution, in addition to control over the geometry. We can obtain this growth estimate from the strict superharmonicity of $p(\cdot, t)$, and so we have the following statement:

\begin{lemma}\label{lem:quantitative_hopf_lemma}
    Let $r_0, c_0, C_0 > 0$ and $\alpha\in (0,1)$. Suppose $u$ is a positive $C^2$ solution to $\Delta u \leq -c_0 < 0$ on $B_{r_0}(0)\cap \{ x : x_d > C_0|x'|^{1 + \alpha} \}$ with $u(0) = 0$, where we write $x = (x', x_d)$. Then there exist $\varepsilon, \delta > 0$ depending only on $d, \alpha, r_0, c_0, C_0$ such that for all $h\in (0,\delta)$, we have $u(h e_n) \geq \varepsilon h$.
\end{lemma}
\begin{proof}
First, we let $C_1 = C_0 + 1$, to give additional separation from the boundary when we are away from $0$, and let $U_r = B_r(0)\cap \{ x : x_d > C_1 |x'|^{1 + \alpha}\}$ for $0 < r < r_0$. For a given $r$, we will decompose $\partial U_r$ as $\Gamma_1 \cup \Gamma_2$, where $\Gamma_2 := \partial B_r(0)\cap \{ x : x_d \geq C_1|x'|^{1 + \alpha} \}$ is a spherical cap, and $\Gamma_1 = B_r(0)\cap \{ x : x_d = C_1 |x'|^{1 + \alpha} \}$.

\begin{figure}[h]
    \centering
    \includegraphics[scale=0.4]{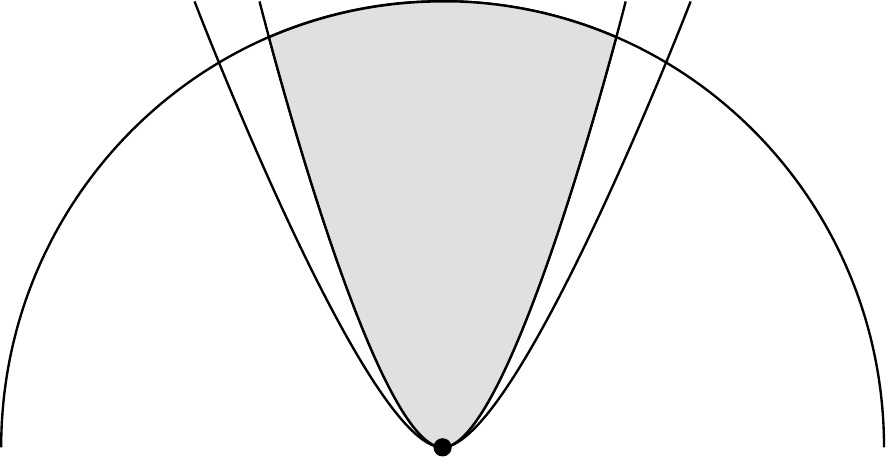}
    \caption{the region $U_r$ (shaded)}
    \label{fig:hopf}
\end{figure}

The proof of the Hopf lemma proceeds by perturbing $u$ by a function $v$ constructed specially for the domain and applying the comparison to the result. In particular, if for some $v, \varepsilon > 0, \delta > 0$ we have
\[ \begin{cases}
    \Delta(u - \varepsilon v) \leq 0 \hbox{ on } U_\delta \\
    u - \varepsilon v \geq 0 \hbox{ on } \partial U_\delta \\
    \partial_d v(0) = 1
\end{cases} \]
then the comparison principle implies that $u - \varepsilon v \geq 0$ on $U_\delta$, and the control over the derivative of $v$ at 0 implies the result.

Borrowing from the proof of the Hopf lemma for $C^{1,\alpha}$ domains in \cite{li2007hopf}, we let 
\[ v(x) = x_d + \frac{2C_1}{\alpha}(\alpha + d - 1)x_d^{1 + \alpha} - 2C_1|x|^{1 + \alpha} \]
We note by direct computation that $\partial_d v(0) = 1$ and
\[ \Delta v(x) = 2C_1(1 + \alpha)(\alpha + d - 1)\left(x_d^{\alpha - 1} - |x|^{\alpha - 1}\right)\]
In particular, $\Delta v(x) \geq 0$ on $U_r$ for any $r$.

Thus, we reduce to verifying the boundary condition, which requires using the particular behavior of $v$ on $\Gamma_1$, and choosing $\delta$ and $\varepsilon$ appropriately to control $v$ on $\Gamma_1$ and $\Gamma_2$. Following this plan, we first check the boundary inequality on $\Gamma_1$, defined above, where we have $x_d = C_1 |x'|^{1 + \alpha}$. Along a curve $x(t) = (te', C_0 t^{1 + \alpha})$ for a unit vector $e'\in \R^{d-1}$, we have
\[ v(x(t)) = C_1 t^{1 + \alpha} + \frac{2C_0^2}{\alpha}(\alpha + d - 1)t^{(1 + \alpha)^2} - 2C_0(t^2 + C_1 t^{2(1 + \alpha)})^{(1 + \alpha)/2} \]
For small $t > 0$, we can drop the higher order terms and see that $v$ grows like $C_1 t^{1 + \alpha} - 2C_1 t^{1 + \alpha} = -C_1 t^{1 + \alpha}$. In particular, there exists $\delta_0 = \delta_0(d, \alpha, C_0)$ such that if $\delta \leq \delta_0$, then $v(x) \leq 0$ on $B_\delta(0) \cap\{ x : x_d = C_1 |x'|^{1 + \alpha} \}$. We will set $\delta = \min(\delta_0, r_0/2)$ for the rest of the proof.

Next, we handle the boundary inequality on $\Gamma_2$, the spherical cap defined by $\Gamma_2 = \partial B_\delta(0)\cap \{ x : x_d \geq C_1 |x'|^{1 + \alpha} \}$. By compactness, $\Gamma_2$ has positive distance to $\{ x : x_d \geq C_0 |x'|^{1 + \alpha} \}$. Let $D = D(d, \alpha, \delta, C_0)$ denote this distance. By shrinking $D$ to $\frac{r_0}{2}$ if necessary, we get that $u$ is defined on a ball of radius $D$ at each point of $\Gamma_2$. The paraboloid on the ball of radius $D$ with zero boundary data and Laplacian $-c_0$ is a subsolution to $u$, so we conclude that $u(x) \geq \frac{c_0}{2}D^2$ on $\Gamma_2$. Let $M = \max(1, \max_{\Gamma_2} v)$, and then we can take
\[ \varepsilon = \frac{c_0 D^2}{2M} \]
to get $u - \varepsilon v \geq 0$ on $\Gamma_2$, with $\varepsilon$ depending on all of the parameters in the statement of the lemma. 

On $\Gamma_1$, we have $u\geq 0$ and $v\leq 0$, so $u - \varepsilon v \geq 0$. On $\Gamma_2$, we chose $\varepsilon$ so that $u - \varepsilon v \geq 0$. Thus, $u - \varepsilon v \geq 0$ on $\partial U_\delta$, so we conclude.

\end{proof}

With the Hopf lemma estimate and the regularity of the boundary near regular points, we can conclude that $p(\cdot, T(x))$ has linear growth at $x$, locally uniformly on $R$.

\begin{prop}\label{prop:regular_linear_growth}
    $R$ can be covered with neighborhoods $V$ with the following property: there exist parameters $C, c, r_0 > 0$ such that for any $x\in V$,
    \[ c \leq r^{-1}\sup_{B_r(x)} p(\cdot, T(x)) \leq C \]
    In particular, $|\nabla p(x, T(x))| \sim 1$ on $V$, for implicit constants depending on $V$.
\end{prop}
\begin{proof}
    To obtain the lower bound, we apply Lemma \ref{lem:quantitative_hopf_lemma}. Thus, we must show that the parameters $r_0, c_0, C_0$ from the statement of the lemma can be chosen locally uniformly on $R$. By Lemma \ref{lem:nutrient_lower_bound}, we can choose $c_0 > 0$ locally uniformly in time so that $\Delta p = -n \leq -c_0 < 0$, using the assumption that $n_0$ is bounded away from 0. By Proposition \ref{prop:regular_point_standard_boundary_regularity}, for any $x_0\in R$, we can find a neighborhood $V$ of $x_0$ in which we have a uniform $C_0$ so that each $R_t$ which intersects the neighborhood does so as a graph with $C^{1,\alpha}$ seminorm controlled by $C_0$. By taking $r_0$ so that $B_{2r_0}(x_0)\subset V$, we can use $r_0, C_0$ as the parameters for all points in $B_{r_0}(x_0)$.
    
    To obtain the upper bound, we first use Proposition \ref{prop:regular_point_standard_boundary_regularity}. This gives us a neighborhood $V$ and a parameter $r_1 > 0$ for which $B_{r_1}(x)\cap \partial \Omega_{T(x)}$ is a $C^{1,\alpha}$ graph, with uniform control over the $C^{1,\alpha}$ seminorm. In particular, this regularity implies that there is an $r_2 > 0$ such that for all $x\in V$ and all $r\in (0, r_2]$, we have $x - r\nu(x)\in \Omega_{T(x)}$. Then the mean value theorem along the path $x - r\nu(x)$ implies that there exists some $r$ for which
    \[ \nabla p(x, x - r\nu(x))\cdot -\nu(x) = \frac{p(x - r_2 \nu(x), T(x))}{r_2} \]
    We recall from the proof of Lemma \ref{lem:w_regularity} that $p\in L^\infty(\R^d\times [0, \tau])$, for any $\tau\in (0,\infty)$, using that the patch has bounded support and comparing to a sufficiently large paraboloid supersolution. Thus, we have
    \[ |\nabla p(x, x - r\nu(x))\cdot \nu(x)| \leq Cr_2^{-1} \]
    Using the regularity of the boundary, the $L^\infty$ bound on the pressure, and the $L^\infty$ bound on the nutrient from Lemma \ref{lem:n_regularity}, we can invoke boundary Schauder estimates to have $p(\cdot, T(x))$ uniformly $C^{1,\alpha}$ on $B_{r_1}(x)\cap \Omega_{T(x)}$ for each $x\in V$. Thus, we can transfer our bound to the boundary:
    \[ |\nabla p(x, T(x))| = |\nabla p(x, T(x))\cdot \nu(x)| \leq |(\nabla p(x, T(x)) - \nabla p(x - r\nu(x)))\cdot\nu(x)| + |\nabla p(x - r\nu(x))\cdot \nu(x)| \leq Cr + Cr_2^{-1} \]
    This also gives the upper bound on the linear growth of $p(\cdot, T(x))$ near $x$, so we conclude.
\end{proof}

Having established nondegeneracy of the pressure, we get improved regularity of $T$ on $R$.

\begin{cor}\label{cor:regular_points_t_lipschitz}
    $T\in C^{0,1}_{\loc}(R)$. In other words, $T$ attains its optimal regularity on $R$ in light of Remark \ref{rem:regular_points_nondifferentiable}, barring additional assumptions on $\Sigma$.
\end{cor}
\begin{proof}
    This follows from the boundary regularity from Proposition \ref{prop:regular_point_standard_boundary_regularity} and the linear nondegeneracy of the pressure from Proposition \ref{prop:regular_linear_growth}, via a method similar to the comparison arguments of Section 4. Using these properties, we can construct a radial subsolution initially supported on an annulus in $\Omega_t$ which expands at a constant rate, near any $x_0\in R$. From this argument, we get
    \begin{equation}\label{eq:t_one_sided_lip}
        (T(x) - T(x_0))_+ \leq C|x - x_0|
    \end{equation}
    for $x_0\in R$ and $x$ sufficiently close to $x_0$. Since the boundary regularity and linear growth rate are uniform for $x_0$ restricted to a compact $K\subset R$, the one-sided bound \eqref{eq:t_one_sided_lip} holds uniformly for $x, x_0\in K$, and so we conclude that $T$ is locally Lipschitz on $R$.
\end{proof}

Finally, we turn to refining the statement of Proposition \ref{prop:regular_point_standard_boundary_regularity}. The first step will be to use the linear growth of the pressure gradient to control the boundary in the Hausdorff metric.

\begin{lemma}\label{lem:regular_hausdorff_distance}
    For any $x\in R$, there exists parameters $r, \delta > 0$ such that for $t_1, t_2\in (T(x) - \delta, T(x) + \delta)$, we have 
    \[ \sup_{y_1\in R_{t_1}\cap B_r(x)} \inf_{y_2\in R_{t_2}\cap B_r(x)} |y_1 - y_2| \sim |t_1 - t_2| \]
    In other words, $D(R_{t_1}\cap B_r(x), R_{t_2}\cap B_r(x))\sim |t_1 - t_2|$, where $D$ denotes Hausdorff distance. Here, all parameters and implicit constants depend on $x$.
\end{lemma}
\begin{proof}
    Let $B$ be a ball compactly contained in $R$, and let $y_1, y_2\in B$ with $T(y_1) < T(y_2)$. Since $T\in C^{0,1}_{\loc}(R)$ by Corollary \ref{cor:regular_points_t_lipschitz}, we have $|T(y_1) - T(y_2)| \leq C(B)|y_1 - y_2|$. We get the reverse bound by an analogous argument using the boundary regularity of Proposition \ref{prop:regular_point_standard_boundary_regularity} and using Proposition \ref{prop:regular_linear_growth}'s upper bound on the linear growth of $p(\cdot, t_1)$ away from $y_1$. That is, following the approach of Section 4, we can construct a supersolution supported outside a ball in the exterior of $\Omega_{T(y_1)}$ near $y_1$, such that the supersolution expands at a constant rate. From that argument, we get that there exists a point $\tilde{y_2}\in R_{t_2}$ with $|y_1 - \tilde{y_2}| \leq C(B)|T(y_1) - T(y_2)|$. Letting $t_1 = T(y_1), t_2 = T(y_2)$, and restricting to the case where $|t_1 - t_2|$ is sufficiently small to guarantee that the $\tilde{y_2}$ from before is in $B$, we have shown that $D(R_{t_1}\cap B, R_{t_2}\cap B) \sim |t_1 - t_2|$ for implicit constants depending on $B$, from which we can obtain the original statement.
\end{proof}

Using the previous result, we can now state our final improved form of Proposition \ref{prop:regular_point_standard_boundary_regularity}.

\begin{prop}\label{prop:regular_final_boundary_regularity}
    $R$ is covered by neighborhoods $V$ with the following property: there exists $r > 0$, a coordinate system $(x', x_n)$, and a locally defined function $f(x', t)$ such that for each $x\in V$, $\Omega_{T(x)}\cap B_r(x)$ is the lower graph $\{ y\in B_r(x) : y_n \leq f(y', T(x)) \}$. Moreover, $f$ is uniformly $C^{1,1}$ in space and $C^{0, 1}$ in time.
\end{prop}
\begin{proof}
   First, by Proposition \ref{prop:unit_normal_holder}, we note that the coordinate system in Proposition \ref{prop:regular_point_standard_boundary_regularity} can be chosen locally uniformly, so that we get $r > 0$ and the family $f(x', t)$ which are uniformly $C^{1,1}$ in space by the Lipschitz regularity of $T$. 
   
   Thus, it remains only to check that the regularity in time follows from our control over the Hausdorff distance. Fix $x'$ and $t_1, t_2$, and write $x_1 = (x', f(x', t_1))$, $x_2 = (x', f(x', t_2))$. Then for $\varepsilon = |f(x', t_1) - f(x', t_2)| = |x_1 - x_2|$, $B_\varepsilon(x_1)$ contains the point $\tilde{x}\in \partial \Omega_{T(x_2)}$ which minimizes the distance to $x_1$. By Lemma \ref{lem:regular_hausdorff_distance}, after possibly shrinking our neighborhood, $|x_1 - \tilde{x}| \leq C|t_1 - t_2|$. In particular, $|\tilde{x}' - x'| \leq C|t_1 - t_2|$, so $|f(\tilde{x}', t_2) - f(x', t_2)| \leq C|t_1 - t_2|$, for some larger $C$ given by the $C^1$ spatial regularity of $f$. Then
   \[ |f(x', t_1) - f(x', t_2)| \leq |f(x', t_1) - f(\tilde{x}', t_2)| + |f(\tilde{x}', t_2) - f(x', t_2)| \leq C|t_1 - t_2|   \]
   which completes the proof.

\end{proof}

\subsection{Singular points}

Now we proceed to analysis of the singular set, with the goal of controlling singular points in dimension. First, we will show that the blowup profile at singular points varies continuously along the spacetime interface. The main tool will be a uniform approximation result that we prove in the appendix, which will allow us to make use of the uniform-in-time spatial continuity of $\eta$.

\begin{prop}\label{prop:singular_hessian_continuous}
    $x\mapsto D^2w(x, T(x))$ is continuous on $\Sigma$.
\end{prop}
\begin{proof}
    First, by Lemma \ref{lem:w_regularity} and Lemma \ref{lem:eta_regularity}, we have that $w$ is locally spacetime Lipschitz and $\eta$ is $C^{0,\alpha}$ in space locally uniformly in time. We also have that $T$ is locally $C^{0,\alpha}$ on $\mathcal{O}$. Thus, for some $C > 0$, we can restrict to a spacetime neighborhood of the interface where all of these norms are bounded by $C$ ( in the case of $T$, in the sense of the neighborhood's projection into space).
    
    Let $\varepsilon > 0$. By Lemma \ref{lem:uniform_sing_blowup_scale}, there exists a scale depending only on the modulus of continuity of $1 - \eta$ in (\ref{eq:w_obst_eqn}), such that the quadratic blowup uniformly approximates $w$ near singular points at that scale; concretely, there is a $\delta = \delta(\varepsilon / 3, C) > 0$ such that if $x_0$ is a singular point in our neighborhood with blowup $q_0(x) = \frac{1}{2}x\cdot D^2w(x_0)x$ (recentered at 0), then
    \[ \|\delta^{-2}w(x_0 + \delta x, T(x_0)) - q_0\|_{C^1(B_1)} < \frac{\varepsilon}{3} \]
    In particular, if $x_1$ is another singular point in the same neighborhood with blowup $q_1$, then we have
    
    \begin{equation}
        \|q_0 - q_1\|_{L^\infty(B_1)} \leq \frac{2\varepsilon}{3} + \delta^{-2}\|w(x_0 + \delta x, T(x_0)) - w(x_1 + \delta x, T(x_1))\|_{L^\infty(B_1)}
    \end{equation}
    
    From the regularity of $T$ and $w$ on our chosen neighborhood, it follows that
    
    \begin{equation}
        \|q_0 - q_1\|_{L^\infty(B_1)} \leq \frac{2\varepsilon}{3} + C_1\delta^{-2}|x_1 - x_0|^\alpha
    \end{equation}
    for some $C_1$ depending only on $C$. Then it is clear that for $|x_1 - x_0|$ sufficiently small, we have $\|q_0 - q_1\|_{L^\infty(B_1)} < \varepsilon$. By equivalence of norms on $\R^{d\times d}$, this gives $|D^2 w(x_0) - D^2 w(x_1)| \leq O(\varepsilon)$, and we conclude.
\end{proof}

We remark that $D^2 w(x, T(x))$ also exists in a one-sided sense for $x\in R$, with $D^2 w(x, T(x)) = \nu(x)\nu(x)^T$ for $\nu$ as in Proposition \ref{prop:unit_normal_holder}. Thus, in light of that proposition, $D^2 w(x, T(x))$ is continuous on $R$. However, due to the jump in rank, there is no possibility of continuity from $R$ to $\Sigma^{k}$ when $k < d-1$.

Using the continuous dependence of the blowup, we can subsequently apply a Whitney extension argument to obtain that singular points are contained in $C^1$ manifolds. Following the approach in \cite{petrosyan}, we introduce the following lemma:

\begin{lemma}[Whitney's extension theorem]\label{lem:whitney}
    Let $K\subset \R^d$ be compact, and suppose we have a function $f:K\to \R$ and a family of degree $m$ polynomials $p_x$ indexed over $K$. If
    \begin{enumerate}[(i)]
        \item $p_{x_0}(x_0) = f(x_0)$ for each $x_0\in K$
        \item $|D^k p_{x_0}(x_1) - D^kp_{x_1}(x_1)| = o(|x_0 - x_1|^{m-k})$ for $x_0, x_1\in K$ and $0\leq k \leq m$.
    \end{enumerate}
    Then $f$ extends to a $C^m$ function on $\R^d$ such that $f(x) = p_{x_0}(x) + o(|x - x_0|^m)$ for all $x_0\in K$.
\end{lemma}

\begin{prop}\label{prop:singular_manifold}
    Near a point in $\Sigma^k$, $\Sigma$ is locally contained in a $C^1$ manifold of dimension $k$. In particular, $\Sigma$ is contained in countably many $C^1$ submanifolds of dimension $d-1$.
\end{prop}

\begin{proof}
We fix a compact $K\subset \Sigma$ for the proof. We will apply Lemma \ref{lem:whitney} to extend the zero function on $K$, with second order Taylor polynomial $q_{x_0}$ at $x_0\in K$ given by the quadratic blowup of $w$ at each point. Specifically, this results in 
\[ q_{x_0}(x) = \frac{1}{2}(x-x_0)\cdot D^2 w(x_0, T(x_0))(x-x_0) \hbox{ where } x_0\in K, x\in \R^d \]
The extension will give us a $C^2$ function $f$ on $\R^d$ such that $\nabla f \equiv 0$ on $K$, after which the implicit function theorem will imply that in a neighborhood of $x_0\in K$, the set $\{ \nabla f = 0 \}$ is contained in a $C^1$ manifold of dimension $\dim\ker D^2 f(x_0) = \dim\ker D^2 w(x_0, T(x_0))$.

Thus, we proceed to verify the assumptions of the lemma. First, we note that by Lemma \ref{lem:eta_regularity}, $\eta$ is uniformly $C^{0,\alpha}$ in space in a spacetime neighborhood of the interface as it passes through $K$ in space. It follows by our uniform approximation result for the quadratic blowup at singular points, Lemma \ref{lem:uniform_sing_blowup_scale}, that there is a modulus of continuity $\sigma$ such that
\begin{equation}\label{eq:singular_blowup_convergence}
    \|r^{-2}w(x_0 + r(x - x_0), T(x_0)) - \frac{1}{2}(x - x_0)\cdot D^2w(x_0, T(x_0)) (x - x_0)\|_{C^1_x(B_1)} \leq \sigma(r)
\end{equation}
for any $x_0\in K$. Then, if we apply this estimate to $x_0, x_1\in K$ with $T(x_0) \leq T(x_1)$, we get
\[ |q_{x_0}(x_1) - q_{x_1}(x_1)| = |\frac{1}{2}(x_1 - x_0)\cdot D^2w(x_0, T(x_0))(x_1 - x_0)| \leq |x_0 - x_1|^2\sigma(|x_0 - x_1|) \]
directly from (\ref{eq:singular_blowup_convergence}) and the fact that $w(x_1, T(x_0)) = 0$. On the other hand, if $T(x_0) > T(x_1)$, we have that $|q_{x_1}(x_0)| \leq |x_0 - x_1|^2 \sigma(|x_0 - x_1|)$ from the above, and
\[ |q_{x_1}(x_0) - q_{x_0}(x_1)| = |\frac{1}{2}(x_1 - x_0)\cdot(D^2w(x_1,T(x_1)) - D^2 w(x_0, T(x_0)))(x_1 - x_0)| \leq o(|x_1 - x_0|^2)  \]
by the continuity of the Hessian from Proposition \ref{prop:singular_hessian_continuous}. This completes the $k = 0$ case of (ii) in Lemma \ref{lem:whitney}.

The verification of the $k = 1$ case of (ii) is similar, using the derivative bound from (\ref{eq:singular_blowup_convergence}). In general, we have $\nabla q_{x_0}(x) = D^2 w(x_0, T(x_0))(x - x_0)$. Then, when $T(x_0) \leq T(x_1)$, we get
\[ |\nabla q_{x_0}(x_1) - \nabla q_{x_1}(x_1)| = |D^2 w(x_0, T(x_0))(x_1 - x_0)| \leq |x_1 - x_0|\sigma(|x_1 - x_0|) \]
directly from (\ref{eq:singular_blowup_convergence}) and the fact that $\nabla w(x_1, T(x_0)) = 0$. On the other hand, if $T(x_0) > T(x_1)$, we have that 
\[ |\nabla q_{x_1}(x_0) + \nabla q_{x_0}(x_1)| = |\frac{1}{2}(D^2w(x_0, T(x_0)) - D^2w(x_1, T(x_1)))(x_1 - x_0)| \leq o(|x_1 - x_0|) \]
again by Proposition \ref{prop:singular_hessian_continuous}, giving that 
\[ |\nabla q_{x_0}(x_1) - \nabla q_{x_1}(x_1)| = |\nabla q_{x_0}(x_1)| \leq |\nabla q_{x_0}(x_1) + \nabla q_{x_1}(x_0)| + |\nabla q_{x_1}(x_0)| \leq o(|x_1 - x_0|) \]

Finally, the $k = 2$ case of (ii) in Lemma \ref{lem:whitney} is exactly continuity of the Hessian, from Proposition \ref{prop:singular_hessian_continuous}. Thus, the conditions of the lemma are satisfied, which completes the proof.

\end{proof}

\begin{figure}[h]
    \centering
    \includegraphics[scale=0.2]{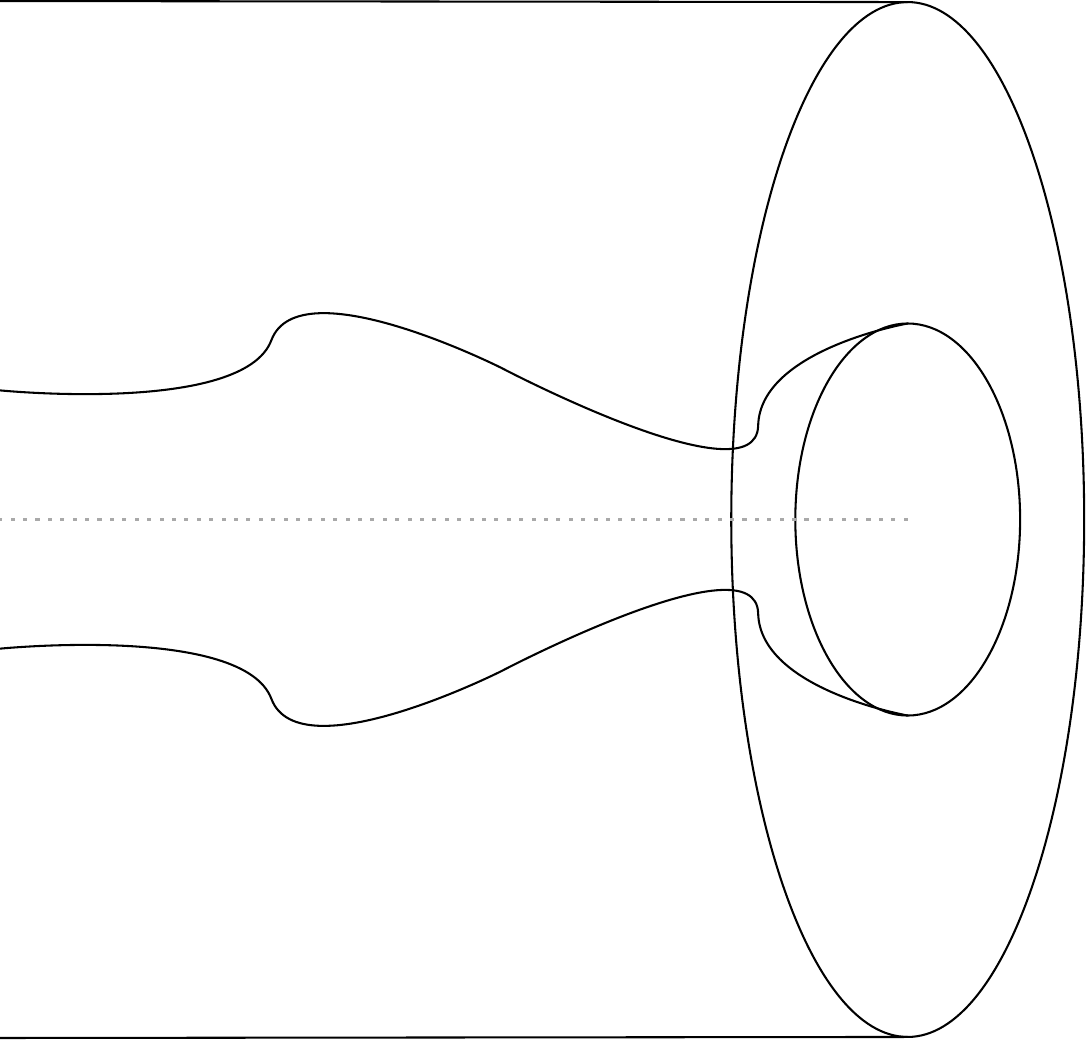}
    \caption{A cylindrical patch with a nearly cylindrical hole. As the hole contracts, singular points are expected to occur near the axis (dotted). Due to variation in the diameter of the hole, however, singular points may occur at different times. Proposition \ref{prop:singular_manifold} confirms that we have the expected spatial regularity for the set $\Sigma\subset \R^d$ of points which are singular at any time.}
    \label{fig:almost_cylinder}
\end{figure}

We stress that this result is for $\Sigma$, and not for the corresponding subset $\mathrm{Graph}_T(\Sigma)$ of the spacetime interface, where we use the notation introduced in \eqref{eq:graph}. Since $T$ is only known to be H\"{o}lder continuous on $\Sigma$, we obtain the weaker result that $\mathrm{Graph}_T(\Sigma)$ is locally contained in $C^{0,\alpha}$ manifolds of dimension $d-1$. Nevertheless, we are able to apply this to establish control in Hausdorff dimension over the interface. The Hausdorff dimension of the spatial interface has previously been studied in \cite{pqm} for a similar problem, where the obstacle problem satisfied by $w$ was applied to conclude that $\partial \Omega_t$ has locally finite $(d-1)$-dimensional Hausdorff measure for each $t$. Using the Lipschitz regularity of $T$ near regular points and our spatial control over singular points, we can study the spacetime interface $\mathrm{Graph}_T(\mathcal{O})$ for the first time and show that it has the expected Hausdorff dimension $d$. We summarize the consequences of Proposition \ref{prop:singular_manifold} with the following statements.

\begin{cor}\label{cor:hausdorff_measure}
    We have, for $\alpha$ as in Theorem \ref{thm:hitting_time_holder}:
    \begin{enumerate}[(i)]
        \item $\partial \Omega_t$ has finite $(d-1)$-dimensional Hausdorff measure.
        \item $\Sigma$ has locally finite $(d-1)$-dimensional Hausdorff measure. In particular, $\Sigma_t$ has zero $(d-1)$-dimensional Hausdorff measure for all but countably many $t$ in $(0,\infty)$, and for a.e. $t\in (0, \infty)$, $\Sigma_t$ has Hausdorff dimension at most $d - 1 - \alpha$.
        \item $\mathrm{Graph}_T(\mathcal{O})$ has Hausdorff dimension $d$, and decomposes as $\mathrm{Graph}_T(R)\cup \mathrm{Graph}_T(\Sigma)$, where the first set is relatively open with locally finite $d$-dimensional Hausdorff measure, and the second set has locally finite $(d - \alpha)$ Hausdorff measure.
    \end{enumerate}
\end{cor}
\begin{proof}
    We use the fact that $\rho\in L^\infty_t BV_x([0, \tau]; \R^d)$, and thus $\Omega_t$ is a set of finite perimeter. In particular, we can consider the reduced boundary $\partial^* \Omega_t$, which has finite $(d-1)$ Hausdorff measure and contains $R_t$, by the local regularity of the boundary at those points. On the other hand, $\Sigma_t$ is locally contained in a $C^1$ manifold of dimension $(d-1)$. In fact, since $\Sigma_t$ is compact, $\Sigma_t$ is contained in a bounded $C^1$ manifold, which then also has finite $(d-1)$ measure, so we have 
    \[ \mathcal{H}^{d-1}(\partial \Omega_t) \leq \mathcal{H}^{d-1}(\partial^* \Omega_t) + \mathcal{H}^{d-1}(\Sigma_t) < \infty \]

    Since $\Sigma$ is locally contained in $C^1$ manifolds of dimension $d-1$, the $(d-1)$-dimensional Hausdorff measure on $\Sigma$ is locally finite. In particular, it is also $\sigma$-finite, which implies that there cannot be uncountably many $t$ for which $\Sigma_t$ has positive $(d-1)$ measure. The improvement to $(d-1-\alpha)$ dimension at a.e. time follows from a geometric measure theory lemma of \cite{figalli_generic}. Since $\Sigma$ has dimension $d-1$ and $T\in C^{0,\alpha}_{\loc}(\mathcal{O})$, Corollary 7.8 of that paper directly gives the result.

    For the final statement, we use the general result that the graph of a $C^{0,\alpha}$ function on a set of Hausdorff dimension $s$, for $\alpha\in (0, 1]$ and $s \geq 0$, has Hausdorff dimension at most $s + 1 - \alpha$. In particular, since $T$ is locally Lipschitz on $R$, which is open in $\R^d$, and locally $C^{0,\alpha}$ on $\Sigma$, which is Hausdorff dimension at most $d-1$, we get the local control in Hausdorff measure for the graphs of those sets. Then we can write $\mathrm{Graph}_T(\mathcal{O})$ as a countable union of sets of Hausdorff dimension at most $d$, so we conclude.
\end{proof}

We remark that a natural open question is whether the space-time interface has locally finite $d$-dimensional Hausdorff measure near singular points. This would follow, for example, if $T$ were known to be uniformly Lipschitz near $\Sigma$.

\subsection{Speculative Results}

We finish our treatment of the singular set by noting that stronger generic control over $\Sigma_t$ is possible with slightly stronger regularity than currently known: namely, when $T$ is Lipschitz.

\begin{prop}\label{coarea_generic}
    If $T\in C^{0,1}_{\loc}(\mathcal{O})$, then $\Sigma_t$ has $(d-2)$-Hausdorff measure 0 for a.e. $t\in (0,\infty)$.
\end{prop}
\begin{proof}
Here we follow an argument by \cite{monneau}, originally applied to a hitting time for the constant Laplacian obstacle problem with a time-varying condition on the fixed boundary. Since $T$ is Lipschitz, Proposition 4.6 in \cite{monneau} implies that for any compact subset $K$ of $\Sigma$,
    \[ \limsup_{x,y\in K, |x-y|\to 0} \frac{|T(x) - T(y)|}{|x - y|} = 0 \]
    Then from the coarea formula, we have
    \[ \int_{\Sigma^d} |\nabla T_{|_{\Sigma^d}}|\,d\mathcal{H}^{d-1} = \int_0^\infty \mathcal{H}^{d-2}(T_{|_{\Sigma^{d-1}}}^{-1}(t))\,dt = \int_0^\infty \mathcal{H}^{d-2}(\Sigma^{d-1}_t)\,dt \]
    Then the integrand on the left is 0, so $\Sigma^{d-1}_t$ has $(d-2)$-Hausdorff measure 0 for a.e. $t$. On the other hand, $\Sigma^k_t$ has $(d-2)$-Hausdorff measure 0 for $k < d-2$ and all $t$, while $\Sigma^{d-2}_t$ has positive $(d-2)$-Hausdorff measure for at most countably many $t$, so the result follows.
\end{proof}

We finish our treatment of the regular set by investigating the regularity improvement possible under the assumption that no singular points occur at some time. As suggested by Remark \ref{rem:regular_points_nondifferentiable}, an assumption of this form is required to go beyond the regularity established in Proposition \ref{cor:regular_points_t_lipschitz}. The idea here will be to apply the regularity of $\nu$ from Proposition \ref{prop:unit_normal_holder} in conjunction with global Schauder estimates to prove time regularity of $p$ and higher spatial regularity of $T$.

As a preliminary step, we show the relationship between $\nabla p$ and $\nabla T$.

\begin{lemma}\label{lem:hitting_time_differentiable}
    Suppose that for some open $U\subset R$, we have that $\nabla p$ is continuous in spacetime on $(U\times (t_0, t_1))\cap \overline{\{ (x,t) : \rho(x,t) = 1 \}}$ for some $t_0, t_1$ with $\inf T(U) < t_0 < t_1 < \sup T(U)$. Then $T$ is continuously differentiable on $U\cap T^{-1}((t_0, t_1))$ with $\nabla T(x) = -\frac{\nabla p(T(x), x)}{|\nabla p(T(x), x)|^2}$.
\end{lemma}
\begin{proof}
Let $e$ be a vector with positive component in the inward normal direction to $\Omega_{T(x)}$ at $x$; that is, with $e\cdot \nu(x) < 0$. Then, we have
\[ \nabla w(T(x), x + he) = \sgn_+(T(x) - T(x + he))\int_{T(x + he)}^{T(x)} \nabla p(t, x + he)\,dt \]
If we divide both sides by $h$ and let $h\to 0$, then the left side converges to $(\nu(x)\cdot e)\nu(x)$, from the quadratic blowup. If the right side is nonzero, we rewrite it as
\[ \frac{T(x) - T(x + he)}{h}\nabla p(T(x), x) + \frac{1}{h}\int_{T(x + he)}^{T(x)} \nabla p(t, x + he) - \nabla p(T(x), x)\,dt \]
Using the Lipschitz continuity of $T$ from Proposition \ref{cor:regular_points_t_lipschitz} and the spacetime continuity of $\nabla p$, the second term vanishes as $h\to 0$. As we have already seen, $\nabla p$ cannot vanish on the interface due to the Hopf lemma, so in the limit, we get
\[ (\nu(x)\cdot e)\nu(x) = -\partial_e T(x)\nabla p(T(x), x) \]
Here, $\partial_e T(x)$ refers to the one-sided derivative of $T$ at $x$ in direction $e$. Since $\nabla p(T(x), x)$ has the same direction as $-\nu(x)$, we get that $\partial_e T(x) = \frac{\nu(x)\cdot e}{|\nabla p(T(x), x)|}$.

Then, it is an elementary result that a continuous function on $\R$ with continuous left derivative is differentiable. Applying it here, we get that $T$ has all two-sided directional derivatives, and we can read from the formula that we must have
\[ \nabla T(x) = \frac{\nu(x)}{|\nabla p(T(x), x)|} = -\frac{\nabla p(T(x), x)}{|\nabla p(T(x),x)|^2} \]
\end{proof}

\begin{prop}
    If for some interval $(t_0, t_1)$, we have that $\Sigma_t$ is empty for every $t\in (t_0, t_1)$, then $T\in C^{1,1/2 - \varepsilon}_{\loc}(\Omega_{t_1}\setminus \overline{\Omega_{t_0}})$ for every $\varepsilon\in (0, \frac{1}{2})$.
\end{prop}
\begin{proof}
    From the previous lemma, we need to show spacetime continuity of $\nabla p$ to establish differentiability of $T$. We do this using the $C^{1,\alpha}$ global Schauder estimates, which applied to a function $u$ on a $C^{1,\alpha}$ domain $\Omega$, give that
    \[ \|u\|_{C^{1,\alpha}(\overline{\Omega})} \leq C(\|\Delta u\|_{L^\infty(\Omega)} + \|u\|_{C^{1,\alpha}(\partial \Omega)}) \]
    for some $C$ which depends only on $\alpha$ and $\Omega$. In our case, $C$ will actually be locally uniform in $t$ for $\Omega_t$, since the global Schauder estimates are proved by patching interior and boundary estimates, and we can cover the boundary with finitely many balls in which it evolves as a uniformly $C^{1,1}$ graph for some time interval.

    Then, specifically, we will apply the Schauder estimate to $p(t) - p(s)$ on $\Omega_s$, for $t > s$, to bound $\|\nabla p(t) - \nabla p(s)\|_{L^\infty(\Omega_s)}$ in terms of $|t - s|$. Since $n\in C^{0,1-}_t L^\infty_x$ by Lemma \ref{lem:n_regularity}, the work will be in controlling $p(t)$ on $\partial \Omega_s$. Intuitively, since $p(t)$ and its tangential derivative vanish on $\partial\Omega_t$, we expect that if the free boundary has not rotated too much between times $s$ and $t$, then these should be close to 0 on $\partial \Omega_s$. We make this quantitative using Proposition \ref{prop:unit_normal_holder}.

    First, since we have a locally uniform in time bound on $\nabla p$, we get by the radial supersolution that $D(\partial \Omega_s, \partial \Omega_t) \leq C|t - s|$ for some locally uniform in time constant, where $D$ denotes Hausdorff distance as in Definition \ref{def:reifenberg_vanishing}. Then, since $p(t)$ vanishes on $\partial \Omega_t$, we integrate along shortest-distance paths and use the gradient bound again to conclude that $\|p(t) - p(s)\|_{L^\infty(\partial \Omega_s)} = \|p(t)\|_{L^\infty(\partial \Omega_s)} \leq C|t-s|$ for some uniform $C$.

    Next, we bound the tangential part of $\nabla p(t)$ on $\partial \Omega_{s}$. Recalling our definition of $\nu(x)$ as the outward unit normal to $\Omega_{T(x)}$ at $x$, we denote the projection onto the tangential part as $P^\perp_{\nu(x)}$. Then for $x\in \partial \Omega_{s}$, we let $\tilde{x}\in \partial \Omega_{t}$ be the distance minimizer so that $|x - \tilde{x}| \leq C|t-s|$, and we have
    \begin{align*}
        |P^\perp_{\nu(x)}\nabla p(x)| &\leq |(P^\perp_{\nu(\tilde{x})} - P^\perp_{\nu(x)})\nabla p(\tilde{x})| + |P^\perp_{\nu(x)}(\nabla p(\tilde{x}) - \nabla p(x))|
    \end{align*}
    where all pressures are at time $t$, and we use that the tangential derivative of $p(t)$ on $\partial \Omega_{t}$ vanishes. The first term is controlled by the continuity of $\nu$ and our uniform bound on the pressure gradient, so by Proposition \ref{prop:unit_normal_holder}, it contributes $C|t - s|^{1/2}$. The second term is controlled by the $C^{1,\alpha}$ regularity of $p$, so it contributes $C|t - s|^\alpha$. Thus, choosing $\alpha > \frac{1}{2}$, we have $\|p(t) - p(s)\|_{C^1(\partial \Omega_{s})} \leq C|t - s|^{1/2}$.

    Finally, we improve this to H\"older by interpolation. Specifically, since $p(\tau)$ is $C^{1,\alpha}$ on $\Omega_{\tau}$, uniformly in $\tau$, for any $\alpha\in (0,1)$, we have
    \[ \frac{|\nabla p(t,x) - \nabla p(t,y)|}{|x-y|^\alpha} \leq C(\alpha) \]
    We want the H\"older seminorm of $\nabla p(t)$ on $\partial \Omega_s$ to be small, so at small scales, we rearrange this to
    \[ \frac{|\nabla p(t,x) - \nabla p(t,y)|}{|x-y|^\beta} \leq C(\alpha)|x-y|^{\alpha - \beta} \]
    for $\beta < \alpha$ to be chosen. For large scales, we use the $L^\infty$ bound for the tangential part of $\nabla p(t)$ on $\partial \Omega_s$, which gives $\frac{C|t - s|^{1/2}}{|x-y|^\beta}$ on the right hand side. Optimizing, the critical scale is $|x-y|\sim |t - s|^{1/2\alpha}$, and the $C^{1,\beta}$ seminorm will scale as $C(\alpha)|t-s|^{\frac{\alpha - \beta}{2\alpha}}$. In particular, this shows that by choosing $\alpha$ close to 1 and $\beta$ close to 0, we can get arbitrarily close to $\frac{1}{2}$, so for each $\varepsilon > 0$, we have some $\beta > 0$ and some $C$ such that
\[ \|p(t) - p(s)\|_{C^{1,\beta}(\partial \Omega_s)} < C|t - s|^{1/2 - \varepsilon} \]
    Then, combining this with the regularity of the nutrient from Lemma \ref{lem:n_regularity}, we get $\|p(t) - p(s)\|_{C^{1,\beta}(\overline{\Omega}_{s})} \leq C|t - s|^{1/2 - \varepsilon}$ from the boundary Schauder estimate. We conclude that for points $(x,t)$ and $(y,s)$ in the region with $t > s$, we have
    \begin{equation}\label{eq:grad_p_spacetime_cts}
       |\nabla p(x,t) - \nabla p(y,s)| \leq |\nabla p(x,t) - \nabla p(y,t)| + |\nabla p(y,t) - \nabla p(y,s)| \leq C|x-y|^{\alpha} + C|t-s|^{1/2 - \varepsilon} 
    \end{equation}
    for any fixed $\alpha\in (0,1)$ and $\varepsilon > 0$, which proves the spacetime continuity of $\nabla p$.

    Then, by the lemma, $\nabla T(x) = -\frac{\nabla p(x, T(x))}{|\nabla p(x, T(x))|^2}$, where we have $|\nabla p(x, T(x))|$ locally uniformly bounded away from 0 by the Hopf lemma. In particular, it follows from (\ref{eq:grad_p_spacetime_cts}) and the Lipschitz continuity of $T$ that
    \[ |\nabla p(x, T(x)) - \nabla p(y, T(y))| \leq C|x - y|^{1/2 - \varepsilon} \]
    so we obtain that $\nabla T\in C^{0,1/2 - \varepsilon}_{\loc}$ directly from its formula.
\end{proof}

This implies that the free boundary has $C^{2, 1/2 - \varepsilon}$ regularity at the relevant times. Subsequently, the regularity of $\nu$ in Proposition \ref{prop:unit_normal_holder} can be improved using second order approximations to the free boundary, leading to a minor improvement in the H\"{o}lder exponent.

\section{Appendix: Obstacle problem with $C^{0,\alpha}$ source}

In this section, we collect several known facts about obstacle problems with H\"{o}lder continuous data. Many results for the model obstacle problem with constant source carry over to the H\"{o}lder continuous case with minor modifications, as noted in \cite{caffarelli98} and \cite{weiss}. We will cite several results from \cite{blank}, \cite{monneau}, and \cite{colombo}, which offer careful treatments of this topic.

Let us consider solutions $u$ to the obstacle problem
\begin{equation}\label{eq:obst}
    \Delta u = (1 + f)\chi_{\{u > 0\}}
\end{equation}
on $B_1$, where $f$ is known a priori to vanish on the free boundary $\Gamma(u) = \partial \{ u > 0\}$. When $f$ is sufficiently regular, this equation has similar local behavior to the model case where $f\equiv 0$; in particular, we make the assumption $f\in C^{0,\alpha}(B_1)$ for some $\alpha\in (0,1)$. Regarding notation, we write $\Omega(u) = \{ u > 0 \}$, $\Lambda(u) = \{ u = 0\}$, and in many cases we will refer to a tuple $(u,f)$ as the solution to (\ref{eq:obst}). We take $\lambda = \inf_{B_1} 1 + f, \mu = \sup_{B_1} 1 + f$, and we will have the standing assumption that $\lambda > \frac{1}{2}$, which holds if $0\in \Gamma(u)$ and $[f]_{C^{0,\alpha}(B_1)}$ is sufficiently small.

On the free boundary, we have $u = 0, \nabla u = 0$, so we should expect $u$ to have quadratic growth away from the free boundary in the positive set. Of course, $u$ is not regular enough to admit a second order Taylor expansion due to the jump in the second derivatives along the free boundary, but nevertheless, we recover several results to the same effect:

\begin{lemma}[Quadratic nondegeneracy, \cite{blank} Thm. 2.1]\label{lem:obst_nondegen}
If $0\in \overline{\Omega(u)}$, then for all $r < 1$,
    \[ \sup_{B_r} u \geq \frac{\lambda}{2d}r^2 \]
\end{lemma}
\begin{lemma}[Quadratic bound, \cite{blank} Thm. 2.4]\label{lem:obst_quadratic_bound}
If $0\in \Gamma(u)$, then for all $r < \frac{1}{2}$,
    \[ \sup_{B_r} u \leq C(d)\mu r^2 \]
\end{lemma}
\begin{lemma}[Regularity up to the free boundary, \cite{blank} Thm. 2.3]\label{lem:obst_regularity}
    If $0\in \Gamma(u)$, then $\|u\|_{C^{1,\beta}(B_{1})} \leq C(d, \beta) \mu$ for all $\beta\in (0, 1)$.
\end{lemma}

We adopt the following notation for the quadratic rescalings:
\begin{equation}\label{eq:quadratic_rescaling}
    u_r(x) = r^{-2}u(rx)
\end{equation}
\begin{equation}
    u_0(x) = \lim_{r\to 0^+} u_r(x), \hbox{ provided this limit exists }
\end{equation}
The combined results of lemmas \ref{lem:obst_nondegen}, \ref{lem:obst_quadratic_bound}, and \ref{lem:obst_regularity} can be used to derive Lemma \ref{lem:quad_blowup_cpt}: the 
$C^{1,\beta}$ compactness of the quadratic blowup sequence $(u_r)$. As we discuss in the main paper with Lemma \ref{lem:caffarelli_dichotomy}, this compactness improves to convergence of the blowup sequence when $f\in C^{0,\alpha}$, and we classify points as regular or singular based on the blowup profile. In general, regular points can be identified at finite scales, using criteria such as Lemma \ref{lem:modified_regular_pt_criterion}.

The quadratic blowup proves to be a key tool in understanding local behavior of the free boundary. For regular points, we use comparison and stability results to show flatness of the free boundary, which leads to regularity of the boundary. For singular points, we use monotonicity formulas and compactness results to show that they can be locally contained in $C^1$ manifolds. Since the treatment of these cases diverges considerably, we will split them into the next two sections.

\subsection{Regular points}

In this section, we review several results from \cite{blank} which connect the regularity of the free boundary at regular points and the scale at which this regularity is achieved with the regularity of the source term and the scale at which the zero set becomes large. The regularity of the free boundary can be summarized as follows:
\begin{lemma}[\cite{blank} Thm. 7.2]
    If $f\in C^{0,\alpha}$ with $\alpha\in (0, 1]$, then in a neighborhood of a regular point, the free boundary is a $C^{1,\alpha}$ graph.
\end{lemma}

In light of this result, we allow the case $\alpha = 1$ for the rest of this subsection.

For the main paper, we require a quantified version of this lemma, in order to apply it uniformly to the family of obstacle problems $w(\cdot, t)$. Thus, we will retrace Blank's approach in this section while keeping track of its dependencies.

\begin{definition}\label{def:reifenberg_vanishing}
    Let $S\subset \R^d$ be a compact set. We define the modulus of flatness,
    \[ \theta(r) = \sup_{0 < \rho \leq r} \sup_{x\in S} \inf_L \frac{D(L\cap B_\rho(x), S\cap B_\rho(x))}{\rho} \]
    where the inner infimum is over all hyperplanes $L$ containing $x$, and $D$ denotes Hausdorff distance:
    \[ D(A,B) = \max(\sup_{x\in A} \dist(x, B), \sup_{y\in B} \dist(y, A)) \]
    We say that $S$ is $\delta$-Reifenberg flat if there exists $R$ such that $\theta(r) \leq 2\delta$ for all $r < R$, and Reifenberg vanishing if $\theta(r)\to 0$.
\end{definition}

\begin{lemma}[\cite{blank} Theorem 6.7]\label{lem:flat_implies_smooth}
    Let $S$ be a compact Reifenberg vanishing set with modulus of flatness $\theta$ satisfying $\int_0^1 \frac{\theta(r)}{r}\,dr < \infty$. Then there exist constants $C_0, C_1$ such that if $\int_0^\rho \frac{\theta(r)}{r}\,dr < C_0$, then there exists a coordinate system in which $S\cap B_{\rho/2}$ is the graph of a $C^1$ function $g$, such that $\nabla g$ is continuous with modulus of continuity $C_1 \int_0^r \frac{\theta(s)}{s}\,ds$.
\end{lemma}

\begin{lemma}[\cite{blank} Theorem 7.1]
    Suppose that $u$ solves the obstacle problem $\Delta u = f\chi_{\{u > 0\}}$ in $B_1$ with $\lambda \leq f \leq \mu$ and $f$ Dini continuous with modulus $\sigma$. Then if $0$ is a regular point and the free boundary is $\delta$-Reifenberg flat in $B_{3/4}$ for some sufficiently small $\delta$, then the modulus of flatness of the free boundary inside $B_{1/2}$ is controlled by $C\sigma(r)$.
\end{lemma}
In particular, these two results imply the $C^{1,\alpha}$ regularity of the free boundary near regular points when $f\in C^{0,\alpha}$ with $f(0) = 1$, at a scale depending on $[f]_{C^{0,\alpha}}(B_1)$ and the scale at which the $\delta$ Reifenberg flatness is achieved. For this, we have another result from Blank
\begin{lemma}[\cite{blank} Theorem 6.4]
    Let $\varepsilon\in (0, \frac{1}{4})$, and suppose we have $f$ with $\lambda \leq f \leq \mu$ and $u$ a solution to the obstacle problem $\Delta u = f\chi_{\{u > 0\}}$ in $B_1$. If $\mu - \lambda$ is sufficiently small, then there exist constants $r_0, \tau, \delta\in (0,1)$ depending on $d, \mu, \lambda, \varepsilon$ for which the following holds:
    
    If for some $t \leq r_0$, we have
    \[ \frac{|B_t\cap \{u = 0\}|}{|B_t|} > \varepsilon, \]
    then $\overline{B_{\tau t}}\cap \partial\{u > 0\}$ is $\delta$-Reifenberg flat.

    Moreover, as $\mu - \lambda \to 0$, $\delta\to 0$. In particular, if $f$ is continuous, $\delta$ can be taken to be arbitrarily small (with all parameters now also depending on the modulus of continuity of $f$).
    
\end{lemma}

To be more precise, suppose $u$ solves the obstacle problem $\Delta u = f\chi_{\{u > 0\}}$ in $B_1$ with $f$ taking values in $[\lambda, \mu]$, and write $u_c$ for the solution to the obstacle problem $\Delta u_c = c\chi_{\{u_c > 0\}}$ such that $u_c|_{\partial B_1} = u$. Then $\{ u_\lambda = 0\} \subset \{ u = 0\} \subset \{u_\mu = 0 \}$. Moreover, if $0$ is a regular point for $u$, then there exists $c\in [\lambda, \mu]$ such that $0$ is a regular point for $u_c$. Then Blank's uniform stability theorem for regular points (\cite{blank} Theorem 5.4) gives that in $B_{1/2}$, there is a universal $C$ such that for any $c'$, we have $d(FB(u_c), FB(u_{c'})) \leq C|c - c'|$ where $FB(v)$ denotes the free boundary of $v$. Then we get flatness of $FB(u)$ by trapping it between $FB(u_\lambda)$ and $FB(u_\mu)$ and using the stability and $C^{1,\alpha}$ regularity of the constant-source free boundaries. As we zoom in, the $C^{1,\alpha}$ seminorm goes to 0, and if $f$ is continuous, $|\mu - \lambda|\to 0$, so we can get $\delta$-Reifenberg flatness with arbitrarily small $\delta$. In particular, the $C^{1,\alpha}$ seminorm is uniformly bounded, depending only on the scale at which the density of the zero set is sufficiently large, and the rate at which $|\mu - \lambda|\to 0$ depends only on the modulus of continuity of $f$. Thus, we can replace the hypothesis of $\delta$-Reifenberg flatness in Lemma \ref{lem:flat_implies_smooth}, to conclude:
\begin{lemma}\label{lem:regular_point_boundary_estimate}
    Suppose that $(u, f)$ solve (\ref{eq:obst}) with $f\in C^{0,\alpha}$ for $\alpha\in (0, 1]$. Then there exist $r_0, \varepsilon_0$ such that if $0$ is a free boundary point and $\frac{|B_t\cap \{u = 0\}|}{|B_t|} > \varepsilon_0$ for some $t \leq r_0$, then there exists $r = r(t,[f]_{C^{0,\alpha}})$ such that the free boundary is a $C^{1,\alpha}$ graph in $B_r$, with $C^{1,\alpha}$ seminorm controlled by $[f]_{C^{0,\alpha}(B_1)}$ and $r$.
\end{lemma}

\subsection{Singular points}

In this section, we use monotonicity formulas to study the continuity of the blowup limit at singular points and the rate of convergence for the blowup sequence. First, we have the following result:

\begin{lemma}[\cite{colombo}, Thm 5]
    Restricted to the singular set, $D^2u$ is continuous with a logarithmic modulus of continuity. In particular, in a neighborhood of a singular point, the singular set is contained in a $C^{1,\log}$ manifold of dimension $\dim \ker D^2 u$.
\end{lemma}

The result of \cite{colombo} is obtained using an epiperimetric inequality to control the Weiss monotonocity formula, introduced in \cite{weiss}. For simplicity, we will instead consider the related Monneau monotonicity formula, at the cost of the explicit logarithmic modulus of continuity. The following result is part of the proof of \cite{monneau} Theorem 1.9.

\begin{lemma}\label{lem:monneau_mf}
    Define
\begin{equation}\label{eq:monneau_mf}
    \Xi_u^q(r) = r^{-(d+3)}\int_{\partial B_r} (u-q)^2 = \int_{\partial B_1} (u_r - q)^2
\end{equation}
    where $u$ solves (\ref{eq:obst}) with $0$ as a singular point, and $q$ is a quadratic form $q(x) = \frac{1}{2}x\cdot Qx$ with $Q\geq 0, \tr Q = 1$. Then
    \[ \frac{d}{dr}\Xi_u^q(r)\geq -Cr^{\alpha - 1} \]
    where $C = C([f]_{C^{0,\alpha}})$. In particular, the limit $\Xi_u^q(0+)$ exists.
\end{lemma}
As a corollary, we can show that near a singular point, $u$ approximates a global solution at a uniform scale. This extends Lemma 13 of \cite{caffarelli98} to the $C^{0,\alpha}$ source case.

We break the proof into two steps. First, we will prove the following slightly weaker claim:

\begin{lemma}
    Let $u$ solve \eqref{eq:obst} with 0 as a singular point. Then for every $\varepsilon > 0$, there exists $\delta = \delta(\varepsilon, [f]_{C^\alpha})$ such that $\|u_{\delta} - q\|_{C^1(B_1)} < \varepsilon$, for some $q$ of the form $q(x) = \frac{1}{2}x\cdot Qx$ with $Q$ a positive semidefinite matrix of trace 1. Here, we use the notation defined in \eqref{eq:quadratic_rescaling}.
\end{lemma}
\begin{proof}
Suppose for contradiction that we can find a sequence $v^k$ solving $\Delta v^k = (1 + f^k)\chi_{\{v^k > 0\}}$ on $B_1$ such that the result fails along the sequence $v^k_{1/k}$. That is, for each $k$, and every $q$ of the form in the statement of the lemma,
\begin{equation}\label{eq:sing_approximation_failure}
    \|v^k_{1/k} - q\|_{C^1(B_1)} > \varepsilon
\end{equation}
Then the sequence $v^k_{1/k}$, defined on the expanding balls $B_k$, converges along a subsequence on all compact sets in $C^1$, to some $v^\infty$. Since the $f^k$ are uniformly $C^\alpha$, the sequence $f^k(\frac{x}{k})$ converges locally uniformly to 0. It follows that $v^\infty$ is a nonnegative solution to $\Delta v^\infty = \chi_{\{v^\infty > 0\}}$ on $\R^d$ with $v^\infty(0) = 0$.

Next, we show that the zero set of $v^\infty$ has empty interior. Suppose otherwise, and we have some $B_r(x)\subset B_1$ such that $v^\infty\equiv 0$ on $B_r(x)$. This implies that $v^k_{1/k}$ is $o(1)$ on $\partial B_r(x)$ as $k\to \infty$, and an application of the nondegeneracy bound Lemma \ref{lem:obst_nondegen} yields that $v^k_{1/k}\equiv 0$ on $B_{r/2}(x)$ for $k$ sufficiently large.

Then, since the zero set of $v^\infty$ has empty interior, all points in the zero set are singular free boundary points, and $v^\infty$ is twice differentiable at those points. In particular, we get that $v^\infty$ solves $\Delta v^\infty = 1$. As noted in the remarks after Lemma 13 of \cite{caffarelli98}, this, along with the quadratic growth estimate, Lemma \ref{lem:obst_quadratic_bound}, implies that $v^\infty$ is a quadratic polynomial of the type in the statement above. Thus, taking $q = v^\infty$, we get a contradiction to \eqref{eq:sing_approximation_failure} for $k$ sufficiently large, which completes the proof.
\end{proof}

\begin{lemma}\label{lem:uniform_sing_blowup_scale}
    For every $\varepsilon > 0$, there exists $\delta = \delta(\varepsilon, [f]_{C^{0,\alpha}(B_1)})$ such that $\|u_\delta - u_0\|_{C^1(B_1)} < \varepsilon$, where we use the notation of \eqref{eq:quadratic_rescaling}. Equivalently, we have $\|u - u_0\|_{C^1(B_\delta)} < o(\delta^2)$ as $\delta\to 0$, where the little $o$ depends only on $[f]_{C^{0,\alpha}}(B_1)$.
\end{lemma}
\begin{proof}
Let $q$ be the quadratic form given by the previous lemma. We apply Monneau's monotonicity formula, $\Xi_u^q$, defined in \eqref{eq:monneau_mf}. The derivative bound from Lemma \ref{lem:monneau_mf} gives
\[  \Xi_u^q(0) - \Xi_u^q(\delta) \leq C\delta^\alpha \]
We have $\|u_\delta - q\|_{C^1(B_1)} < \varepsilon$ from the lemma, so it follows that $\|u_\delta - q\|_{L^2(\partial B_1)} < C\varepsilon$ for some dimensional constant. Then $\|u_0 - q\|_{L^2(\partial B_1)}^2 \leq C\varepsilon^2 + C\delta^\alpha$. But now we recall that $u_0, q$ are both quadratic forms, and so by equivalence of norms on $\R^{d\times d}$, there exists a dimensional constant for which $\|u_0 - q\|_{C^1(B_1)} \leq C\|u_0 - q\|_{L^2(\partial B_1)}$. Combining our estimates for $\|u_0 - q\|_{C^1(B_1)}$ and $\|u_{\delta} - q\|_{C^1(B_1)}$, we conclude.
\end{proof}

\bibliography{ref_CJK}
\bibliographystyle{alpha}
\end{document}